\documentclass[sn-mathphys-num]{sn-jnl}

\usepackage{graphicx}%
\usepackage{multirow}%
\usepackage{amsmath,amssymb,amsfonts}%
\usepackage{amsthm}%
\usepackage{mathrsfs}%
\usepackage[title]{appendix}%
\usepackage{xcolor}%
\usepackage{textcomp}%
\usepackage{manyfoot}%
\usepackage{booktabs}%
\usepackage{algorithm}%
\usepackage{algorithmicx}%
\usepackage{algpseudocode}%
\usepackage{listings}%
\usepackage{rotating}

\usepackage{setspace}
\usepackage{makecell}
\usepackage{cleveref}
\usepackage{subcaption}

\RequirePackage{multirow}
\RequirePackage{adjustbox}
\RequirePackage{colortbl}
\RequirePackage{import}
\RequirePackage{hyperref}
\RequirePackage{stackrel}
\RequirePackage{wrapfig}
\RequirePackage{xspace} 
\RequirePackage{bm}

\newcommand{\Span}[1]{\text{Span}\left\{#1\right\}}

\DeclareMathOperator*{\argmin}{arg\,min}                   
\DeclareMathOperator*{\argmax}{arg\,max}                   
\renewcommand{\t} {^{\top}}                                
\newcommand{\norm} [2][]{\left\|#2\right\|_{#1}}           
\newcommand{\diag} [1]  {{\rm diag\!}\left( #1 \right)}    
\newcommand{\kron}{\otimes}                                

\newcommand{\bfGamma}{{\boldsymbol{\Gamma}}}

\newcommand{\bfTheta}{{\boldsymbol{\Theta}}}

\newcommand{\bfSigma}{{\boldsymbol{\Sigma}}}

\newcommand{\bfOmega}{{\boldsymbol{\Omega}}}

\newcommand{\bfepsilon}{{\boldsymbol{\epsilon}}}

\newcommand{\bfmu}{{\boldsymbol{\mu}}}

\newcommand{\bfxi}{{\boldsymbol{\xi}}}

\newcommand{\bfA}{{\bf A}}
\newcommand{\bfB}{{\bf B}}

\newcommand{\bfD}{{\bf D}}

\newcommand{\bfF}{{\bf F}}
\newcommand{\bfG}{{\bf G}}
\newcommand{\bfH}{{\bf H}}
\newcommand{\bfI}{{\bf I}}

\newcommand{\bfL}{{\bf L}}
\newcommand{\bfM}{{\bf M}}

\newcommand{\bfP}{{\bf P}}
\newcommand{\bfQ}{{\bf Q}}
\newcommand{\bfR}{{\bf R}}
\newcommand{\bfS}{{\bf S}}

\newcommand{\bfU}{{\bf U}}
\newcommand{\bfV}{{\bf V}}
\newcommand{\bfW}{{\bf W}}

\newcommand{\bfY}{{\bf Y}}

\newcommand{\bfb}{{\bf b}}

\newcommand{\bfe}{{\bf e}}

\newcommand{\bfg}{{\bf g}}

\newcommand{\bfu}{{\bf u}}
\newcommand{\bfv}{{\bf v}}

\newcommand{\bfx}{{\bf x}}
\newcommand{\bfy}{{\bf y}}
\newcommand{\bfz}{{\bf z}}
\newcommand{\bfzero}{{\bf0}}

\newcommand{\calK}{\mathcal{K}}

\newcommand{\calN}{\mathcal{N}}
\newcommand{\calO}{\mathcal{O}}
\newcommand{\calP}{\mathcal{P}}

\newcommand{\calR}{\mathcal{R}}

\newcommand{\bbR}{\mathbb{R}}

\newcommand{\jmc}[1]{\textcolor{black}{#1}}

\theoremstyle{thmstyleone}%
\theoremstyle{thmstyletwo}%

\theoremstyle{thmstylethree}%

\usepackage{color}
\newcommand{\rev}[1]{\textcolor{black}{#1}}
\usepackage{makecell}
\usepackage{booktabs}

\begin{document}

\title[Efficient sampling approaches based on generalized Golub-Kahan methods for large-scale hierarchical Bayesian inverse problems]{Efficient sampling approaches based on generalized Golub-Kahan methods for large-scale hierarchical Bayesian inverse problems}


\author[1]{\fnm{Elle} \sur{Buser}}\email{ebuser@emory.edu}
\author*[1]{\fnm{Julianne} \sur{Chung}}\email{jmchung@emory.edu}

\affil[1]{\orgdiv{Department of Mathematics}, \orgname{Emory University}, \orgaddress{\street{400 Dowman Drive}, \city{Atlanta}, \postcode{30322}, \state{Georgia}, \country{USA}}}


\abstract{Uncertainty quantification for large-scale inverse problems remains a challenging task. For linear inverse problems with additive Gaussian noise and Gaussian priors, the posterior is Gaussian but sampling can be challenging, especially for problems with a very large number of unknown parameters (e.g., dynamic inverse problems) and for problems where computation of the square root and inverse of the prior covariance matrix are not feasible.  Moreover, for hierarchical problems where several hyperparameters that define the prior and the noise model must be estimated from the data, the posterior distribution may no longer be Gaussian, even if the forward operator is linear.  Performing large-scale uncertainty quantification for these hierarchical settings requires new computational techniques. In this work, we consider a hierarchical Bayesian framework where both the noise and prior variance are modeled as hyperparameters.  Our approach uses Metropolis-Hastings independence sampling within Gibbs where the proposal distribution is based on generalized Golub-Kahan methods. We consider two proposal samplers, one that uses a low-rank approximation to the conditional covariance matrix and another that uses a preconditioned Lanczos method. Numerical examples from seismic imaging, dynamic photoacoustic tomography, and atmospheric inverse modeling demonstrate the effectiveness of the described approaches. }

\keywords{hierarchical Bayes, Gibbs sampler, inverse problems, uncertainty quantification, Krylov methods}


\pacs[MSC Classification]{65F22, 65M32, 62F10}

\maketitle
\section{Introduction}
\label{sec:intro}
Inverse problems arise in many scientific applications, where the main goal is to use collected measurements or observations to estimate some underlying unknown parameters of physical models. We focus on inverse problems in imaging, where the unknown parameters represent detailed spatial or spatiotemporal reconstructions of physical properties such as images of attenuation coefficients in X-ray tomography or spatiotemporal maps of greenhouse gas emission fluxes in atmospheric inverse modeling \cite{chung2024computational}. In dynamic inverse problems, the unknown parameters can change over time, easily resulting in millions of unknown parameters.  For example, in dynamic atmospheric inverse modeling, the goal is to estimate thousands of emission fluxes at 3-hourly time intervals over multiple months or years.  Obtaining such estimates is computationally challenging, and recent works in the field of computational inverse problems have addressed various theoretical and computational advancements (e.g., developing improved reconstruction algorithms that enable faster reconstructions at higher resolutions with higher accuracy), e.g.,~\cite{cho2022computationally,pasha2023computational,chung2018efficient}.  
Many of these approaches rely on sophisticated tools from optimization and numerical linear algebra for obtaining reconstructions.  However, to provide quantification of uncertainty about the solutions of inverse problems, we follow a Bayesian interpretation of inverse problems.

In a Bayesian formulation, the parameters of interest and the observed data are modeled as random variables, and any prior knowledge or lack thereof (e.g., uncertainty in the parameters) is encoded in the prior distribution and any noise or measurement error is encoded in the likelihood function (along with the forward process).  Contrary to deterministic approaches, where a single solution is provided, the Bayesian approach provides a distribution of plausible solutions in the form of samples from \jmc{the} posterior probability distribution.  Given the observation data, Bayes' law allows uncertainty quantification (UQ) via incorporation of prior knowledge about the unknown parameters (in the form of the prior) and the likelihood. Good references on Bayesian or statistical approaches to inverse problems and computational UQ include \cite{calvetti2023bayesianbook,bardsley2018computational,calvetti2007introduction}.  

However, there are various challenges that have hindered the extension of many of these approaches to the large-scale problems of interest. For hierarchical Bayesian approaches where the prior and/or likelihood distributions depend on additional (hyper-)parameters, hyperpriors must be incorporated.  This usually results in complicated posterior distributions that do not have a closed form, thereby requiring expensive approximation techniques \cite{ghattasinfinitebayes,BuiThanh2014,flath2011fast}.  Moreover, even when it is possible to derive a closed-form for the posterior distribution, drawing samples from the posterior distribution can \jmc{be} computationally expensive.

\paragraph{Overview of contributions}
For sampling the posterior distribution in hierarchical Bayesian \jmc{inverse} problems, we focus on Markov chain Monte Carlo (MCMC) methods, particularly Gibbs sampling and its variants.  For large-scale inverse problems, the main computational bottleneck of these MCMC routines is the repeated sampling from high-dimensional Gaussian random variables, which requires a symmetric factorization of a large, and often dense, covariance matrix that is changing at each iteration.  We seek to reduce the computational burden of repeated sampling by using Metropolis-Hastings within Gibbs, with proposal samplers based on generalized Golub-Kahan methods. Similar to the approach described in \cite{brown2018low}, one approach we consider is to use a proposal distribution based on a low-rank approximation of the prior-preconditioned Hessian. We exploit generalized Golub-Kahan approximations for independence sampling, where the added benefits are that more general prior covariance matrices can be included (since we only require matrix-vector \jmc{(mat-vec)} multiplications with the prior covariance matrix) and we can reuse Krylov matrices across MCMC iterations (since they are independent of the hyperparameters).
For these independence samplers, we derive explicit formulas for the acceptance rates and demonstrate the computational benefits of their use for hierarchical Bayesian inverse problems in a wide range of applications \jmc{and} for problems with thousands of unknown parameters.

\paragraph{Related work}
For hierarchical Bayesian inverse problems, efficient optimization techniques have been considered for \jmc{maximum a posteriori (MAP)} estimation, where many previous works use an iterative alternating scheme or Krylov-based iterative methods \cite{calvetti2015hierarchical,calvetti2019hierachical,lindbloom2025efficient,lindbloom2025priorconditioned}. A unified theoretical framework and reconstruction error bounds for MAP estimates are provided in a recent survey paper \cite{sanz2025hierarchical}. 
\jmc{For UQ, sampling techniques have been considered for hierarchical Bayesian inverse problems, but they can be challenging, especially for problems with many unknown parameters. }
The most common approach is to use centered algorithms, such as block Gibbs and \jmc{Metropolis-Hastings}-within-Gibbs algorithms, but a potential concern is poor mixing in the presence of strong correlations \cite{agapiou2014analysis,ascolani2024scalability}.  Thus, there have been various approaches that exploit prior normalization for accelerating MCMC for problems with heavy-tailed priors \cite{glaubitz2025efficient,calvetti2024computationally,calvetti2025subspace,agrawal2022variational}.  
Although these works highlight the recent and increased interest in sampling for hierarchical Bayesian inverse problems, many of these approaches are demonstrated for small problems with moderately sized unknowns (e.g., on the order of hundreds) or consider special cases, e.g., where the size of the data space is significantly smaller than the number of unknown parameters \cite{calvetti2025subspace}.  
\jmc{Other existing works} focus on sparsity-promoting priors that assume conditionally Gaussian priors with variances that are mutually independent and distributed according to a generalized gamma hyperprior \rev{\cite{Calvetti2020sparsehyperpriors}}.

Our proposed approach is more closely related to previous approaches that have considered low-rank approximations obtained via randomization within Gibbs sampling \cite{brown2018low, saibaba2021randomized}.  However, since such methods can only provide good approximations for severely low-rank matrices, they are not suitable for large-scale \jmc{linear} inverse problems of interest where \jmc{the forward model matrix} has a slowly decaying spectrum. There is also related work on marginal then conditional samplers that work for image deblurring applications \cite{fox2016fast}.
Our approach also builds on previous work on sampling from high-dimensional Gaussians that was considered in \cite{UQlargebayesian}.  Since these methods were considered for fixed hyperparameters, a simple extension to hierarchical problems would be expensive, since each MCMC iteration would require a large linear solve.  By exploiting the shift invariance property of the generalized Krylov subspaces, we provide a computationally efficient approach for hierarchical sampling where the subspace vectors can be reused across MCMC iterations.  This is a feature that is not possible in other settings (e.g., in the sparsity-promoting setting where the hyperparameters define the diagonals of the covariance matrix). Also, by allowing priors defined using covariance kernel functions, we allow hierarchical sampling for more general smoothness priors.
We mention that although most approaches for sampling from a high-dimensional posterior are based on MCMC techniques \cite{ghanem2017handbook,dashti2015bayesian}, there are related works that exploit machine learning techniques, e.g., using deep posterior sampling \cite{adler2019deep}, or hybrid MCMC algorithms that use emulators and autoencoders for Bayesian UQ \cite{lan2022scaling}. However, these approaches require supervised training data that may not be readily available.

\paragraph{Outline}
The paper is organized as follows. In \Cref{sec:BayesianIP}, we describe a hierarchical Bayesian framework for a general linear inverse problem.  We formulate the posterior density function and review various MCMC samplers for sampling the posterior.  Then in \Cref{sec:methods}, we provide an overview of generalized Golub-Kahan methods and describe two approaches for their use in independence sampling in Metropolis-Hastings within Gibbs. Numerical results for various large-scale image processing applications are provided in \Cref{sec:numerics}. Conclusions and future work are described in \Cref{sec:conclusions}.

\section{Hierarchical Bayesian inverse problem}
\label{sec:BayesianIP}
Consider a linear inverse problem of the form 
\begin{equation}
     \bfb =  \bfA\bfx + \bfe
\end{equation}
where observation data $\bfb\in\bbR^{m}$ are corrupted by measurement error $\bfe\in\bbR^{m}$, $\bfA\in\bbR^{m\times n}$ represents the forward parameter-to-observable map, and $\bfx\in\bbR^{n}$ contains the unknown solution.  In a deterministic inverse problem setting, the goal of the inverse problem is to compute a reconstruction of $\bfx$ (e.g., a point estimate),  given $\bfA$ and $\bfb$.  In a Bayesian setting, the goal is to fully characterize the posterior probability distribution (thereby quantifying uncertainty about the solutions), given assumptions or prior knowledge about the unknowns.  

Given $\lambda, \delta >0$ and $\bfmu \in \bbR^n$, we assume that $\bfe$ and $\bfx$ are independent Gaussian random variables such that 
\begin{equation}
    \label{eq:assumptions_e_x}
\bfe \mid \lambda \sim\calN\left( \bfzero,\lambda^{-1}\bfR \right)\quad  \mbox{and} \quad \bfx \mid \delta \sim\calN\left( \bfmu, \delta^{-1}\bfQ \right)
\end{equation}
where $\bfQ$ and $\bfR$ are symmetric positive definite covariance matrices. With this model, $\bfb \mid \bfx, \lambda \sim \calN(\bfA \bfx,\lambda^{-1}\bfR)$ so that the likelihood is given by
\begin{equation}
\label{eq:likelihood}
\pi(\bfb \mid \bfx, \lambda) \propto \lambda^{m/2} \exp\left( -\frac{\lambda}{2} (\bfb - \bfA \bfx)\t \bfR^{-1} (\bfb - \bfA \bfx) \right)
\end{equation}
and the prior is given by
\begin{equation}
\label{eq:prior}
\pi(\bfx \mid \delta) \propto \delta^{n/2} \exp\left( - \frac{\delta}{2} (\bfx- \bfmu)\t \bfQ^{-1} (\bfx - \bfmu) \right).
\end{equation}
For fixed $\lambda$ and $\delta$, we get a conditional that is also Gaussian.  That is, $\bfx \mid \bfb, \lambda, \delta \sim \calN(\bfx_{\rm cond}, \bfGamma_{\rm cond})$, where 
\begin{equation}
    \label{eq:cond_x}
\bfGamma_{\rm cond} = (\lambda \bfA\t \bfR^{-1} \bfA + \delta \bfQ^{-1})^{-1} \mbox{ and }\bfx_{\rm cond} =  \bfGamma_{\rm cond} (\lambda \bfA\t \bfR^{-1}\bfb + \delta \bfQ^{-1}\bfmu).
\end{equation}
More specifically,
\begin{equation}
\label{eq:conditional}
\pi(\bfx \mid \bfb, \lambda, \delta) \propto  \exp\left(-\frac{\lambda}{2} (\bfb - \bfA \bfx)\t \bfR^{-1} (\bfb - \bfA \bfx) - \frac{\delta}{2} (\bfx - \bfmu)\t \bfQ^{-1} (\bfx-\bfmu) \right).
\end{equation}
The conditional posterior mode is the minimizer of the negative log-likelihood and corresponds to \jmc{the solution of} a general-form Tikhonov problem,
\begin{equation}
\label{eq:Tik}
    \widehat \bfx = \argmax_\bfx \pi(\bfx \mid \bfb, \lambda, \delta) = \argmin_\bfx \frac{\lambda}{2} \|\bfb - \bfA \bfx \|_{\bfR^{-1}}^2 + \frac{\delta}{2} \| \bfx - \bfmu \|_{\bfQ^{-1}}^2,
\end{equation}
where $\norm[\bfM]{\bfx} = \sqrt{\bfx\t\bfM\bfx}$ for \jmc{symmetric positive definite matrix} $\bfM$.

We assume that hyperparameters $\lambda$ and $\delta$ are unknown, and we assume that they are random variables distributed according to some hyperprior. For example, a common assumption is to use gamma hyperpriors defined by
\begin{equation}
\label{eq:hyperpriors}
    \pi(\lambda)\propto \lambda^{\alpha_{\lambda}-1}\exp(-\beta_{\lambda}\lambda) \quad \mbox{and} \quad
    \pi(\delta)\propto \delta^{\alpha_{\delta}-1}\exp(-\beta_{\delta}\delta),
\end{equation}
where $\alpha_\lambda, \alpha_\delta$ and $\beta_\lambda, \beta_\delta$ are given parameters defining the shape and rate of the distributions.
Using Bayes' theorem with assumptions \eqref{eq:likelihood}, \eqref{eq:prior}, and \eqref{eq:hyperpriors}, the (non-Gaussian) joint posterior probability density \jmc{function} is given by,
\begin{align}
    & \pi(\bfx, \lambda,\delta\mid\bfb) \notag \\
    &\propto \lambda^{m/2} \delta^{n/2} \pi(\bfb\mid\bfx,\lambda)\pi(\bfx\mid\delta)\pi(\delta)\pi(\lambda) \notag \\
    &\propto \frac{\lambda^{m/2+\alpha_{\lambda}-1}\delta^{n/2+\alpha_{\delta}-1}}{\left( (2\pi)^{n+m} |\bfQ| |\bfR| \right)^{1/2}}\exp\left( -\frac{\lambda}{2}\norm[\bfR^{-1}]{\bfA\bfx-\bfb}^{2} - \frac{\delta}{2}\norm[\bfQ^{-1}]{\bfx-\bfmu}^{2} -\beta_{\lambda}\lambda -\beta_{\delta}\delta\right) 
    \label{eq:post_nonGaus}
\end{align}
where $\pi(\bfb|\bfx,\lambda)$ and $\pi(\bfx\mid\delta)$ are the likelihood and prior density functions respectively and $| \cdot |$ denotes the determinant.

In a Bayesian framework, $\pi(\bfx, \lambda,\delta\mid\bfb)$ is the solution to the inverse problem.  However, since the full joint posterior \eqref{eq:post_nonGaus} is non-Gaussian, exploring the posterior is more challenging, especially for large-scale problems. There are various approaches \jmc{that can be} used to describe the posterior distribution.  One idea is to compute point estimates such as the \jmc{MAP} estimate, corresponding to the maximum of the posterior density function. 
Obtaining this point estimate requires sophisticated nonlinear optimization techniques, e.g., computing the MAP requires solving
\begin{equation}
\label{eq:MAP}
    \min_{\bfx,\lambda,\delta} \frac{\lambda}{2}\norm[\bfR^{-1}]{\bfA\bfx-\bfb}^{2} + \frac{\delta}{2}\norm[\bfQ^{-1}]{\bfx-\bfmu}^{2} + \beta_{\lambda}\lambda + \beta_{\delta}\delta - (m/2+\alpha_\lambda -1) \log \lambda - (n/2+\alpha_\delta -1) \log \delta.
\end{equation}
\jmc{For UQ,} a common approach is to approximate \eqref{eq:post_nonGaus} by a Gaussian distribution (e.g., by linearizing around the MAP estimate), but such approximations can be poor (e.g., if far from the MAP estimate) and \jmc{may} yield unsatisfactory uncertainty estimates \cite{gelman2013bayesian,ghoshbook,agrawal2022variational}.
\jmc{An alternative approach is} to use Monte Carlo methods for sampling from the \jmc{joint} posterior \eqref{eq:post_nonGaus}, and \jmc{the samples can be used to} obtain summary statistics (e.g., estimation of the posterior mean and variances).
\jmc{For example, for nonlinear inverse problems, stochastic Newton MCMC approaches could be used, where Gaussian proposals are constructed from local gradient information and local (or MAP-based) Hessian information \cite{martin2012stochastic,petra2014computational}. Such methods could be extended to sample from \eqref{eq:post_nonGaus}, where the Hessian matrix of the negative log posterior consists
of the inverse conditional covariance matrix, augmented with two rows and two columns (coming from the 2 hyperparameters), but these approaches would not exploit the directions in which the distribution is Gaussian (e.g., due to linearity of the problem in $\bfx$).}

A common type of MCMC algorithm is the standard block Gibbs approach \cite{geman1984Gibbs,bardsley2018computational}.  The main idea is to alternate sampling from the conditional distributions,
\begin{align}
    \lambda\mid\bfb,\bfx,\delta &\sim \Gamma\big(m/2 + \alpha_{\lambda},\frac{1}{2}\norm[\bfR^{-1}]{\bfA\bfx-\bfb}^{2} + \beta_{\lambda}\big), \\
    \delta\mid\bfb,\bfx,\lambda &\sim \Gamma\big(n/2 + \alpha_{\delta},\frac{1}{2}\norm[\bfQ^{-1}]{\bfx-\bfmu}^{2} + \beta_{\delta}\big), \\
    \bfx\mid\bfb,\lambda,\delta &\sim \calN\left( \bfx_{\rm cond}, \bfGamma_{\rm cond} \right), \label{eq:condGaus}
\end{align}
where $\bfGamma_{\rm cond}$ and $\bfx_{\rm cond}$ are 
defined in \eqref{eq:cond_x} and $\Gamma$ is the gamma distribution defined in \eqref{eq:hyperpriors}.
Here, $\bfx$ is drawn separately from $\lambda$ and $\delta$ to exploit the conditionally conjugate Gaussian distribution \eqref{eq:condGaus}.  
A block Gibbs algorithm for sampling from \eqref{eq:post_nonGaus} is outlined in \Cref{alg:gibbs} \cite{bardsley2018computational}.  
\begin{algorithm}[H]
        \caption{Block Gibbs algorithm for sampling the posterior density \eqref{eq:post_nonGaus} \cite{bardsley2018computational}}\label{alg:gibbs}
        \begin{algorithmic}[1]
        \Require{Number of samples $T$ and burn-in period $T_b$}
        \Ensure{Approximate samples from \eqref{eq:post_nonGaus}: $\left\{\bfx^t, \lambda^t, \delta^t \right\}_{t=T_b + 1}^T$}
        \State Initialize $\bfx^0, \ \lambda^0, \ \delta^0$
        \For{t = 1 to T}
        \State {Compute $\lambda^t \sim \Gamma\big(m/2 + \alpha_{\lambda},\frac{1}{2}\norm[\bfR^{-1}]{\bfA\bfx^{t-1}-\bfb}^{2} + \beta_{\lambda}\big)$}
        \State {Compute $\delta^t \sim \Gamma\big(n/2 + \alpha_{\delta},\frac{1}{2}\norm[\bfQ^{-1}]{\bfx^{t-1}-\bfmu}^{2} + \beta_{\delta}\big)$}
        \State {Compute $\bfx^t \sim \calN(\bfGamma_{\rm cond}^t (\lambda^t \bfA\t \bfR^{-1} \bfb + \delta^t \bfQ^{-1} \bfmu), \bfGamma_{\rm cond}^t )$, \label{alg:gibbsline5} \\ \qquad  \qquad where $\bfGamma_{\rm cond}^t = (\lambda^{t}\bfA\t \bfR^{-1} \bfA + \delta^{t} \bfQ^{-1})^{-1}$}
        \EndFor
        \end{algorithmic}
\end{algorithm}

The Gibbs sampler generates a Markov chain $\{\bfx^t,\lambda^t,\delta^t\}_{t=1}^T$ that converges in distribution to the posterior density $\pi(\bfx,\lambda,\delta\mid\bfb)$ \cite{bardsley2018computational}. However, the computational cost of this Gibbs sampler is prohibitive for problems with large $n$.  This is due to the fact that, even though the conditional distribution $\pi(\bfx\mid \bfb, \lambda,\delta)$ is Gaussian, drawing a sample requires the solution of an $n \times n$ linear system. For problems with small $n$, one could factorize $\bfA$ in an off-line phase, or if $\bfA$ has a specific structure (e.g., in image deblurring problems) then one could exploit a structured factorization using the Fourier transform \cite{nagydeblurring,fox2016fast}.  However, we emphasize that for many inverse problems (e.g., atmospheric inverse modeling and tomography), $\bfA$ is never formed but is \jmc{accessed purely} as a function evaluation, so iterative methods are the common approaches to use.  Early termination of an iterative method when used within a Gibbs sampler would produce an approximate sample from the conditional, hence motivating a Metropolis-Hastings approach \cite{brown2018low}. Moreover, the number of Gibbs samples increases with $n$ since the integrated autocorrelation time of the MCMC chain tends to $\infty$ \cite{saibaba2019efficient}.

For large-scale problems, computing a sample from \eqref{eq:condGaus} may be computationally infeasible, so to address this, previous approaches substitute direct sampling with a Metropolis-Hastings algorithm \cite{Metropolis_Rosenbluth_Rosenbluth_Teller_Teller_1953,Hastings_1970}. We follow this approach and replace the computation of $\bfx^t$ in line 5 of \Cref{alg:gibbs} with an accept-reject step where a sample is drawn from a proposal distribution, $\hat{\pi}_x(\bfx\mid\lambda,\delta,\bfb)$, which is an approximation of the conditional distribution \eqref{eq:condGaus}, and then the sample is accepted with some probability. In \cite{brown2018low}, a proposal distribution based on a low-rank approximation of the prior-preconditioned Hessian was used, where randomized \jmc{singular value decomposition (rSVD)} techniques were used to compute a low-rank approximation.  Such low-rank approximations were considered in \cite{ghattasinfinitebayes} and were combined with marginalization-based MCMC methods in \cite{saibaba2019efficient}. In this paper, we are interested in low-rank independence samplers that are based on the generalized Golub-Kahan bidiagonalization.  These will be discussed in \cref{sec:methods}.

First, we provide an overview of independence sampling.  Let $\bfx \in \bbR^n$ and denote the target density by $h(\bfx)$. The Metropolis-Hastings algorithm generates at iteration $t$ a sample $\bfx^\star$ from a proposal distribution, possibly conditioned on the current state $\bfx^{t-1}$, and sets $\bfx^t = \bfx^\star$ with probability $\min(1,\jmc{\rho})$ where, for fixed $\lambda$ and $\delta$,
$$\jmc{\rho} (\bfx^{t-1}, \bfx^\star) = \frac{h(\bfx^\star) q(\bfx^{t-1} \mid \bfx^\star)}{h(\bfx^{t-1}) q(\bfx^\star \mid \bfx^{t-1})}$$
where $q(\cdot \mid \bfx^{t-1})$ is the density of the proposal distribution. The algorithm generates a Markov chain $\{ \bfx^t \}$ that converges to the target distribution \cite{MCstatmethods}.

A Metropolis-Hastings independence sampler generates proposal states from a density that is independent of the current state of the chain, i.e., the proposal density has the form $q(\bfx^\star \mid \bfx^{t-1}) = g(\bfx^\star)$, and the ratio can now be written as 
\begin{equation}
\label{eq:ratio}
\frac{h(\bfx^\star) g(\bfx^{t-1})}{h(\bfx^{t-1}) g(\bfx^\star)} = \frac{w(\bfx^\star)}{w(\bfx^{t-1})}
\end{equation}
where $w(\bfx;\lambda,\delta) \propto \frac{h(\bfx)}{g(\bfx)}$.  Let $h(\bfx) = \pi (\bfx \mid \bfb, \lambda, \delta)$ from \eqref{eq:condGaus} where the conditional variables have been dropped, and let $g(\bfx)$ be a proposal density function (to be defined later). 

An independence Metropolis-Hastings within Gibbs algorithm for sampling the posterior density \eqref{eq:post_nonGaus} would look similar to \Cref{alg:gibbs} but instead of drawing a sample from the exact distribution in \cref{alg:gibbsline5}, a sample is drawn from some proposal $g(x)$ and accepted with probability $\min(1,\jmc{\rho})$ where $\jmc{\rho} = \frac{w(\bfx^\star)}{w(\bfx^{t-1})}$. If rejected, $\bfx^t = \bfx^{t-1}$.


\section{Generalized Golub-Kahan based proposals for independence sampling}
\label{sec:methods}

The target conditional Gaussian density function is given by
\begin{equation}
\label{eq:target}
h(\bfx) := \frac{1}{\sqrt{(2 \pi)^n |\bfGamma_{\rm cond}|}} \exp \left(-\frac{1}{2} (\bfx - \bfx_{\rm cond})\t \bfGamma_{\rm cond}^{-1} (\bfx - \bfx_{\rm cond}) \right),
\end{equation}
and samples from the distribution $\calN(\bfx_{\rm cond}, \bfGamma_{\rm cond})$ can be generated as $\bfx = \bfx_{\rm cond} + \bfG \bfepsilon$ where $\bfGamma_{\rm cond} = \bfG \bfG\t$ and $\bfepsilon \sim \calN(\bfzero, \bfI)$.  However, it is computationally expensive to compute $\bfx_{\rm cond}$ and the product $\bfG \bfepsilon$.  Works such as \cite{brown2018low,ghattasinfinitebayes} exploit the low-rank structure of the forward operator $\bfA$ and rely on factorizations of the prior covariance $\bfQ^{-1} = \bfL \bfL\t$ to get efficient representations for $\bfGamma_{\rm cond}$ that can exploit the fast decay of singular values.  These approaches may not be computationally feasible.  For example, in atmospheric emissions tomography, the forward model matrices are nowhere near low-rank and, moreover, prior covariance matrices are defined using complicated spherical covariance kernels that allow seasonal changes in variability \cite{cho2022computationally}.
The approach we describe herein can handle these scenarios.

Specifically, we use a Krylov subspace method based on the generalized Golub-Kahan (genGK) bidiagonalization process to approximate $\bfx_{\rm cond}$.  Then we consider two proposal distributions, one which uses the resulting genGK matrices to form an approximation of $\bfGamma_{\rm cond}$ and another which uses preconditioned Krylov sampling.  Similar approximations were considered in \cite{UQlargebayesian,chung2018efficient}, but the main difference in this work is that we will use the genGK approximations to define a proposal distribution $g(\bfx)$ that approximates the target distribution  \eqref{eq:target}.  To the best of our knowledge, the use of genGK low-rank approximations within MCMC sampling approaches has not been explored. We begin in \Cref{sub:genGK} with a brief overview of the genGK approach, followed by details of the genGK approximation to the target distribution in \Cref{sub:sampling_prop}. We provide details regarding the Metropolis-Hastings within Gibbs algorithm with the genGK approximation used for the proposal distribution.  For many problems, the genGK approximate distribution provides a computationally efficient approach for generating proposal samples. However, for problems where the genGK proposal distribution may require very large ranks to obtain a sufficient approximation, we propose in \Cref{sub:MH_genGK} an alternative proposal sampler that is based on efficient preconditioned Krylov methods and consider its use in Metropolis-Hastings within Gibbs sampling.

\subsection{Generalized Golub-Kahan methods}
\label{sub:genGK}
The genGK bidiagonalization process is an iterative Krylov subspace projection method that was developed in \cite{arioli2013generalized} and can be used to efficiently compute general-form Tikhonov solutions \eqref{eq:Tik} \cite{ChungSaibabaHybrid2017}.  The genGK approach is well suited for problems where \jmc{mat-vecs} with $\bfA,$ $\bfA\t$, and $\bfQ$ can be done efficiently, but $\bfQ^{-1}$ or any factorization of $\bfQ$ is expensive. We assume that the inverse and square root of $\bfR$ are computationally feasible, e.g, an identity or diagonal matrix.  For the prior covariance matrix $\bfQ$, we focus on problems where explicit computation of the square root and inverse of the covariance matrix $\bfQ$ are not computationally feasible.  This scenario often arises for large-scale problems where $\bfQ$ is highly structured or is constructed from covariance kernels, possibly on unstructured grids. We focus on the Mat\'ern class of covariance kernels, where \jmc{a mat-vec} multiplication with $\bfQ$ can be performed easily, e.g., in $\calO(n \log n)$ time by exploiting fast Fourier transforms if the solution is represented on a uniform equi-spaced grid. For such covariance matrices, a symmetric factorization of $\bfQ^{-1} = \bfL\t \bfL$ is not available, and thus it is not possible to reformulate the regularization term as $\|\bfL (\bfx - \bfmu)\|_2^2$, i.e., transforming to standard-form.  Instead, we develop methods that work directly with $\bfQ$.

With the change of variables $\bfy = \bfQ^{-1}(\bfx - \bfmu)$, we get the
equivalent problem,
    \begin{equation}
        \min_{\bfy} \frac{\lambda}{2}\norm[\bfR^{-1}]{\bfA\bfQ\bfy - (\bfb - \bfA \bfmu)}^{2} + \frac{\delta}{2}\norm[\bfQ]{\bfy}^{2}, \label{eq:trans}
    \end{equation}
which is projected onto subspaces of increasing dimension in an iterative Krylov projection process.  That is, given matrices $\bfA, \bfR,$ and $\bfQ$ and vectors $\bfb$ and $\bfmu$, 
let $\gamma_1 \bfu_1 = \bfb - \bfA \bfmu$ and $\jmc{\zeta}_1 \bfv_1 = \bfA\t \bfR^{-1} \bfu_1,$ \jmc{where $\gamma_1 = \| \bfb - \bfA \bfmu\|_2$ and $\zeta_1 = \|\bfA\t \bfR^{-1} \bfu_1\|_2$.}  At the $k$th iteration of the genGK process, we generate vectors $\bfu_{k+1}$ and $\bfv_{k+1}$ such that
\begin{equation*}
\gamma_{k+1} \bfu_{k+1} = \bfA \bfQ \bfv_{k} - \jmc{\zeta}_k \bfu_k, \qquad \jmc{\zeta}_{k+1} \bfv_{k+1} = \bfA\t \bfR^{-1}\bfu_{k+1} - \gamma_{k+1} \bfv_k,
\end{equation*}
where after $k$ iterations, we have matrices 
\begin{equation*}
   \bfU_{k+1}=[\bfu_1,\ldots,\bfu_{k+1}], \bfV_{k+1}=[\bfv_1,\ldots,\bfv_{k+1}], \,\, \mbox{and} \,\, \bfB_{k} = 
    \begin{bmatrix}
        \jmc{\zeta}_1 & & \\
        \gamma_2 & \ddots &  \\
        & \ddots & \jmc{\zeta}_k\\
        & & \gamma_{k+1}
    \end{bmatrix}
\end{equation*}
that satisfy the following relationships, in exact arithmetic,
\begin{equation}
    \bfA\bfQ\bfV_{k} = \bfU_{k+1}\bfB_{k} \quad \mbox{and} \quad \bfA\t\bfR^{-1}\bfU_{k+1} = \bfV_{k}\bfB_{k}\t + \jmc{\zeta}_{k+1}\bfv_{k+1}\bfe_{k+1}\t, \label{eq:gk1}
\end{equation}
        with 
        \begin{equation}
    \bfU_{k+1}\t\bfR^{-1}\bfU_{k+1}=\bfI_{k+1} \quad \text{and} \quad \bfV_{k}\t\bfQ\bfV_{k}=\bfI_{k}. \label{eq:gk2}
        \end{equation}
The vector $\bfe_i$ is the $i$th column of the identity matrix of the appropriate size.

The genGK process constructs a basis for the Krylov subspaces,
    $\calR(\bfV_k) = \calK_k(\bfA\t\bfR^{-1}\bfA\bfQ, \bfA\t\bfR^{-1}\bfb)$
and 
  $  \calR(\bfU_k) = \calK_k(\bfA \bfQ \bfA\t \bfR^{-1}, \bfb), $
where $\calK_{k}(\bfM,\bfg) = \Span{\bfg,\bfM\bfg,\ldots, \bfM^{k-1}\bfg}$ and $\calR$ denotes the column space.  At the $k$th iteration, we have an approximate solution for \eqref{eq:Tik}, given by
 $\bfx_{k} = \bfmu + \bfQ \bfV_{k}\bfz_{k}$, where $\bfz_k$ solves the projected problem.  That is,
\begin{equation}
    \min_{\bfy_{k}\in\calR(\bfV_{k})}\frac{\lambda}{2}\norm[\bfR^{-1}]{\bfA\bfQ\bfy_{k}-(\bfb-\bfA \bfmu)}^{2}+\frac{\delta}{2}\norm[\bfQ]{\bfy_{k}}^{2} \quad \Longleftrightarrow \quad \min_{\bfz_{k}\in\bbR^{k}}\frac{\lambda}{2}\norm[2]{\bfB_{k}\bfz_{k}-\gamma_{1}\bfe_{1}}^{2} + \frac{\delta}{2}\norm[2]{\bfz_{k}}^{2} .\label{eq:proj}
\end{equation}

Note that by using the genGK approach, we avoid $\bfQ^{-1}$ and rely on projections of the original problem to obtain a reconstruction in $k$ iterations ($k \ll n$).    
Moreover, by using the genGK relations, we can define oblique projectors, 
$$\calP_{\bfV_k} = \bfV_k \bfV_k\t \bfQ \qquad \mbox{and} \qquad \calP_{\bfU_{k+1}} = \bfU_{k+1} \bfU_{k+1}\t \bfR^{-1},$$
and we can build a low-rank approximation for $\bfA$ as $$\bfA \approx \bfA \calP_{\bfV_k}\t = \bfU_{k+1} \bfB_k \bfV_k\t \equiv \widehat \bfA .$$
Such approximations are related to the generalized singular value decomposition and were used in the context of hyperparameter estimation in \cite{hall2024efficient}. As described in \cite{ChungSaibabaHybrid2017}, the generalized singular values are related to the singular values of $\bfR^{-1/2} \bfA \bfQ^{1/2}$ and tend to exhibit faster decay compared to the singular values of $\bfA$. Next we will use the genGK approximation to build a proposal distribution that can be used within a Metropolis-Hastings within Gibbs approach for sampling from \eqref{eq:post_nonGaus}.  A key feature to note is that the genGK bidiagonalization process does not depend on $\lambda$ and $\delta$.  That is, all of the resulting matrices $\bfV_k, \bfB_k, \bfU_{k+1}$ are independent of $\lambda$ and $\delta$ and can be reused if these parameters change. We will exploit this property in \Cref{sub:sampling_prop} for efficient sampling from the genGK proposal distribution and again in \Cref{sub:MH_genGK} for efficient conditional mean estimation.

\subsection{genGK approximation to the target distribution}
\label{sub:sampling_prop}
Now consider the target conditional distribution \eqref{eq:target} with conditional mean $\bfx_{\rm cond}$ and covariance matrix $\bfGamma_{\rm cond}$.  
Recall that after $k$ iterations of the genGK process, we have $\bfx_k \approx \bfx_{\rm cond}$ and matrices $\bfU_{k+1}, \ \bfV_{k},$ and $\bfB_{k}$ that satisfy the relations in \eqref{eq:gk1} and \eqref{eq:gk2}. Following a similar approach as in \cite{UQlargebayesian}, we aim to use these matrices to approximate $\bfGamma^{1/2}$. Notice that one can factorize the covariance matrix as $\bfGamma_{\rm cond} = \bfGamma_{\rm cond}^{1/2}\bfGamma_{\rm cond}^{1/2}$, where
\begin{equation}
    \bfGamma_{\rm cond}^{1/2} = \delta^{-1/2}\bfQ^{1/2}\left( \bfI + \frac{\lambda}{\delta}\bfQ^{1/2}\bfA\t\bfR^{-1}\bfA\bfQ^{1/2} \right)^{-1/2}.
\end{equation}
Although we write out such a factorization for exposition purposes, we will not compute $\bfQ^{1/2}$ explicitly.  We will use a Lanczos algorithm to access it in a matrix-free fashion, following the procedure outlined in the appendix \Cref{sec:lowrank}, and it will be part of a preprocessing step (that is independent of MCMC sampling).

Prior to sampling, we compute the low-rank representation, $\bfQ^{1/2}\bfV_{k}\bfB_{k}\t\bfB_{k}\bfV_{k}\t\bfQ^{1/2} = \rev{\bfP_k}\bfTheta_k \rev{\bfP_k}\t$.
Then using the genGK approximation $\bfA\t \bfR^{-1}\bfA \approx \widehat \bfA\t \bfR^{-1}\widehat \bfA = \bfV_k \bfB_k \t \bfB_k \bfV_k\t$, we can define the matrix approximation,
\begin{align*}
    \widehat \bfGamma_{\rm cond} & = (\lambda \bfV_k \bfB_k\t \bfB_k \bfV_k\t + \delta \bfQ^{-1})^{-1}\\
    & = \bfQ^{1/2}(\lambda \underbrace{\bfQ^{1/2} \bfV_k \bfB_k\t \bfB_k \bfV_k\t \bfQ^{1/2}}_{\rev{\bfP_k} \bfTheta_k \rev{\bfP_k}\t} + \delta \bfI)^{-1} \bfQ^{1/2}.
\end{align*}
Thus, the square root of the conditional covariance matrix can be approximated as,
\begin{align*}
        \bfGamma_{\rm cond}^{1/2} 
        &= \delta^{-1/2}\bfQ^{1/2}\left( \frac{\lambda}{\delta}\bfQ^{1/2}\underbrace{\bfA\t\bfR^{-1}\bfA}_{\approx \bfV_{k}\bfB_{k}\t\bfB_{k}\bfV_{k}\t}\bfQ^{1/2} +\bfI \right)^{-1/2} \\
        &\approx \delta^{-1/2}\bfQ^{1/2}\left(  \frac{\lambda}{\delta} \underbrace{\bfQ^{1/2}\bfV_{k}\bfB_{k}\t\bfB_{k}\bfV_{k}\t\bfQ^{1/2}}_{\approx \bfP_k\bfTheta_k\bfP_k\t}+ \bfI  \right)^{-1/2}\\
        & = \delta^{-1/2}\bfQ^{1/2}(\bfI - \bfP_{k}\bfD_{k}\bfP_{k}\t) \equiv \widehat{\bfGamma}_{\rm cond}^{1/2}
    \end{align*}
where  $\bfD_{k} \equiv \bfI_{k} - (\bfI_{k} + \frac{\lambda}{\delta}\mathbf{\Theta}_{k} )^{-1/2}$. Finally, we can define the proposal distribution $\calN(\bfx_k, \widehat \bfGamma_{\rm cond})$, i.e.,
\begin{equation}
\label{eq:proposal_genGK}
g_1(\bfx) := \frac{1}{\sqrt{(2 \pi)^n |\widehat \bfGamma_{\rm cond}|}} \exp \left(-\frac{1}{2} (\bfx - \bfx_{k})\t \widehat \bfGamma_{\rm cond}^{-1} (\bfx - \bfx_{k}) \right),
\end{equation}
and draw a sample,
\begin{align}
    \bfx^{\star} &= \bfx_k + \widehat{\bfGamma}_{\rm cond}^{1/2}\bfxi \\
    &= \bfmu+ \bfQ\bfV_{k}\bfz_{k} + \delta^{-1/2}\bfQ^{1/2}\left(\bfI - \bfP_{k}
\bfD_{k}\bfP_{k}\t\right)\bfxi \label{eq:propsample_genGK}
\end{align}
for $\bfxi\sim\calN(\bf0,\bfI)$. 
Theoretical results regarding the accuracy of the approximate distribution to the true conditional distribution can be found in \cite{UQlargebayesian}.  

Next we incorporate the genGK approximation above as a proposal distribution within a Metropolis-Hastings within Gibbs approach to sample from \eqref{eq:post_nonGaus}. We need to investigate the acceptance ratio.  Let $\bfx$ denote the current state of the chain and let $\bfx^\star$ denote the proposed state.  Then from Proposition 1 of \cite{brown2018low}, we have that the acceptance ratio can be computed as $\jmc{\rho}_1 = \frac{w(\bfx^\star)}{w(\rev{\bfx^{t-1}})}$ where $$w(\bfx) = \exp \left(- \frac{1}{2}\bfx\t \left(\bfGamma_{\rm cond}^{-1} - \widehat \bfGamma_{\rm cond}^{-1} \right) \bfx\right).$$
Note that at each iteration, we have the weight from the previous iteration $w(\bfx^{t-1})$ but we must compute the weight $w(\bfx^\star)$.  An efficient implementation of this can be obtained by observing that
\begin{align}
\log w(\rev{\bfx}) & = - \frac{1}{2}\rev{\bfx}\t \left(\bfGamma_{\rm cond}^{-1} - \widehat \bfGamma_{\rm cond}^{-1} \right) \rev{\bfx} \\
& = - \frac{1}{2}\rev{\bfx}\t \left(\lambda\bfA\t \bfR^{-1} \bfA + \delta \bfQ^{-1}  - \bfQ^{-1/2} (\lambda \bfP_k \bfTheta_k \bfP_k\t + \delta \bfI) \bfQ^{-1/2} \right) \rev{\bfx} \\
& = - \frac{\lambda}{2}\rev{\bfx}\t \left(\bfA\t \bfR^{-1} \bfA - \bfQ^{-1/2} \bfP_k \bfTheta_k \bfP_k\t\bfQ^{-1/2} \right) \rev{\bfx}.
\end{align}
Similarly, the quality of the low-rank approximation to the target distribution will be evident from the acceptance ratio. The \jmc{Metropolis-Hastings} within Gibbs algorithm with the genGK approximation used for proposal sampling is summarized in \Cref{alg:gibbs_genGK}.  Notice that mat-vecs with $\bfA$ and $\bfA\t$ are required in the genGK process (i.e., $k$ multiplications with $\bfA$ and $\bfA\t$) and again in the computation of the acceptance ratio.  Also, for each new sample of $\lambda^t$ and $\delta^t$, $\bfB_k$ and $\bfV_k$ are reused for efficient computation of $\bfx_k$, and $\bfP_k$ and $\bfTheta_k$ are reused for the efficient computation of the proposal sample. \rev{From Proposition 3 in \cite{brown2018low}, the subchain from \Cref{alg:gibbs_genGK} has stationary distribution $h(\bfx)=\pi(\bfx\mid\bfb,\lambda,\delta)$ and is uniformly ergodic \cite{MCstatmethods}.}

\begin{algorithm}[H]
        \caption{Metropolis-Hastings within Gibbs with genGK approximation }\label{alg:gibbs_genGK}
           \begin{algorithmic}[1]
        \Require{$\bfA, \bfb, \bfQ, \bfR$, number of genGK iterations $k$}
        \State Run $k$ iterations of genGK to get bidiagonal matrix $\bfB_k$ and basis vectors $\bfV_k$
        \State Initialize $\lambda^{0}, \delta^0,$ and $\bfx^0$
        \State Precompute $\bfP_k$, $\bfTheta_k$
        \For{t = 1 to T}
        \State {Compute $\lambda^t \sim \pi(\lambda \mid \bfb, \bfx^{t-1})$} and $\delta^t \sim \pi(\delta \mid \bfb, \bfx^{t-1})$
        \State {Compute $\bfx_k$ with fixed $\lambda^t,\ \delta^t$}
        \State {Compute proposal sample $\bfx^\star \sim g_1(\bfx)$ as in \eqref{eq:propsample_genGK}}
        \State {Accept $\bfx^t = \bfx^\star$ with probability $\min(1,\jmc{\rho}_1)$ where
        $\jmc{\rho}_1 = \frac{w(\bfx^\star)}{w(\bfx^{t-1})}$}
        \State {Otherwise set $\bfx^t = \bfx^{t-1}$}
        \EndFor
        \end{algorithmic}
\end{algorithm}

\subsection{Proposal sampling using a preconditioned Lanczos method}
\label{sub:MH_genGK}
 Next we consider an alternative proposal distribution, where we generate an approximate sample from $\calN(\bfx_k, \bfGamma_{\rm cond})$ by using a preconditioned Lanczos method. This approach follows that of Method 2 in \cite{UQlargebayesian}. First rewrite the covariance matrix as 
$$
\bfGamma_{\rm cond} = \left(\delta\bfQ^{-1} + \lambda\bfA\t\bfR^{-1}\bfA\right)^{-1} = \lambda\bfQ\bfF^{-1}\bfQ
$$
where $$\bfF = \frac{\delta}{\lambda}\bfQ + \bfQ\bfA\t\bfR^{-1}\bfA\bfQ. $$
Next define $$\bfS_F = \lambda^{-1/2}\bfQ\bfF^{-1/2}$$ which is a square root matrix of $\bfGamma_{\rm cond}$, that is $\bfGamma_{\rm cond}=\bfS_F\bfS_F\t$. Let $\bfG$ be a preconditioner satisfying $\bfG\bfG\t \approx \bfF^{-1}$ which allows us to compute a square root of $\bfF$. Now we have the exact factorization
$$
\bfGamma_{\rm cond} = \bfS_F\bfS_F\t \qquad \bfS_F = \lambda^{-1/2}\bfQ\bfG\t\left(\bfG\bfF\bfG\t\right)^{-1/2}.
$$

Then we draw the sample as
\begin{equation}
\label{eq:method2prop}
\bfx^{\star} = \bfx_k + \bfS_F\bfxi 
\end{equation}
for $\bfxi\sim\calN(\bf0,\bfI)$. It should be noted that while the factorization $\bfGamma_{\rm cond}=\bfS_F\bfS_F\t$ is exact, the mean $\bfx_k$ is the genGK solution computed after $k$ iterations.  Thus, the sample drawn is not from the exact target conditional density function $h(\bfx)$ but a proposal distribution $\calN(\bfx_k,\bfGamma_{\rm cond})$ with a probability density given by
$$
g_2(\bfx) := \frac{1}{\sqrt{(2 \pi)^n | \bfGamma_{\rm cond}|}} \exp \left(-\frac{1}{2} (\bfx - \bfx_{k})\t \bfGamma_{\rm cond}^{-1} (\bfx - \bfx_{k}) \right). 
$$
As shown in the appendix \Cref{sec:accept_norm_precond}, the acceptance ratio for $g_2(\bfx)$ can be evaluated efficiently as
\begin{align*}
\jmc{\rho}_2 (\bfx^{t-1}, \bfx^\star) &= \frac{h(\bfx^\star) g_2(\bfx^{t-1})}{h(\bfx^{t-1}) g_2(\bfx^\star)} \\
&= \exp\left(\left( \bfx^{\star} - \bfx^{t-1} \right)\t \left( \lambda\left(\gamma_1\jmc{\zeta}_1\bfv_1 - \bfV_{k+1}\begin{bmatrix}
    \bfB_k\t \\ \jmc{\zeta}_{k+1}\bfe_{k+1}\t
\end{bmatrix} \bfB_k\bfz_k\right)- \delta\bfV_k\bfz_k\right)\right).
\end{align*}

Lastly, to draw the proposal sample \eqref{eq:method2prop}, the matrix $\bfS_F$ is never actually formed.  Instead applications of the matrix $\bfG\t(\bfG\bfF\bfG\t)^{-1/2}$ on a vector $\bfxi$ are performed using a preconditioned Lanczos approach.  Note that a mat-vec with $\bfF$ requires one mat-vec with $\bfA$ and $\bfA\t$ and two mat-vecs with $\bfQ$.  For computing the acceptance ratio, we can avoid further mat-vecs with $\bfA$ and $\bfA\t$ by exploiting the \jmc{genGK} relationship.
This approach can be computationally expensive, compared to the approach described in
\Cref{sub:sampling_prop}, especially if the forward operator is expensive; however, it serves as a good alternative in cases where the goal is to get a few, very good samples from the proposal density.

\begin{algorithm}[H]
        \caption{Metropolis-Hastings within Gibbs with preconditioned Lanczos}\label{alg:gibbs_precond}
           \begin{algorithmic}[1]
        \Require{$\bfA, \bfb, \bfQ, \bfR, \bfG$, number of genGK iterations $k$}
        \State Run $k$ iterations of genGK to get bidiagonal matrix $\bfB_k$ and basis vectors $\bfV_k$
        \State Initialize $\lambda^{0}, \delta^0,$ and $\bfx^0$
        \For{t = 1 to T}
        \State {Compute $\lambda^t \sim \pi(\lambda \mid \bfb, \bfx^{t-1})$} and $\delta^t \sim \pi(\delta \mid \bfb, \bfx^{t-1})$
        \State {Compute $\bfx_k$ with fixed $\lambda^t,\ \delta^t$}
        \State {Compute proposal sample $\bfx^\star \sim g_2(\bfx)$ as in \eqref{eq:method2prop}} with preconditioner $\bfG$
        \State {Accept $\bfx^t = \bfx^\star$ with probability $\min(1,\jmc{\rho}_2)$}
        \State {Otherwise set $\bfx^t = \bfx^{t-1}$}
        \EndFor
        \end{algorithmic}
\end{algorithm}

\section{Numerical results}
\label{sec:numerics}
In this section, we provide multiple examples to demonstrate the effectiveness and applicability of the proposed sampling approaches. We begin in \Cref{sub:seismic} with a seismic tomography example, where the problem is small enough so that comparisons to existing methods are feasible.  Then in \Cref{sub:atmospheric} we consider an atmospheric tomography example where comparisons to existing methods are not feasible, but the genGK approximation provides good proposal samples.  Finally, in \Cref{sub:PAT} we consider a dynamic photoacoustic tomography problem, where the goal is to reconstruct a spatiotemporal image.  This is a highly underdetermined problem, where the number of unknowns is very large, so the prior is critical.

In all of the examples, three chains, each with $T$ samples, are run in parallel.
We use a wavelet based estimation of the noise covariance parameter $\lambda_{est}$ \cite{donoho1995noising} and an estimate of  $\delta_{est}$
using the genHyBR method \cite{ChungSaibabaHybrid2017} with the weighted \jmc{generalized cross-validation} method on the projected problem. Each chain is initialized using a random number within an interval centered at $\lambda_{est}$ and $\delta_{est}$. \rev{Following \cite{bardsley2018computational}, the parameters defining the shape and rate of the gamma hyperpriors in \Cref{eq:hyperpriors} are set to be $\alpha_\lambda = \alpha_\delta =1 $ and $\beta_\lambda = \beta_\delta = 10^{-4}$.}

To ensure that there is no bias in the results, we remove samples from the initial stage of the MCMC chain. During this stage, called burn-in, the samples are moving from the starting position to a region that has a higher probability of being the target distribution \cite{bardsley2018computational}. For all examples, the first $50\%$ of samples from the initial burn-in stage are removed before any analysis is done on the mean and variance of $\{\bfx^t\}$ or the hyperparameter chains. 

Next, to examine whether the MCMC chain is in equilibrium, we perform the Geweke test on the individual $\lambda$ and $\delta$ chains. Here, we compare the mean of the first $10\%$ of samples, denoted $\mu_{10}$, to the mean of the last $50\%$, denoted $\mu_{50}$, and either accept or reject the null hypothesis, $H_0: \mu_{10}=\mu_{50}$. For $p$ values close to 1, there is strong evidence that the chain is in equilibrium. 
\jmc{As a convergence diagnostic,}
we use the Gelman-Rubin statistic, denoted as $\widehat{R}$, which compares the variance between multiple chains run in parallel. \jmc{A desired value for} $\widehat{R}$ is close to 1 \cite{gelman_rubin_diag_1992}. The cutoff is commonly set at 1.01, but in practice it ranges from 1.003 to 1.1 \cite{Vats_Knudson_2021}.

Finally, we test the independence of the individual $\lambda$ and $\delta$ chains through autocorrelation functions (ACF), as described in \cite{bardsley2018computational}. If the ACF decays to 0 fast enough, this means that there is no correlation between the samples in the individual chains. \rev{To estimate the number of independent samples,} we \rev{compute} the effective sample size (ESS)\rev{, $T_{\rm ESS} =\frac{T}{\tau_{\rm int}},$ where $\tau_{\rm int}$ is the integrated autocorrelation time approximated using Sokal’s adaptive truncated periodogram-estimator \cite{Sokal1997}.}

\subsection{Seismic tomography}
\label{sub:seismic}
\begin{figure}[H]
    \centering
    \includegraphics[width=0.3\linewidth]{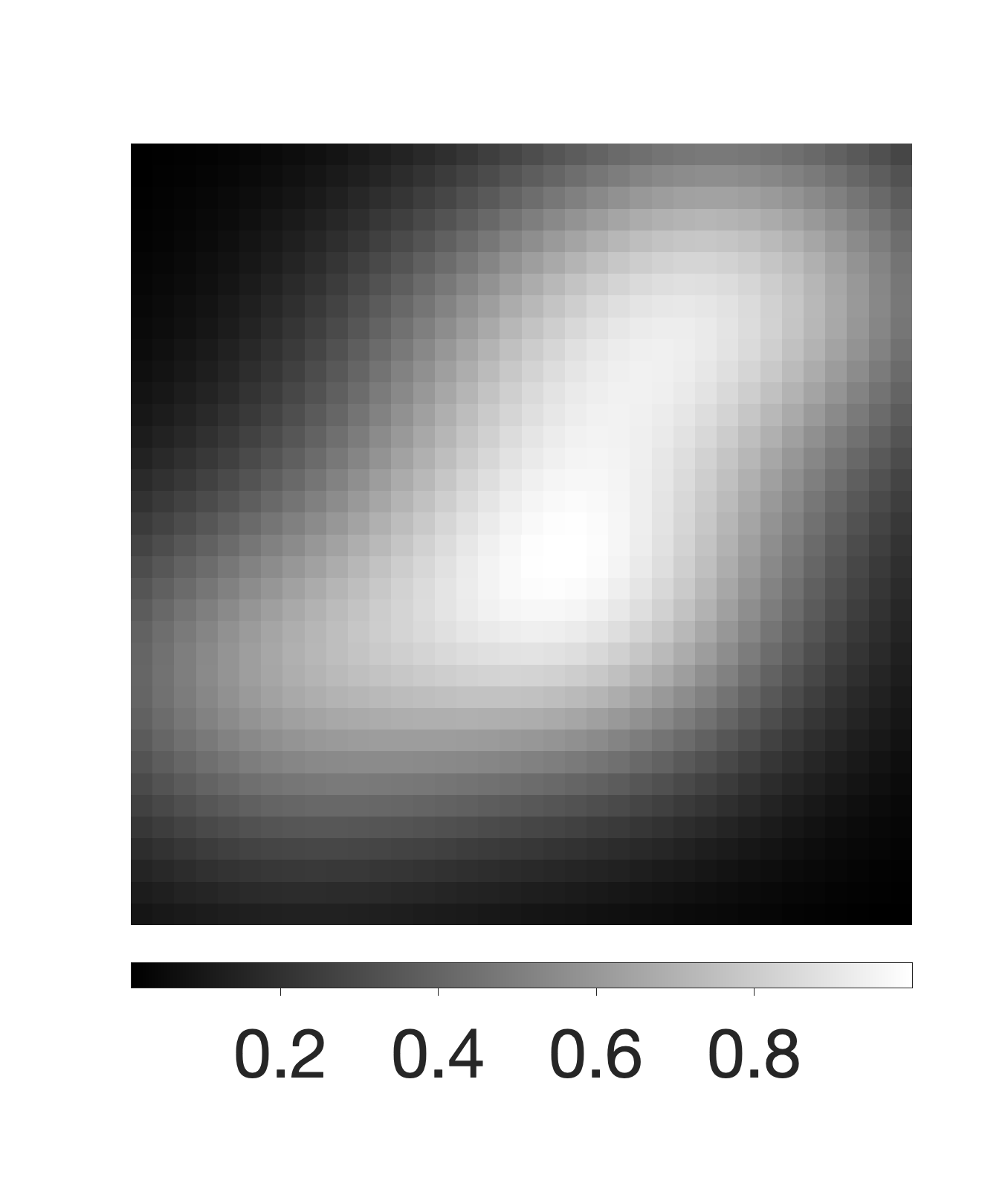}
    \caption{True image for $36\times36$ \texttt{PRseismic} seismic tomography test problem.}
    \label{fig:trueTomo}
\end{figure}

In this example, we use the \texttt{PRseismic} test problem from IRTools \cite{gazzola2019ir} with default settings. The goal is to approximate $\bfx\in\bbR^{1,296}$, the vectorized $36\times 36$ true image shown in \Cref{fig:trueTomo}, given measurements $\bfb\in\bbR^{1,800}$ and the forward model matrix $\bfA\in\bbR^{1,800\times 1,296}$ representing seismic tomography.  The observations contain $2\%$ Gaussian white additive noise, i.e., $\frac{\sigma\norm[2]{\bfxi}}{\norm[2]{\bfA\bfx_{\rm true}}}=0.02$ where $\sigma$ is the standard deviation and $\bfxi\sim\calN\left(\bf0,\bfI\right)$. Using this setup, the hyperparameter associated with the noise, $\lambda$, should be $\frac{1}{\sigma^2}\approx 7.11$. Additionally, the prior covariance matrix $\bfQ$ represents a Mat\'ern kernel with $\nu=1/2$ and $\ell = 1/4$. 

\begin{figure}[bthp]
    \centering
    \begin{tabular}{c c c c}
         & \textbf{genGK} & \textbf{tSVD} & \textbf{rSVD}  \\
        \raisebox{1.7cm}{\textbf{Mean}}
        & \includegraphics[width=.25\textwidth]{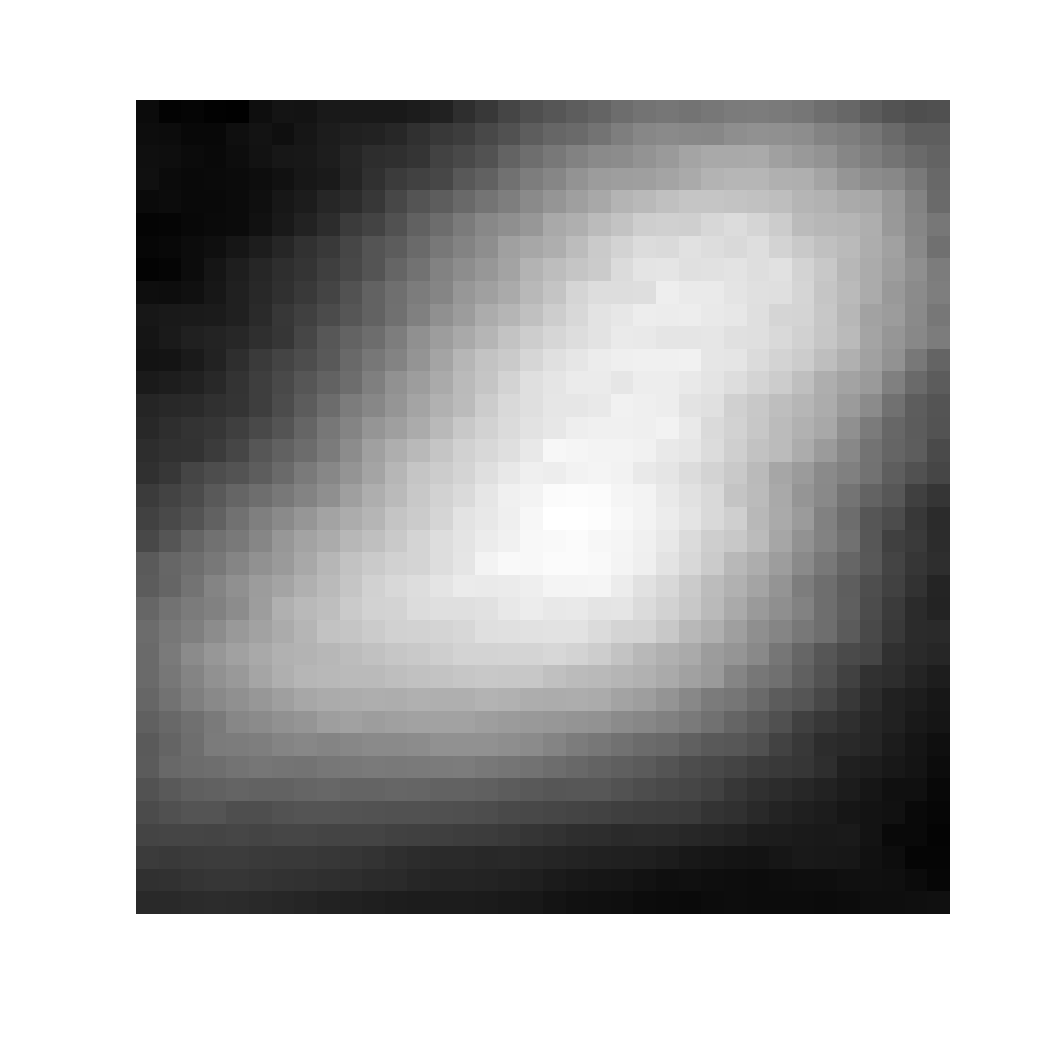}
        & \includegraphics[width=.25\textwidth]{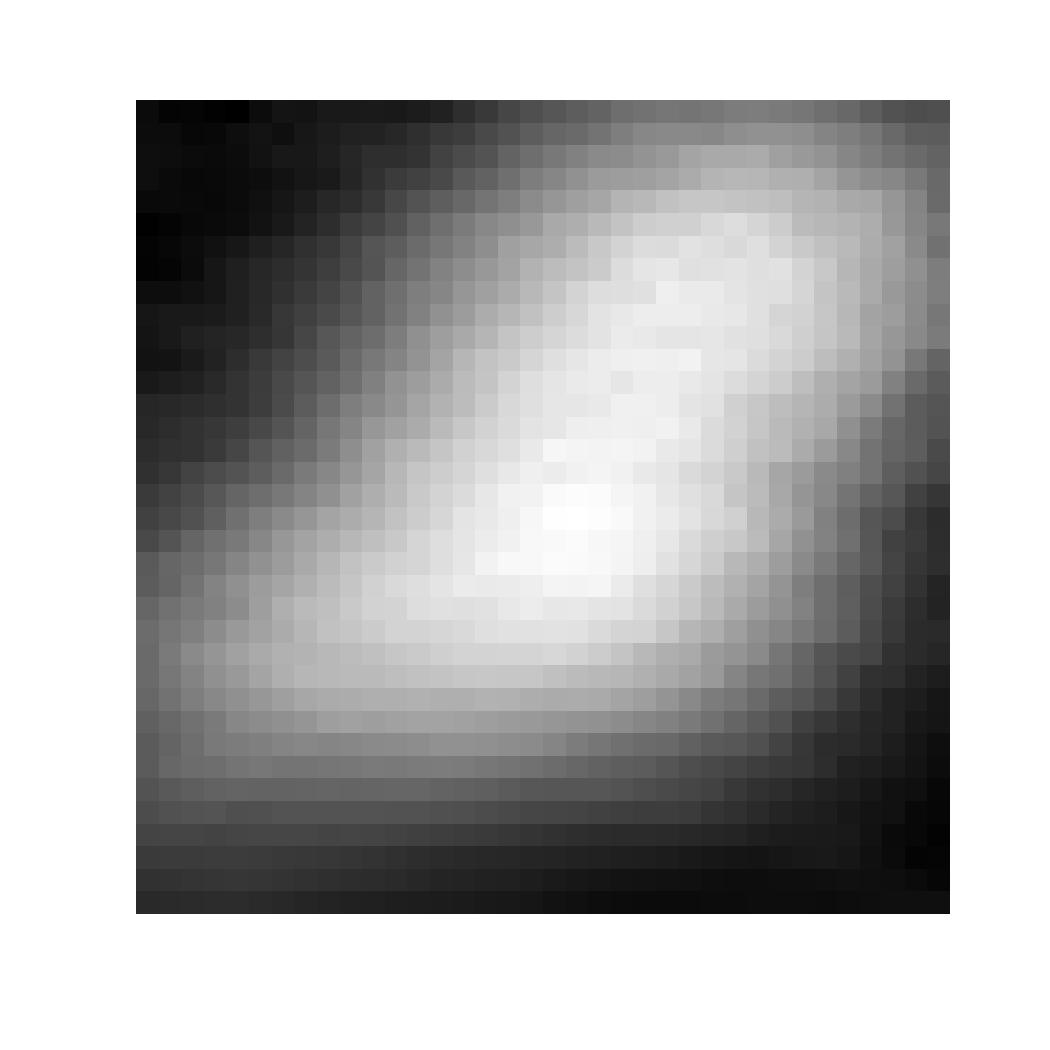}
        & \includegraphics[width=.277\textwidth]{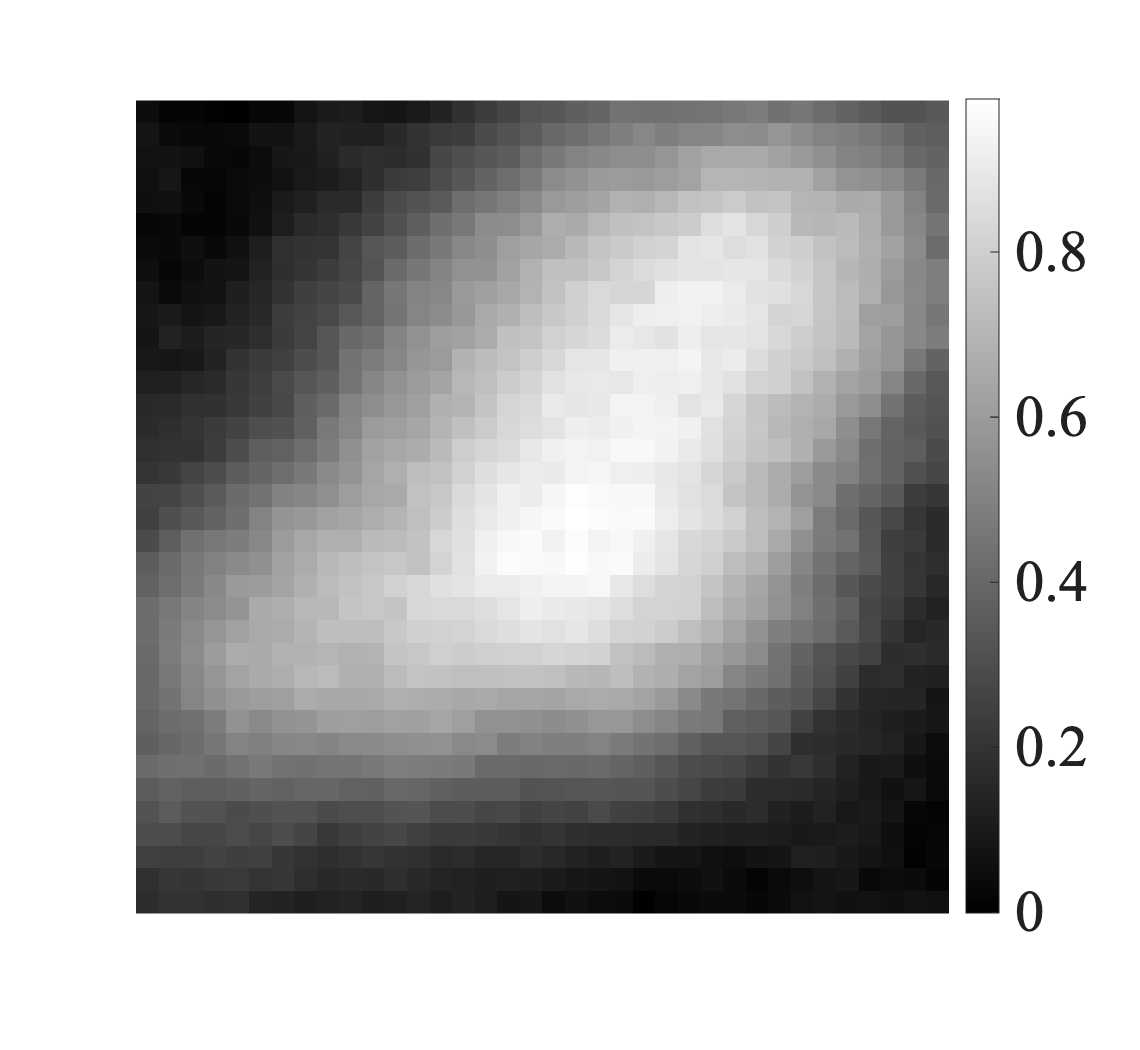}
        \\
        \raisebox{1.7cm}{\textbf{Variance}}
        & \includegraphics[width=.25\textwidth]{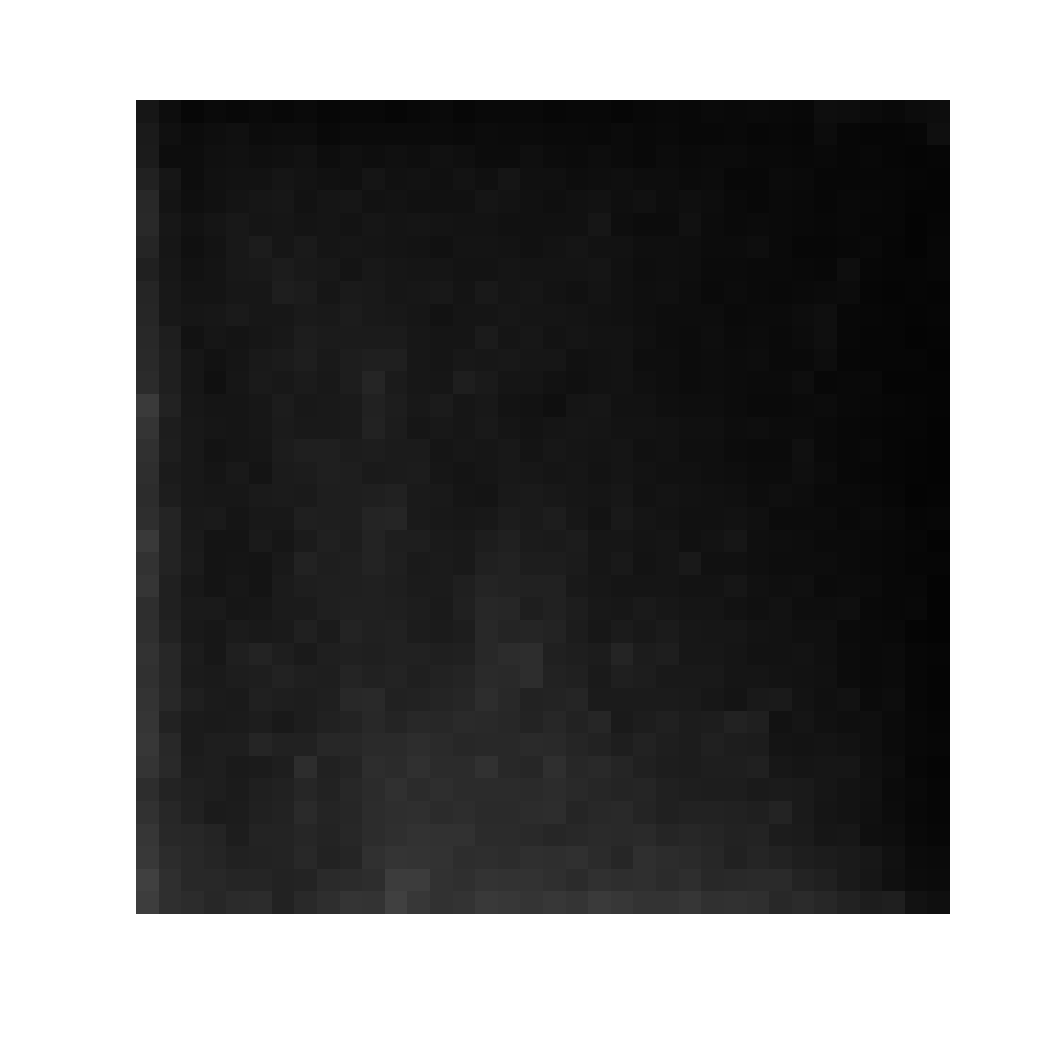}
        & \includegraphics[width=.25\textwidth]{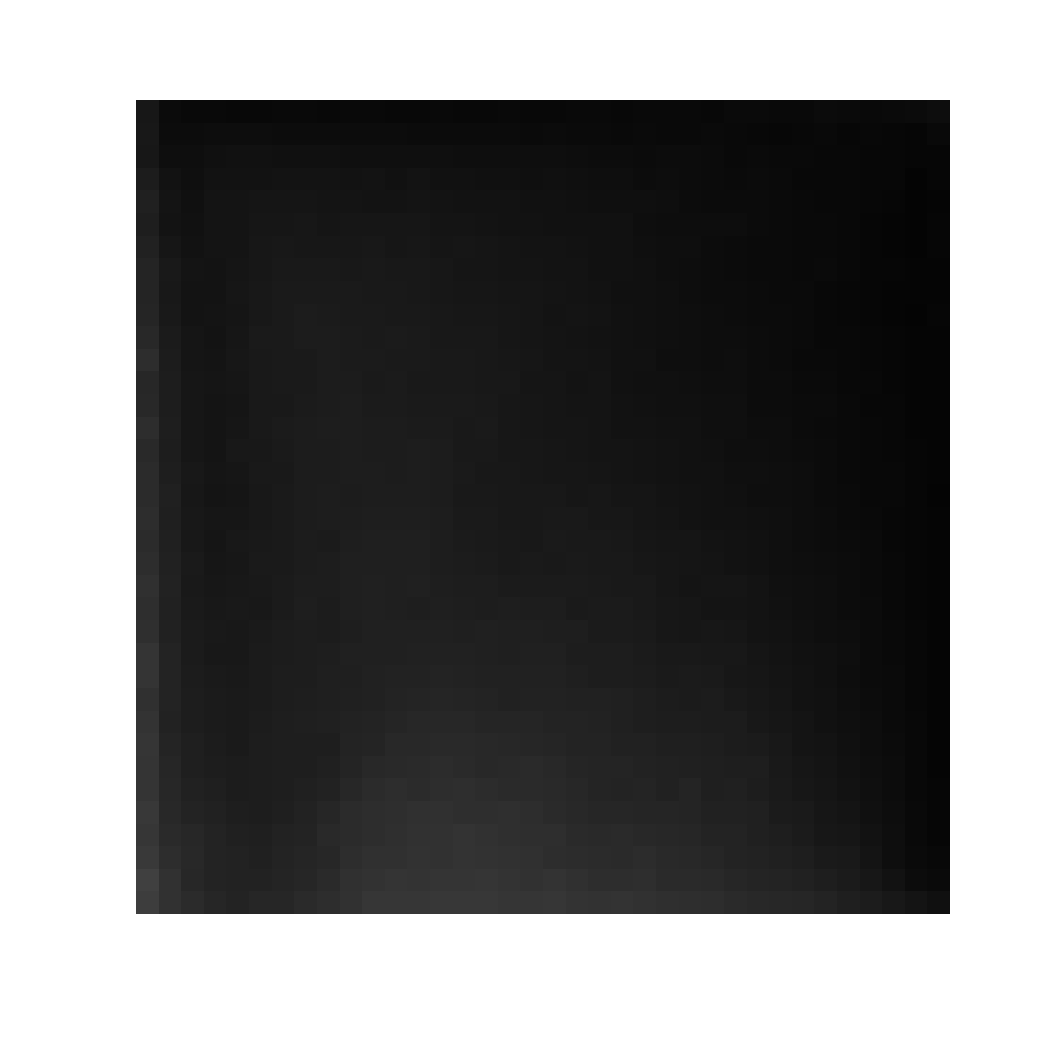}
        & \includegraphics[width=.277\textwidth]{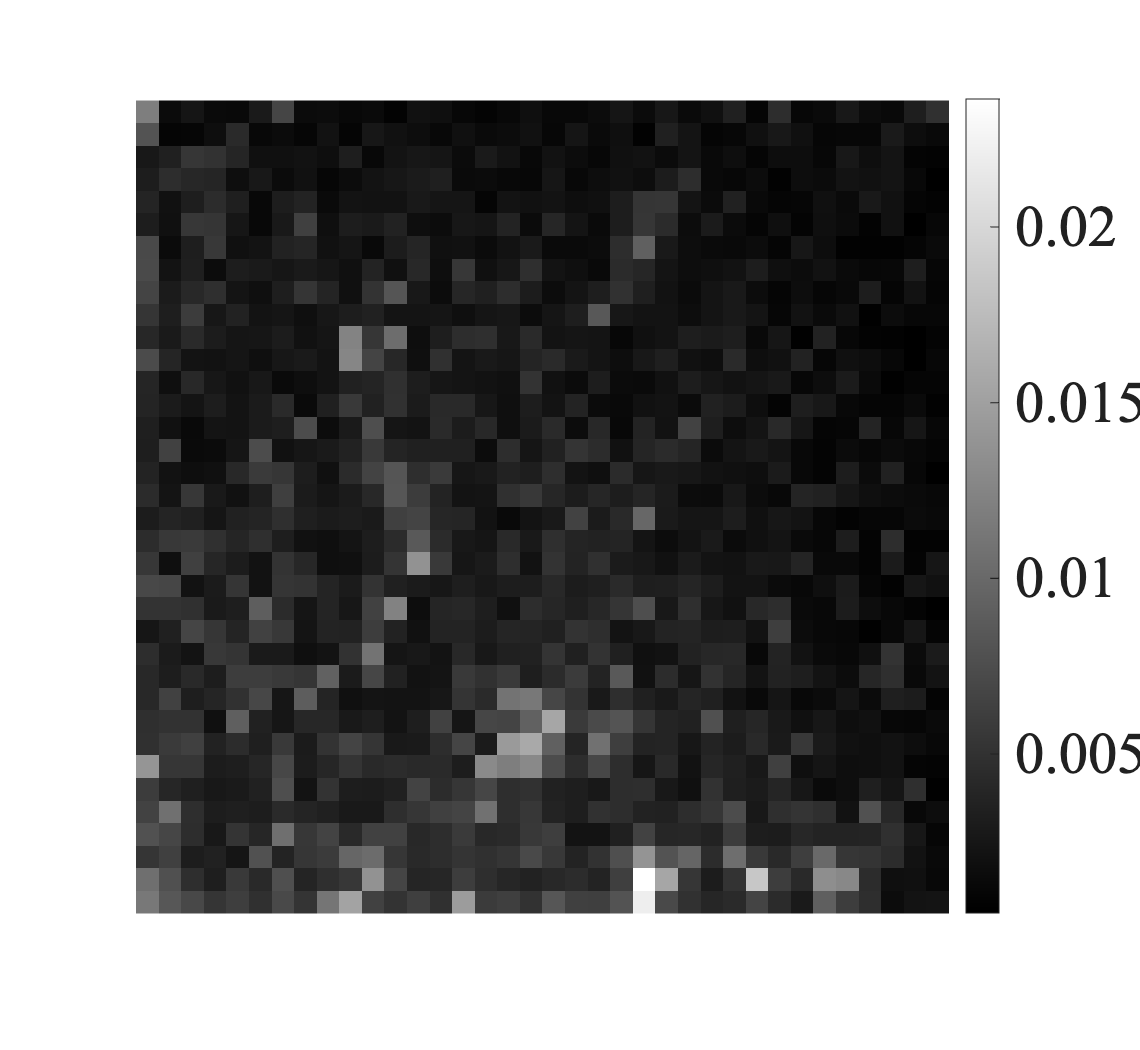}
        \\
        \raisebox{1cm}{$\lambda$} 
        & \includegraphics[width=.25\textwidth]{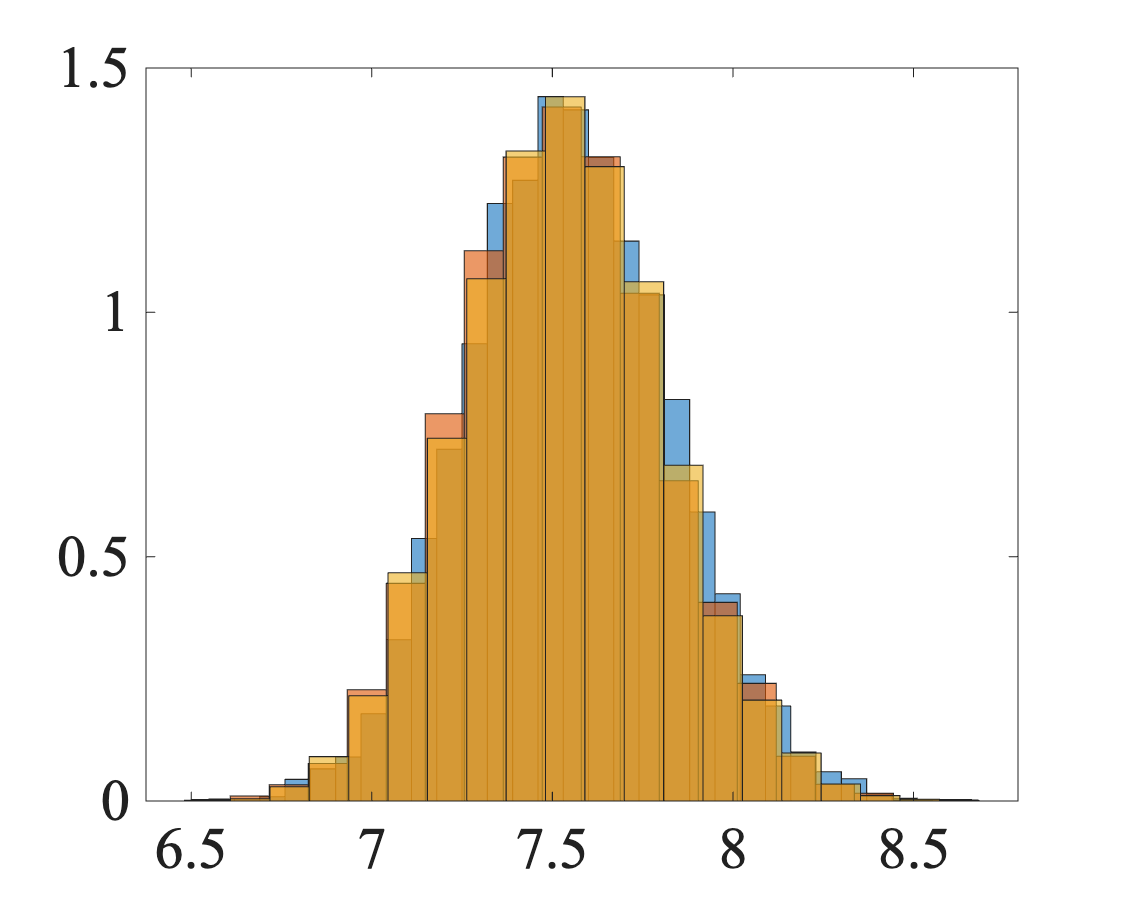}
        & \includegraphics[width=.25\textwidth]{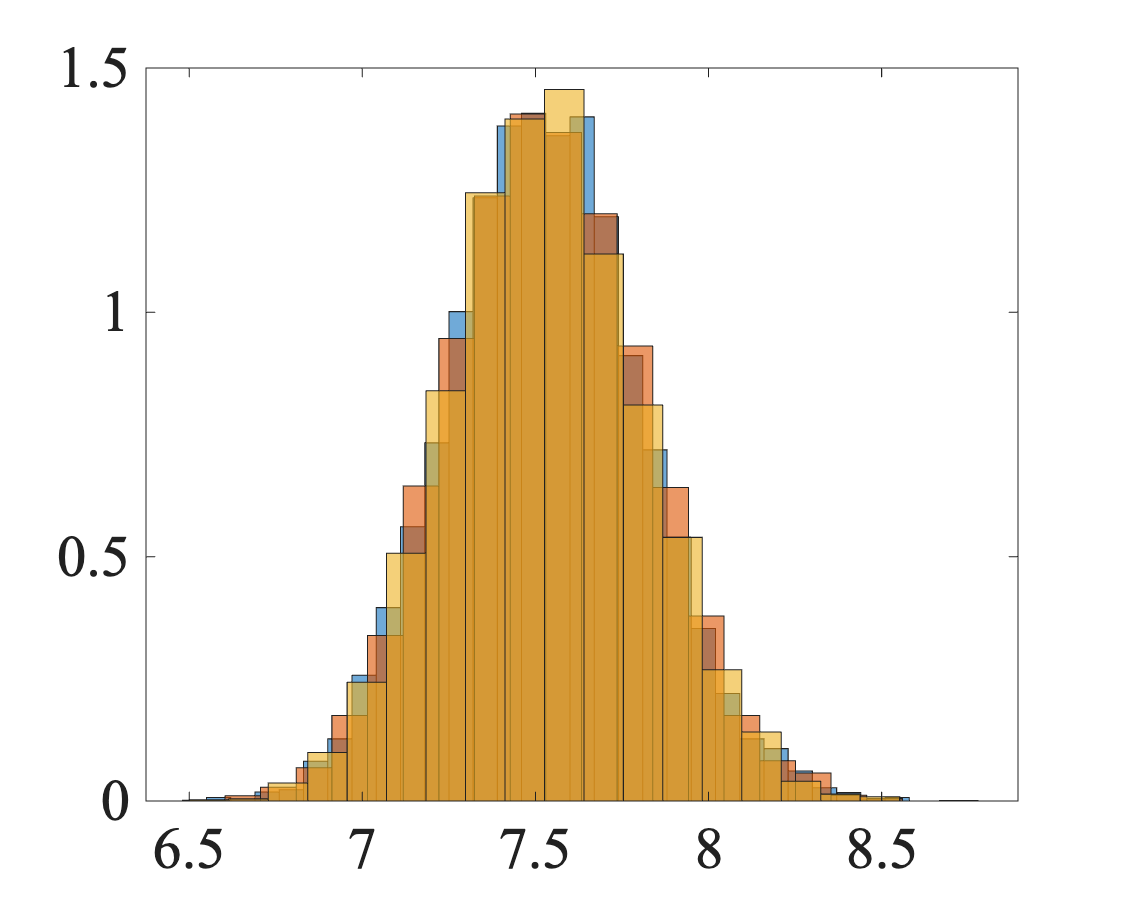}
        & \includegraphics[width=.25\textwidth]{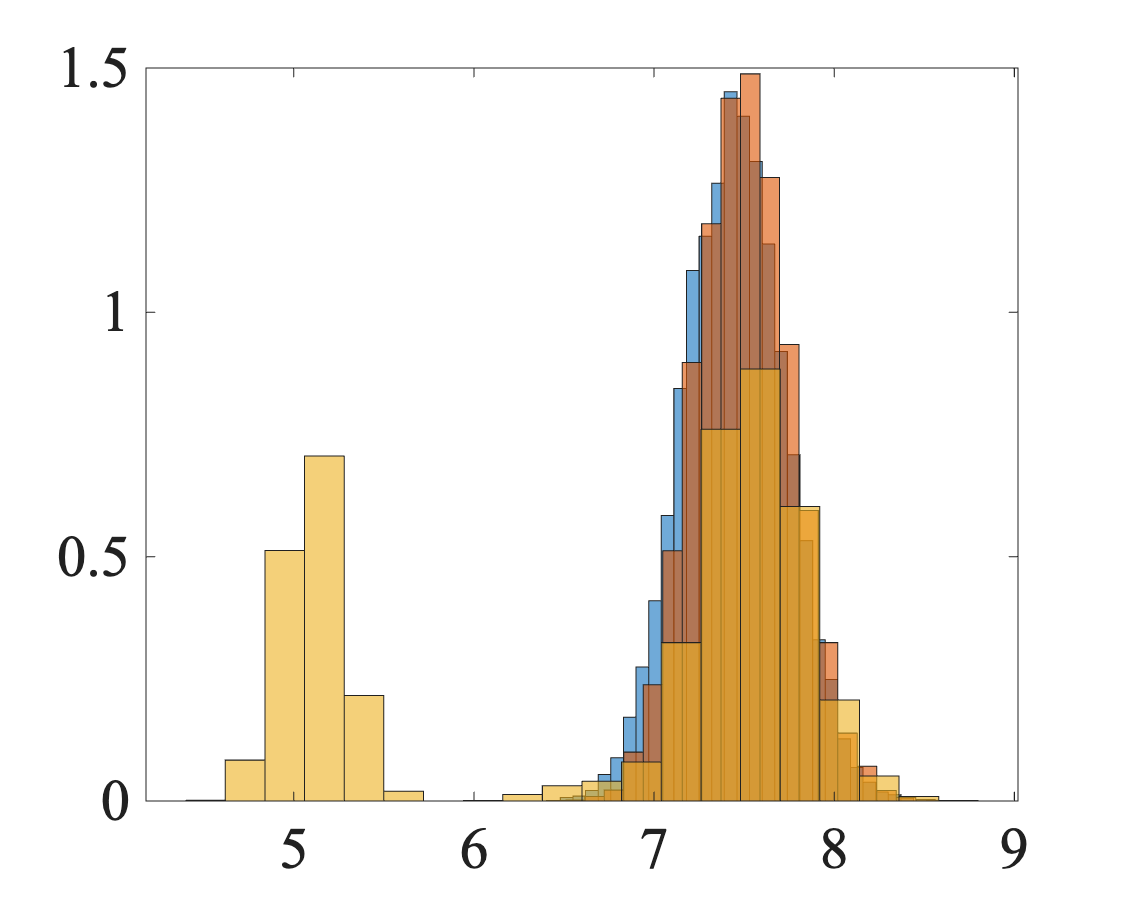}
        \\
        \raisebox{1cm}{$\delta$}
        & \includegraphics[width=.25\textwidth]{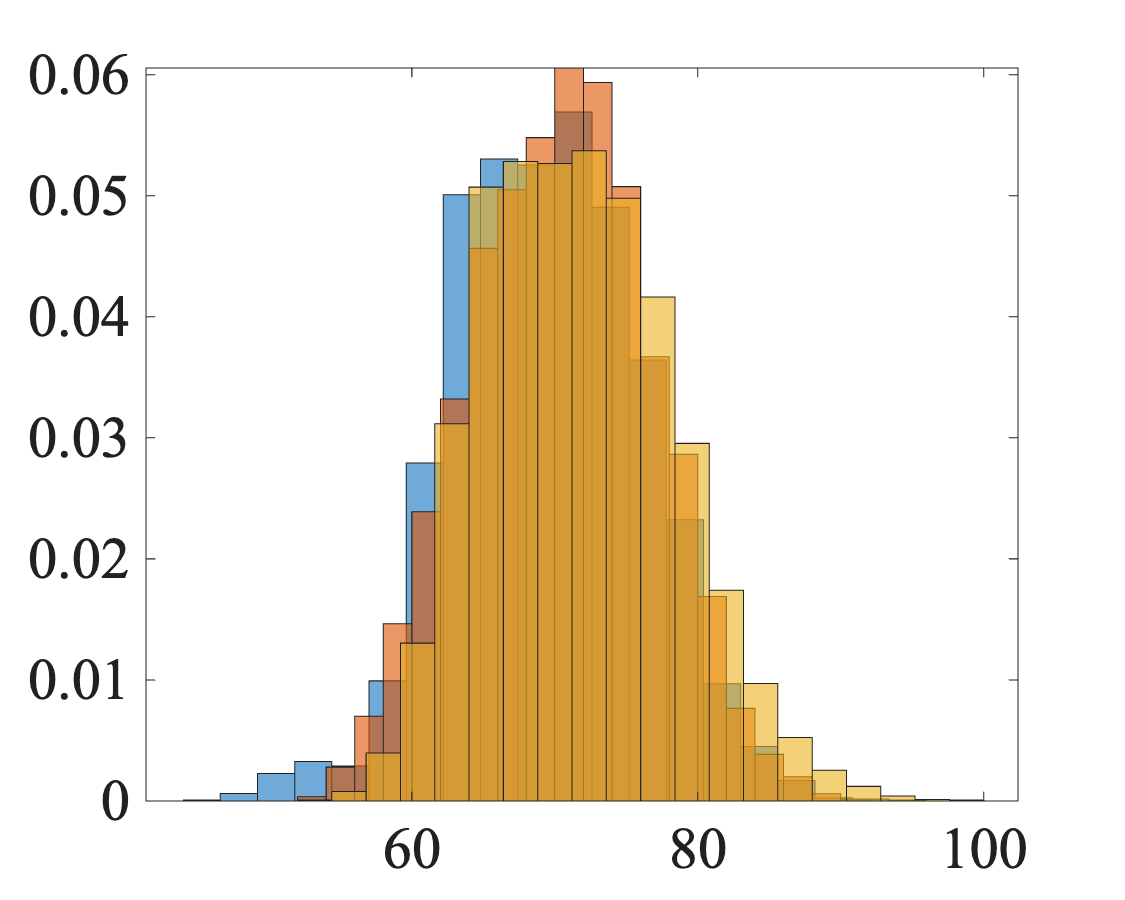}
        & \includegraphics[width=.25\textwidth]{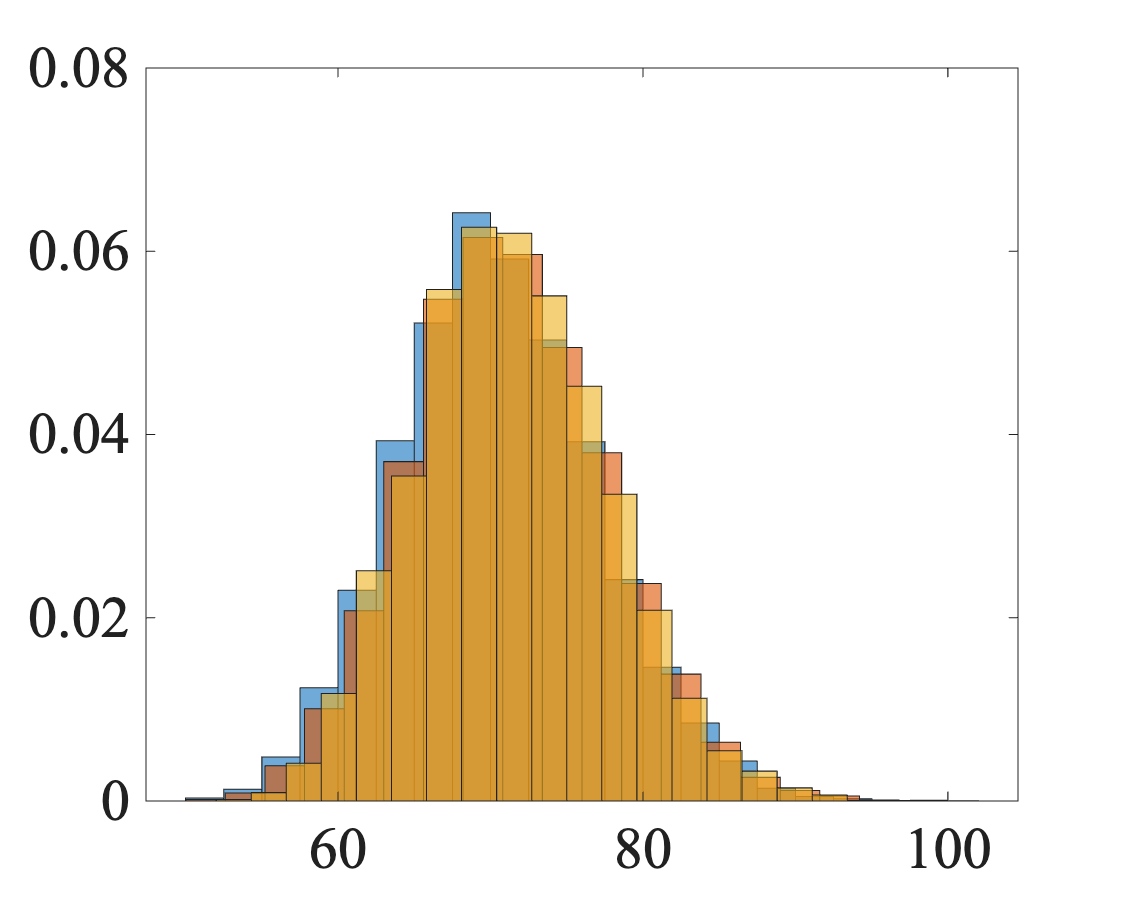}
        & \includegraphics[width=.25\textwidth]{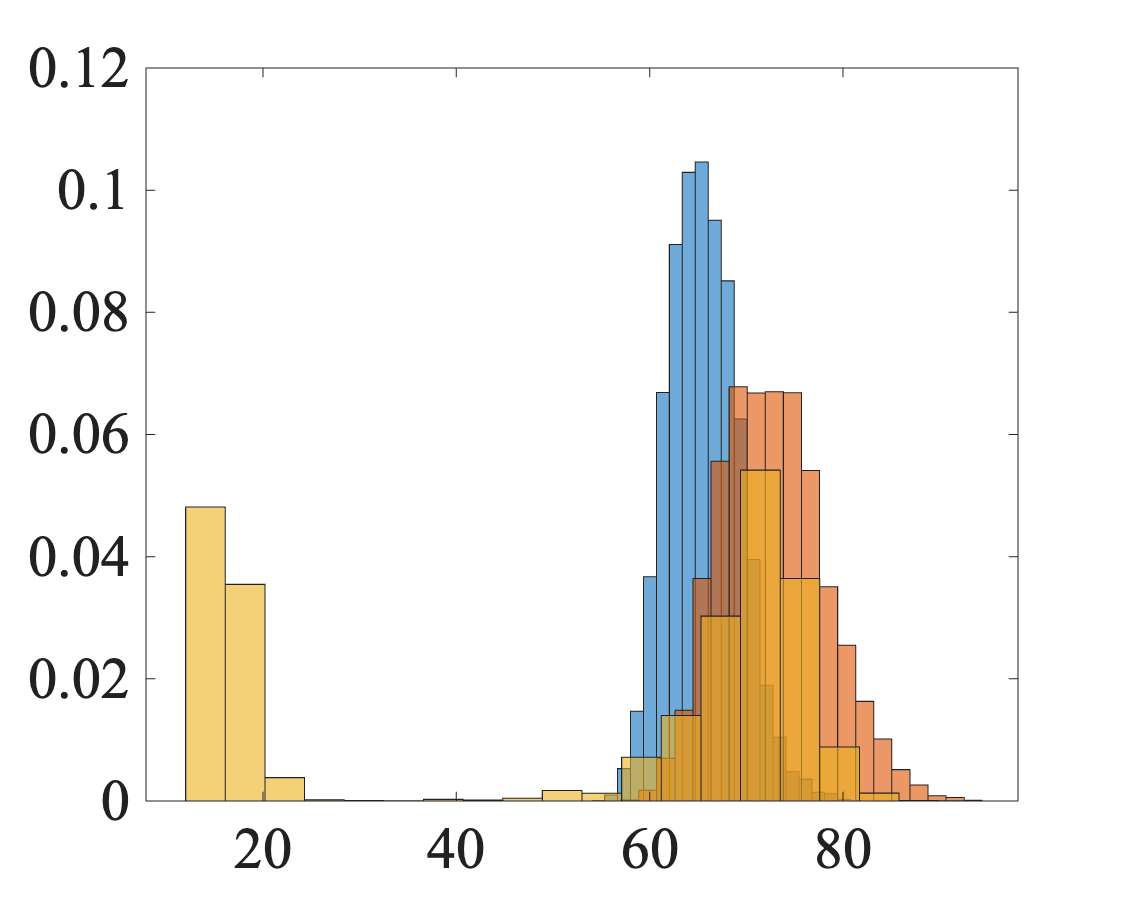}
    \end{tabular}
    \caption{The results for three chains, each with $T=20,000$ samples, on the \texttt{PRseismic} test problem using genGK, tSVD, and rSVD for the approximate covariance matrix $\widehat{\bfGamma}$ with rank $k=500$ in the \jmc{proposal} distribution $g_1(\bfx)$. The mean and variance represent the mean and variance of all accepted samples after burn-in. The $\lambda$ and $\delta$ distributions are normalized histograms containing all draws from $\pi_\lambda$ and $\pi_\delta$ after 50\% burn-in.}
    \label{fig:seismic_500}
\end{figure}

\begin{figure}[bthp]
    \centering
    \begin{tabular}{c c c c}
         & \textbf{genGK} & \textbf{tSVD} & \textbf{rSVD}  \\
        \raisebox{1.7cm}{\textbf{Mean}}
        & \includegraphics[width=.25\textwidth]{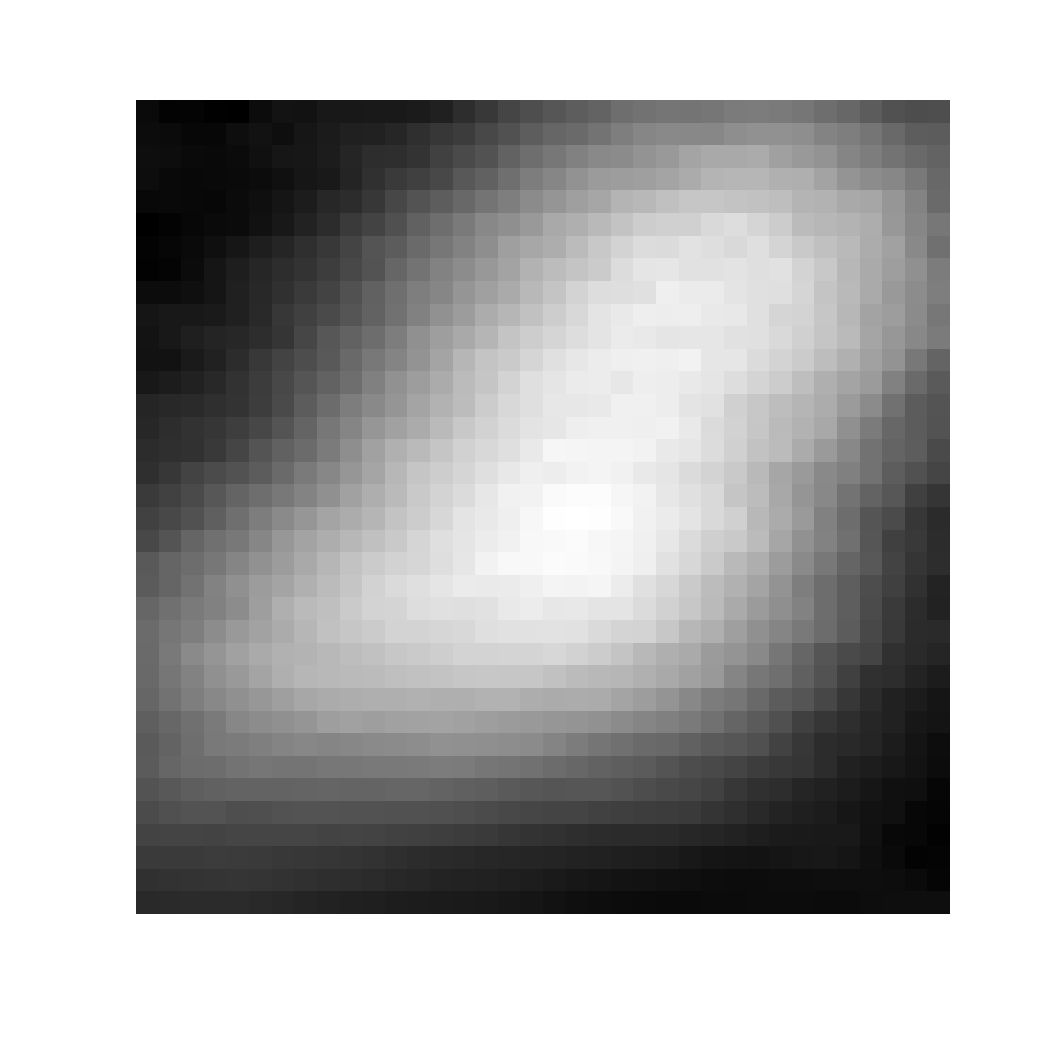}
        & \includegraphics[width=.25\textwidth]{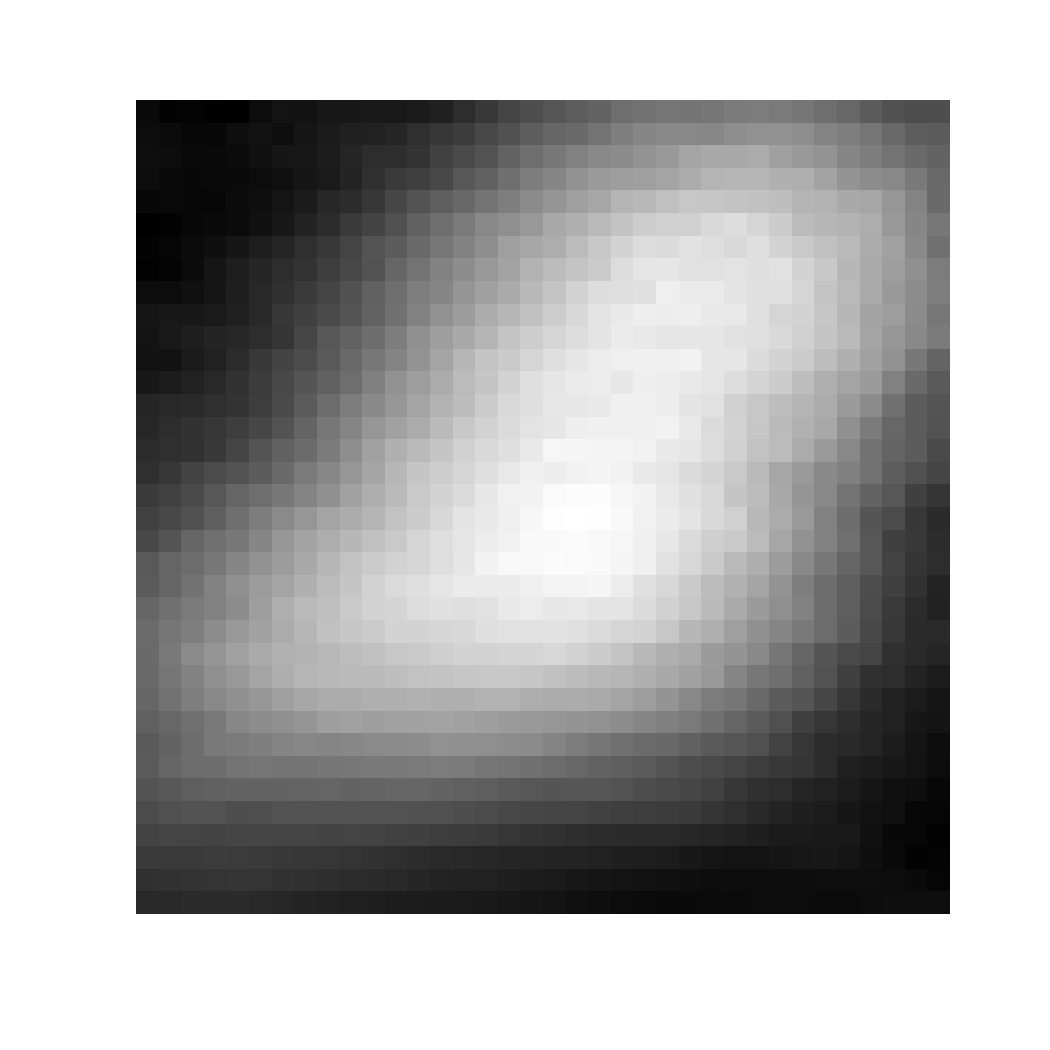}
        & \includegraphics[width=.277\textwidth]{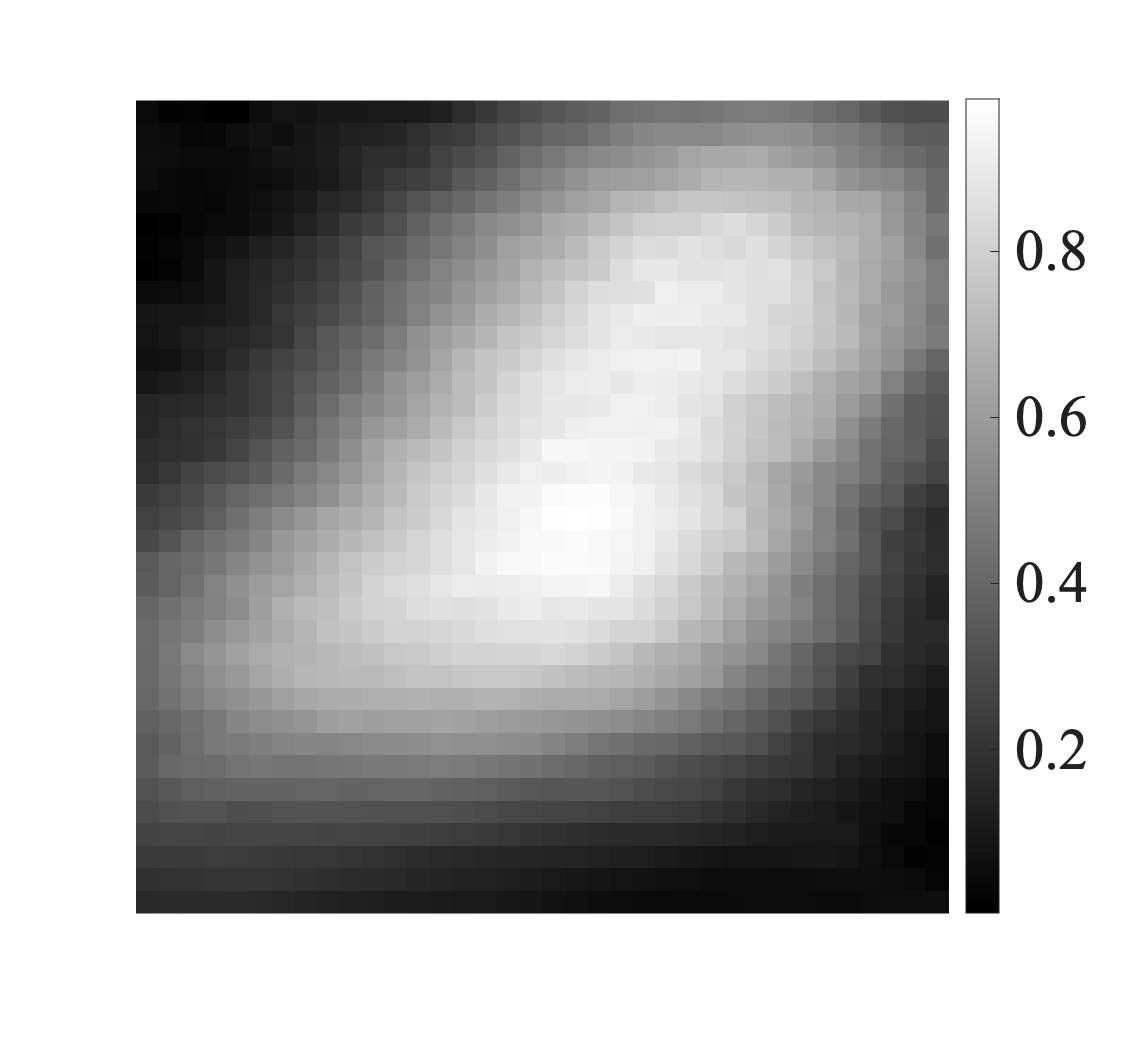}
        \\
        \raisebox{1.7cm}{\textbf{Variance}}
        & \includegraphics[width=.25\textwidth]{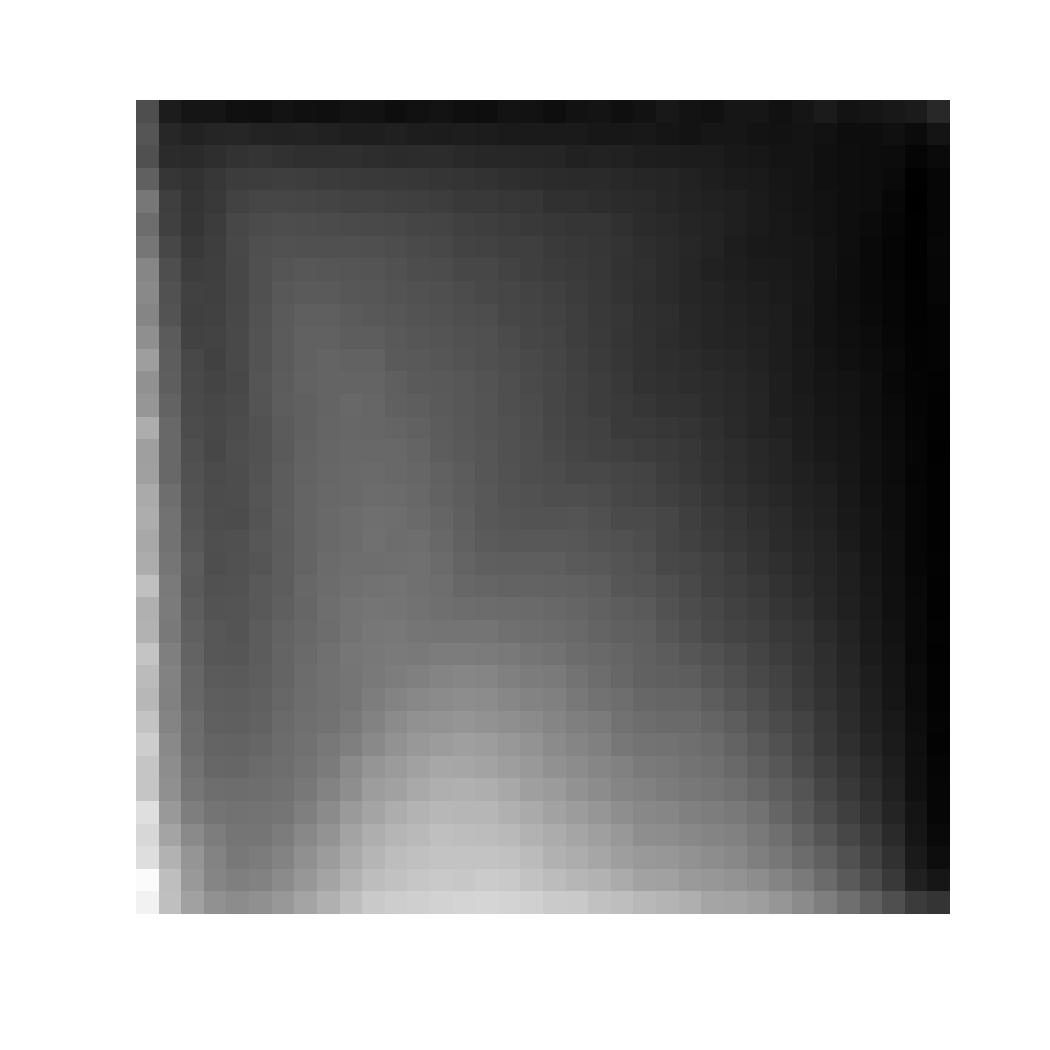}
        & \includegraphics[width=.25\textwidth]{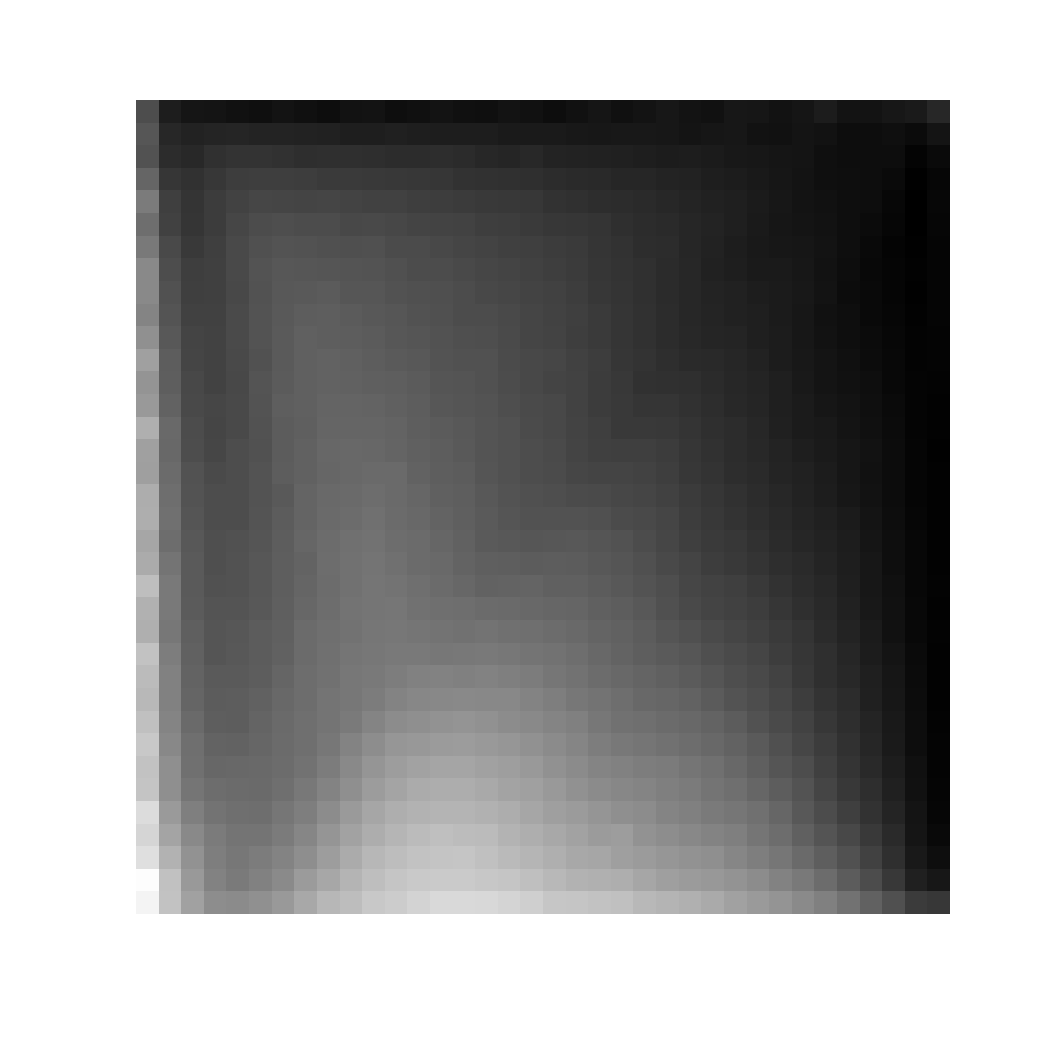}
        & \includegraphics[width=.277\textwidth]{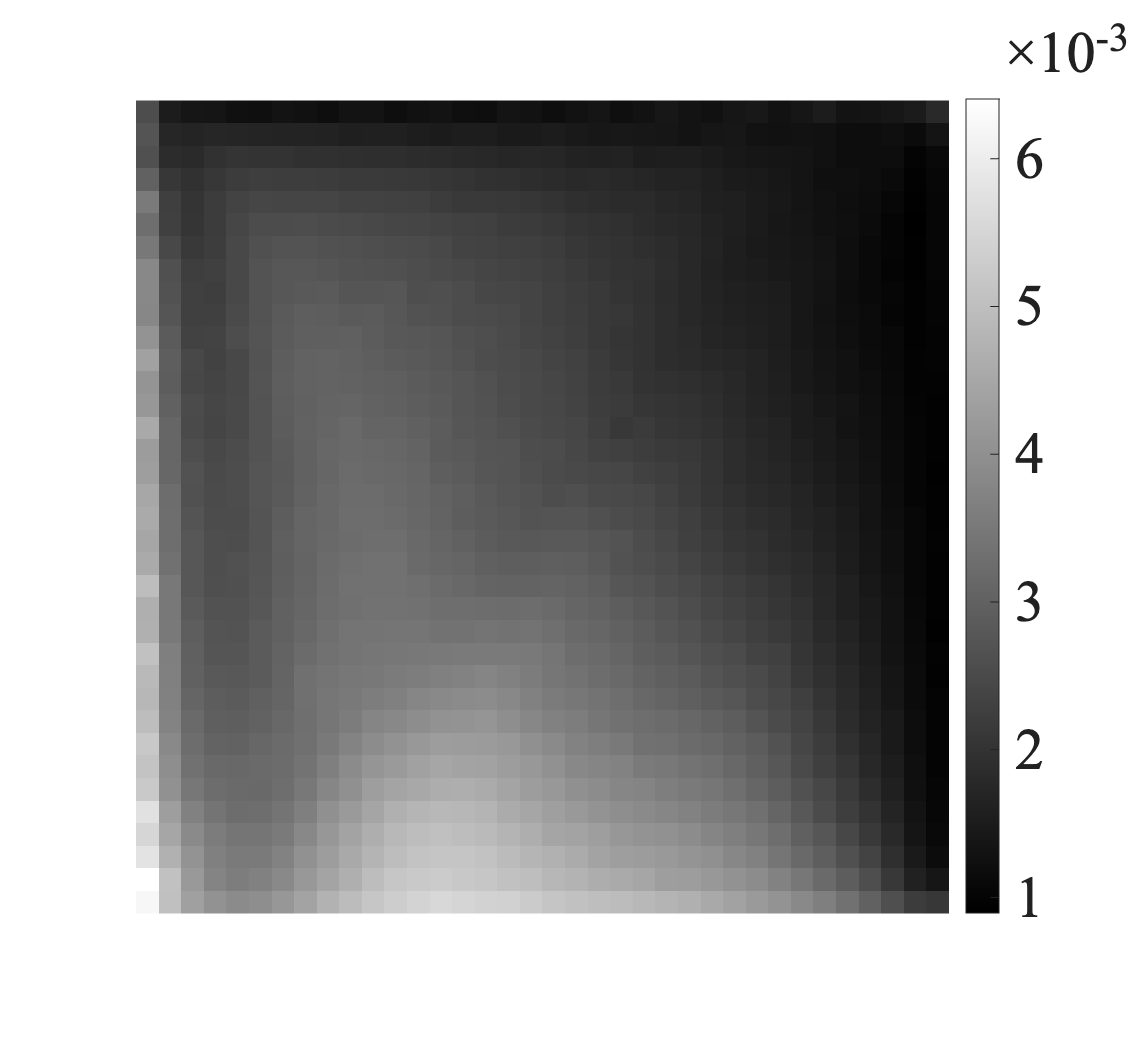}
        \\
        \raisebox{1cm}{$\lambda$} 
        & \includegraphics[width=.25\textwidth]{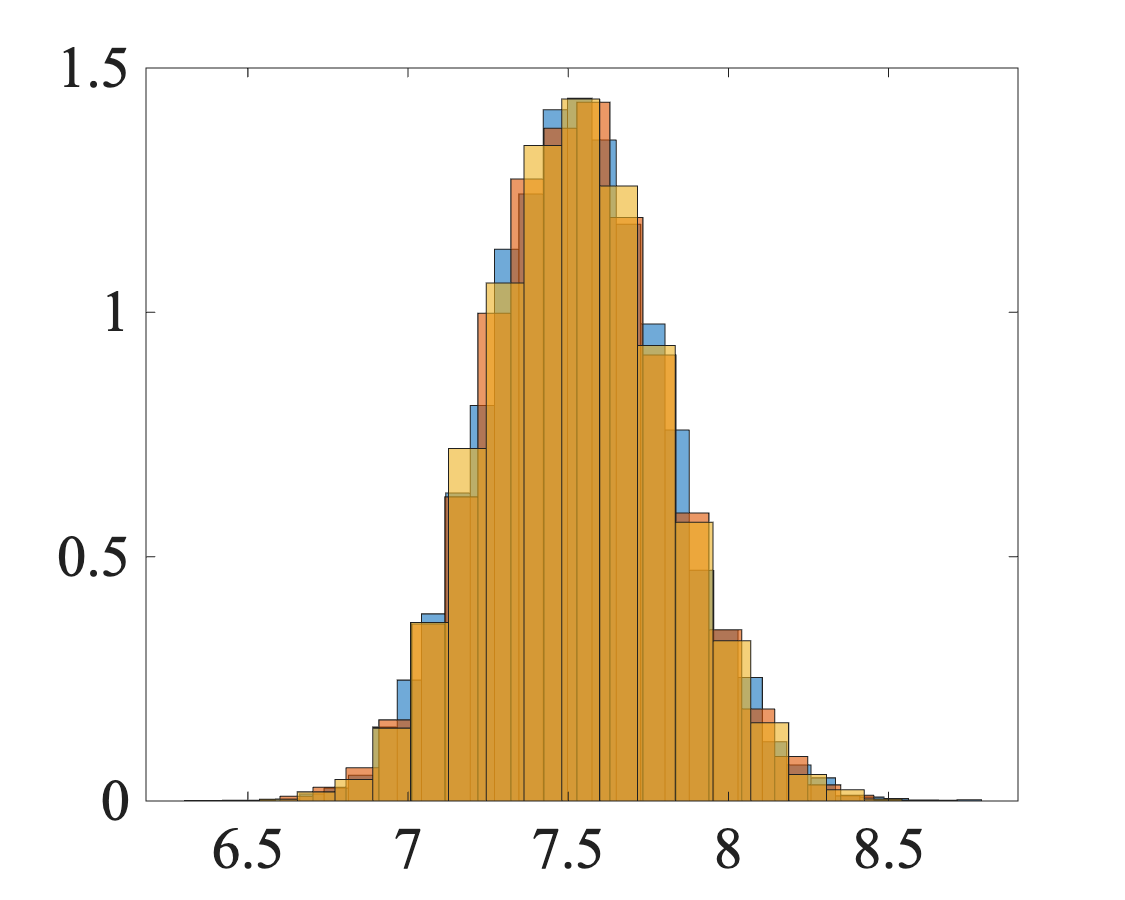}
        & \includegraphics[width=.25\textwidth]{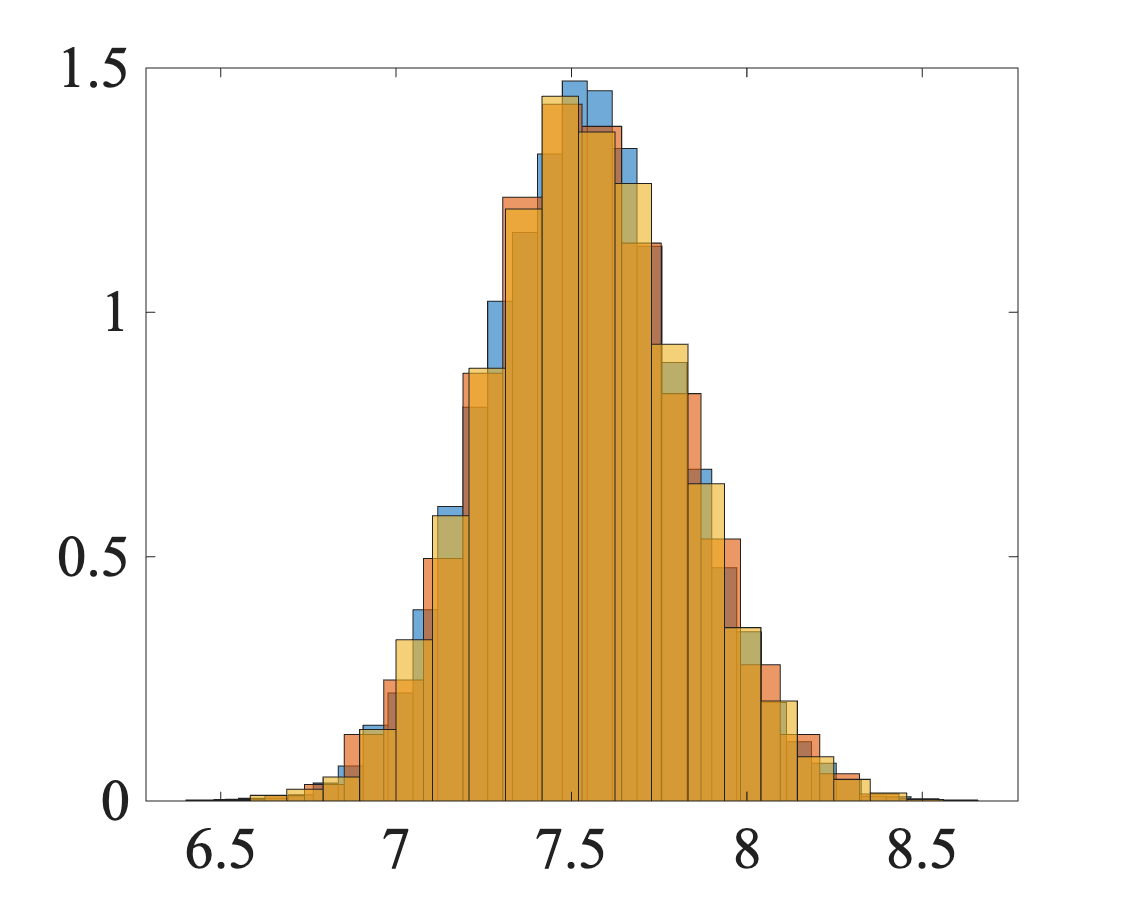}
        & \includegraphics[width=.25\textwidth]{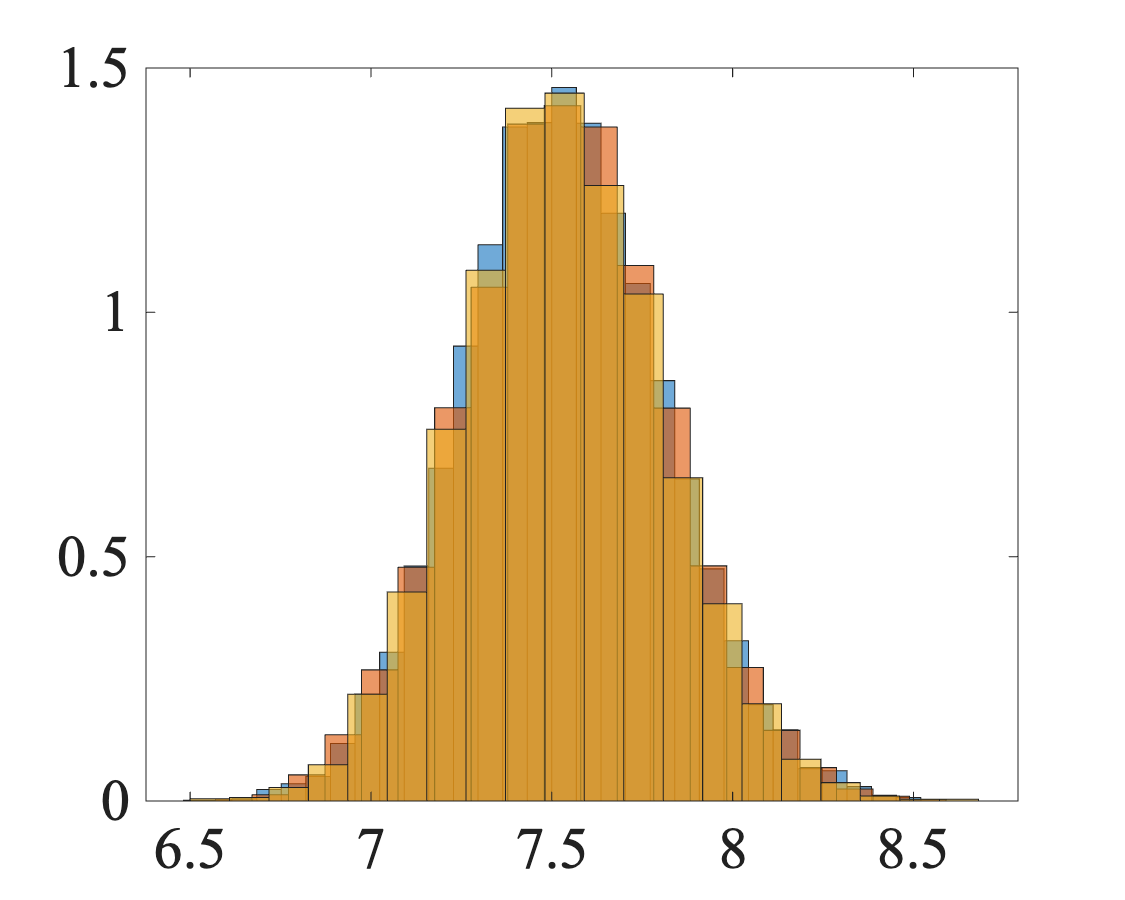}
        \\
        \raisebox{1cm}{$\delta$}
        & \includegraphics[width=.25\textwidth]{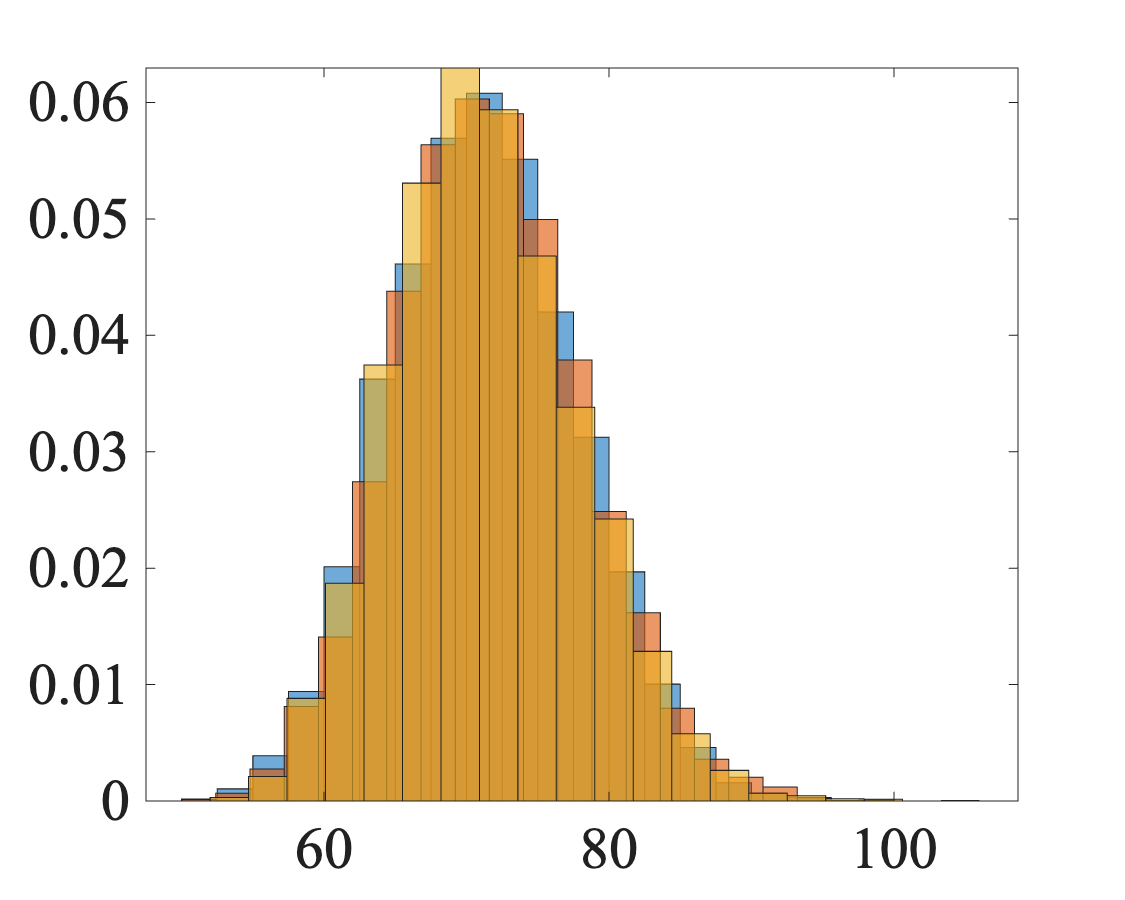}
        & \includegraphics[width=.25\textwidth]{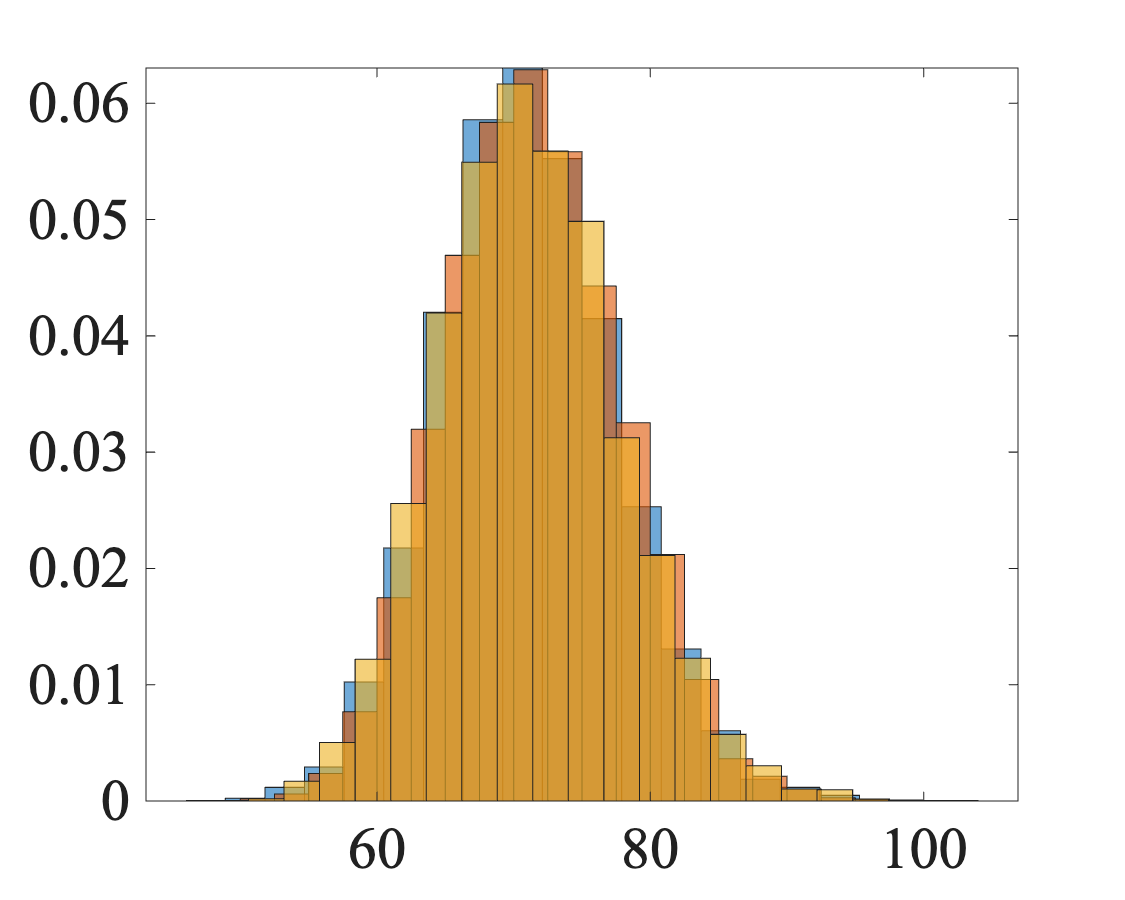}
        & \includegraphics[width=.25\textwidth]{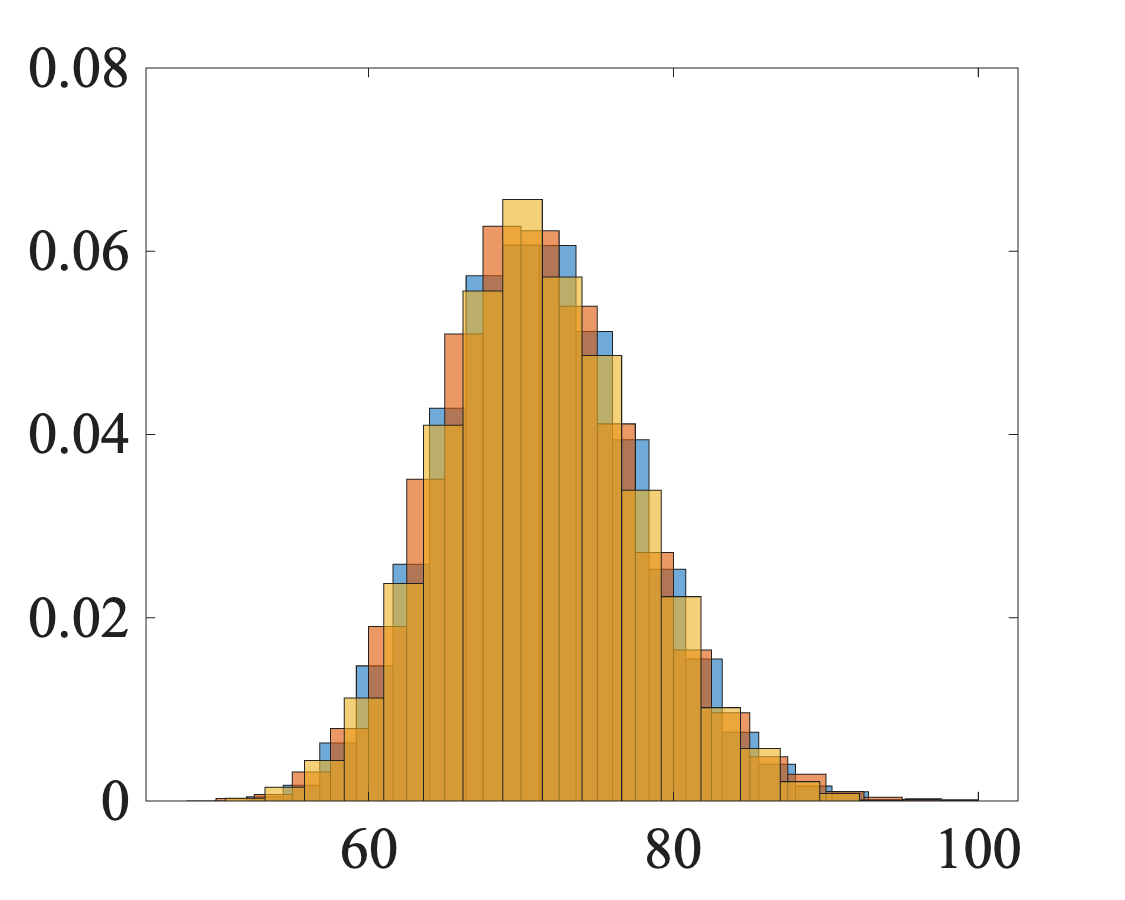}
    \end{tabular}
    \caption{For the \texttt{PRseismic} test problem, we provide results for three chains, each with $T=20,000$ samples, using genGK, tSVD, and rSVD for the approximate covariance matrix $\widehat{\bfGamma}$ with rank $k=1,000$ in the proposed distribution $g_1(\bfx)$. We provide the mean and variance of all accepted samples after 50\% burn-in, along with the normalized histograms from the $\lambda$ and $\delta$ distributions.}
    \label{fig:seismic_1000}
\end{figure}

\begin{figure}[bthp]
    \centering
    \begin{tabular}{cc}
        $k=500$ & $k=1,000$  \\
        \includegraphics[width=.37\textwidth]{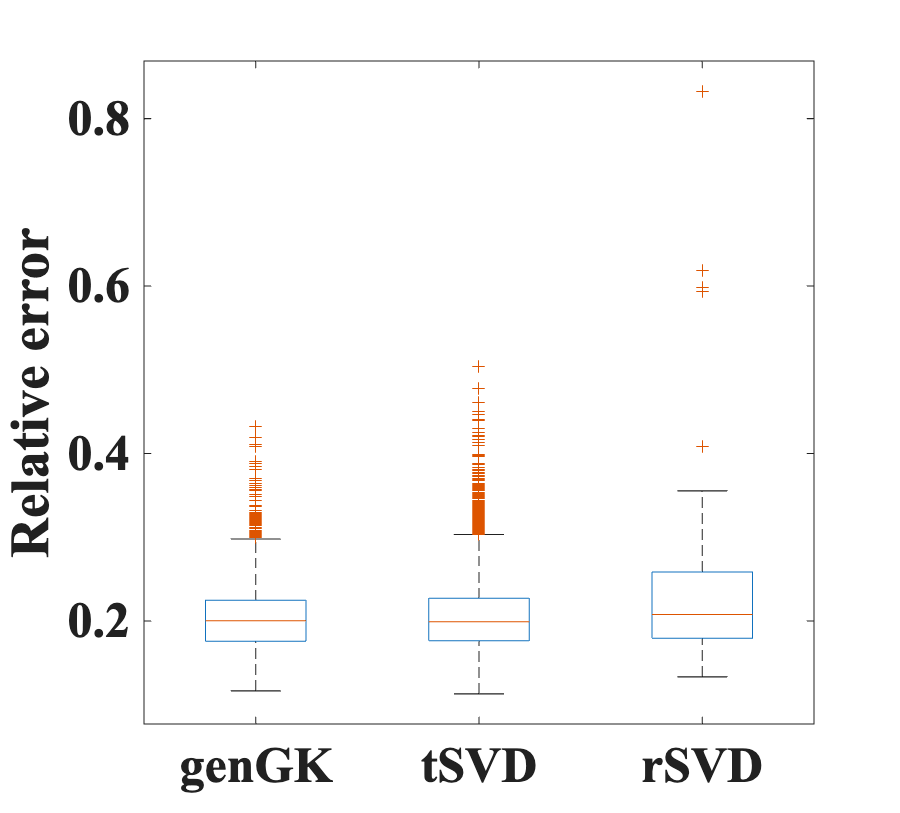} 
        & \includegraphics[width=.37\textwidth]{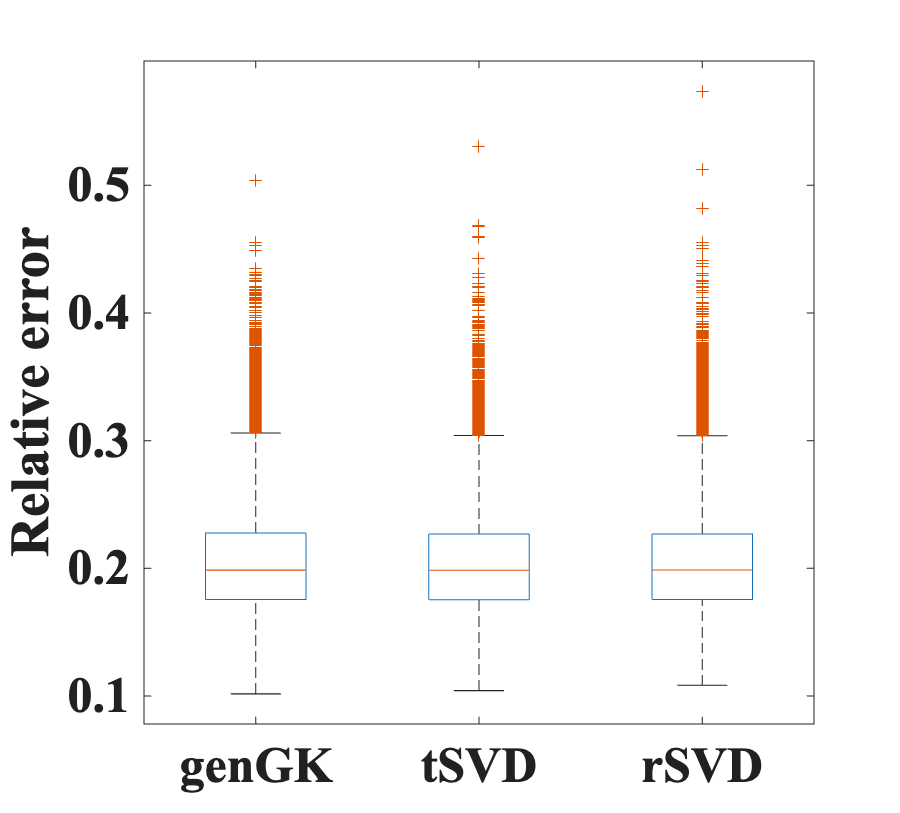}
    \end{tabular}
    \caption{Box-and-whisker plots providing the relative reconstruction errors for all accepted samples for the \texttt{PRseismic} test problem using genGK, tSVD, and rSVD for the approximate covariance matrix $\widehat{\bfGamma}$ with rank $k=500$ (left) and $k=1,000$ (right).}
    \label{fig:boxandwhisker}
\end{figure}

We compute $T=20,000$ samples using the Metropolis-Hastings within Gibbs method with proposal sampling using the genGK approximation (\Cref{alg:gibbs_genGK}) with $k=500$ and $k=1,000$.  
In \Cref{fig:seismic_500} and \Cref{fig:seismic_1000} we compare the mean and variance images for reconstructions, and the normalized distributions of computed hyperparameters $\lambda$ and $\delta$ using the genGK approximation to low-rank proposal sampling given by truncated SVD (tSVD) and rSVD approximations. Details for using these alternative methods can be found in the Appendix, \Cref{sec:svdapprox}. The relative reconstruction error norms between the true solution and each of the accepted samples for genGK, tSVD, and rSVD are provided in box-and-whisker plots for $k=500$ and $k=1,000$ in \Cref{fig:boxandwhisker}. In \Cref{tab:seismic_diag} we provide the acceptance rate, runtime, and diagnostics of the $\lambda$ and $\delta$ chains for each of the low-rank proposal sampling methods.  We observe that as the rank of the approximation $\widehat\bfGamma_{\rm cond}$ increases, $g_1(\bfx)$ becomes closer to the target distribution $h(\bfx)$ and the acceptance rate approaches $100\%$. The acceptance rates found using genGK are comparable to those using the rSVD method in \cite{brown2018low}. As the rank increases, we can also see that 
$\widehat{R}$ is less than 1.01.  For both considered ranks, the Geweke test indicates that the chains are in equilibrium. 
The $\lambda$ chains have $95\%$ confidence intervals of $[7.01,\ 8.03]$ and $[6.99,\ 8.06]$ for $k=500$ and $1,000$ respectively. 

\begin{figure}[bthp]
\centering
\begin{tabular}{ccc}
\textbf{$k$} & \textbf{$\lambda$ Samples} & \textbf{\jmc{ACF}} \\  
\raisebox{1.8cm}{500} & \raisebox{.05cm}{\includegraphics[width=.37\textwidth]{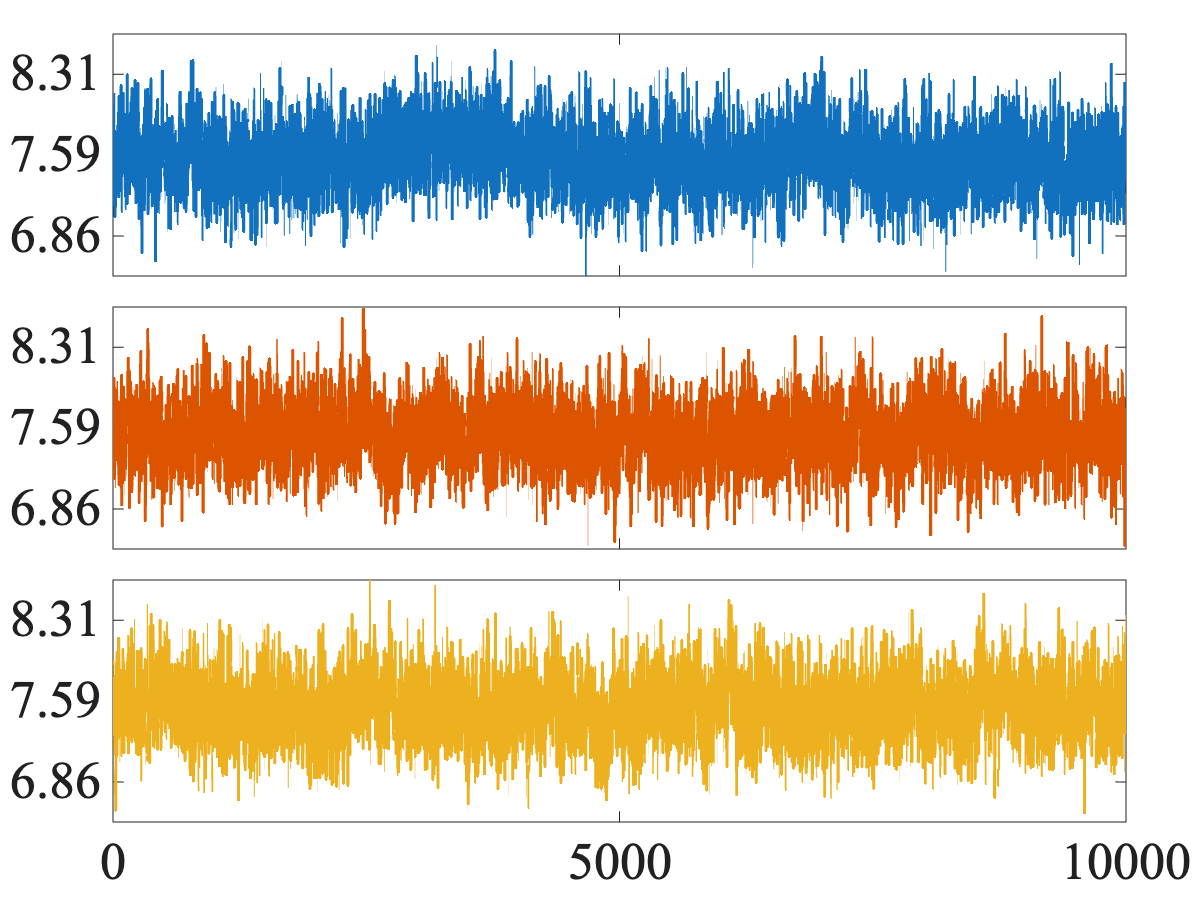}} & {\includegraphics[width=.285\textwidth]{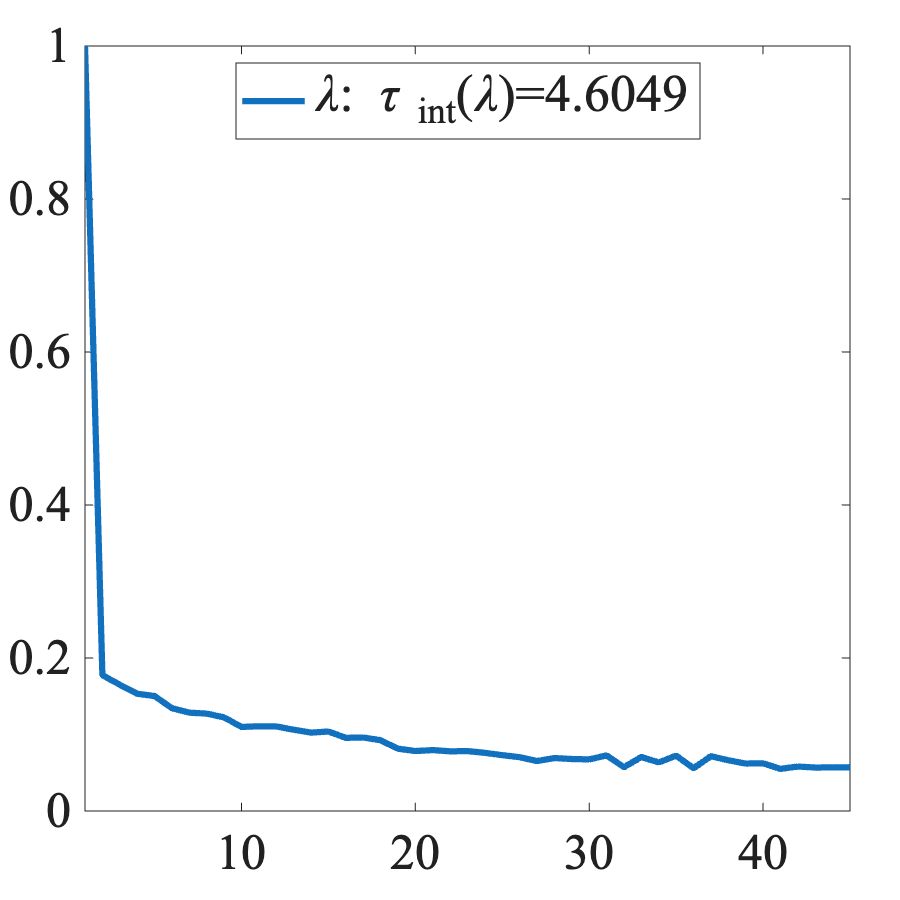}} 
\\
\raisebox{1.8cm}{1,000} &  \raisebox{.05cm}{\includegraphics[width=.37\textwidth]{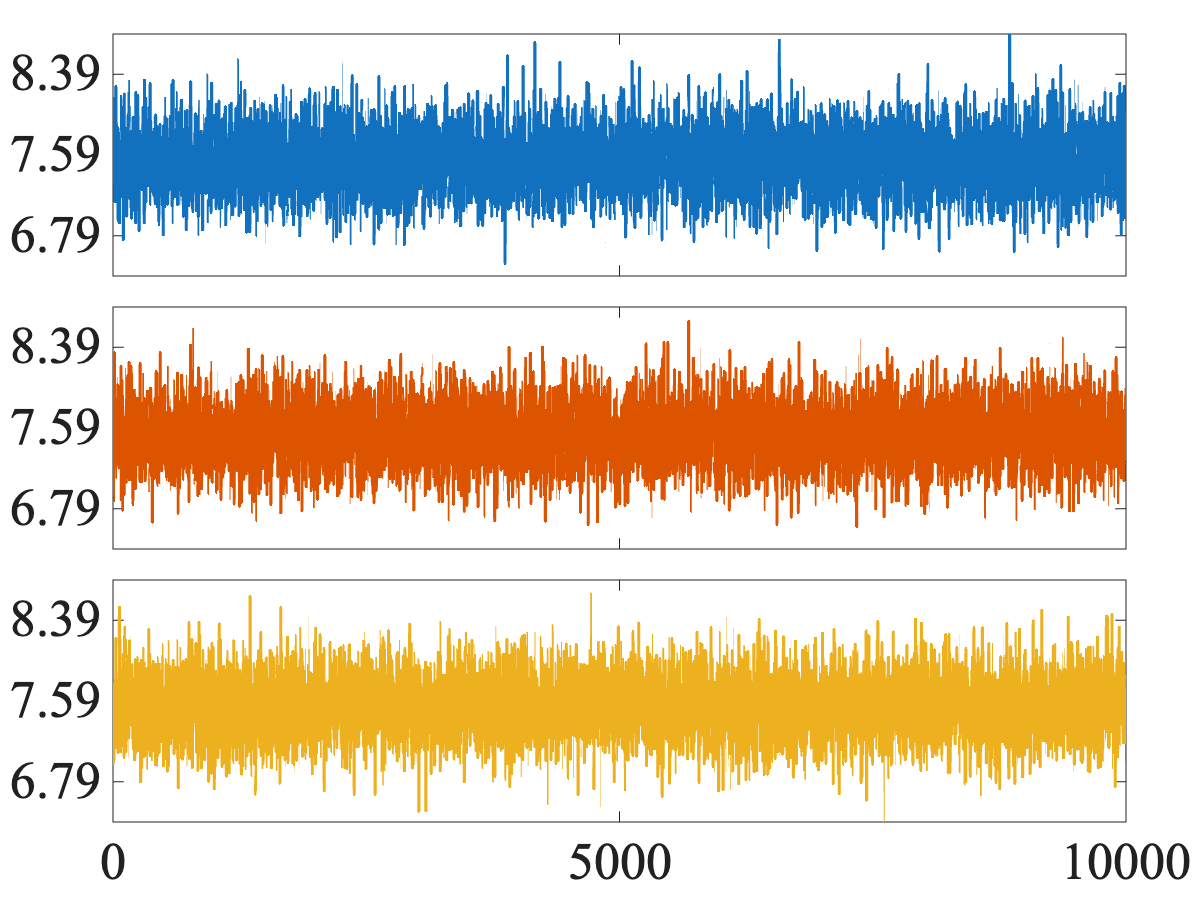}} & {\includegraphics[width=.285\textwidth]{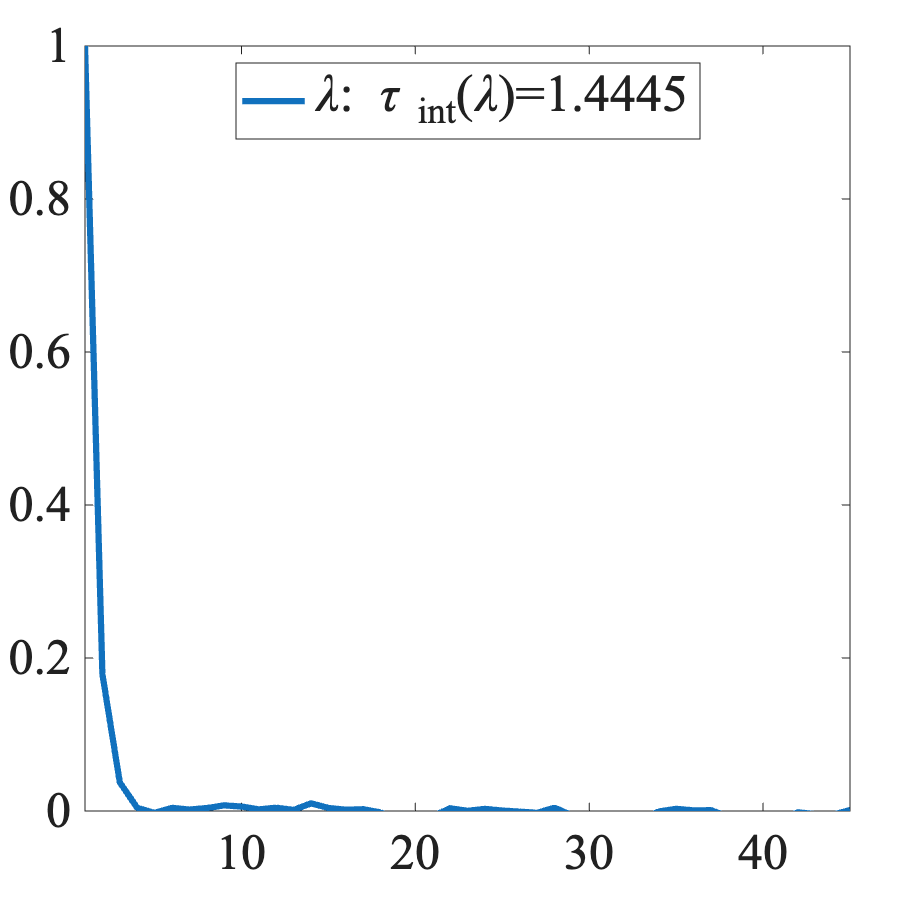}} 
\\
\end{tabular}
\caption{We provide three $\lambda$ chains for $T=20,000$ samples using \Cref{alg:gibbs_genGK} on the \texttt{PRseismic} test problem, along with the \jmc{ACF}. Here, $k$ denotes the rank of $\widehat\bfGamma$. }

\label{tab:lamchain}
\end{figure}

\begin{figure}[bthp]
\centering
\begin{tabular}{ccc}
\textbf{$k$} & \textbf{$\delta$ Samples} & \textbf{\jmc{ACF}} \\  
\raisebox{1.8cm}{500} &  \raisebox{.05cm}{\includegraphics[width=.37\linewidth]{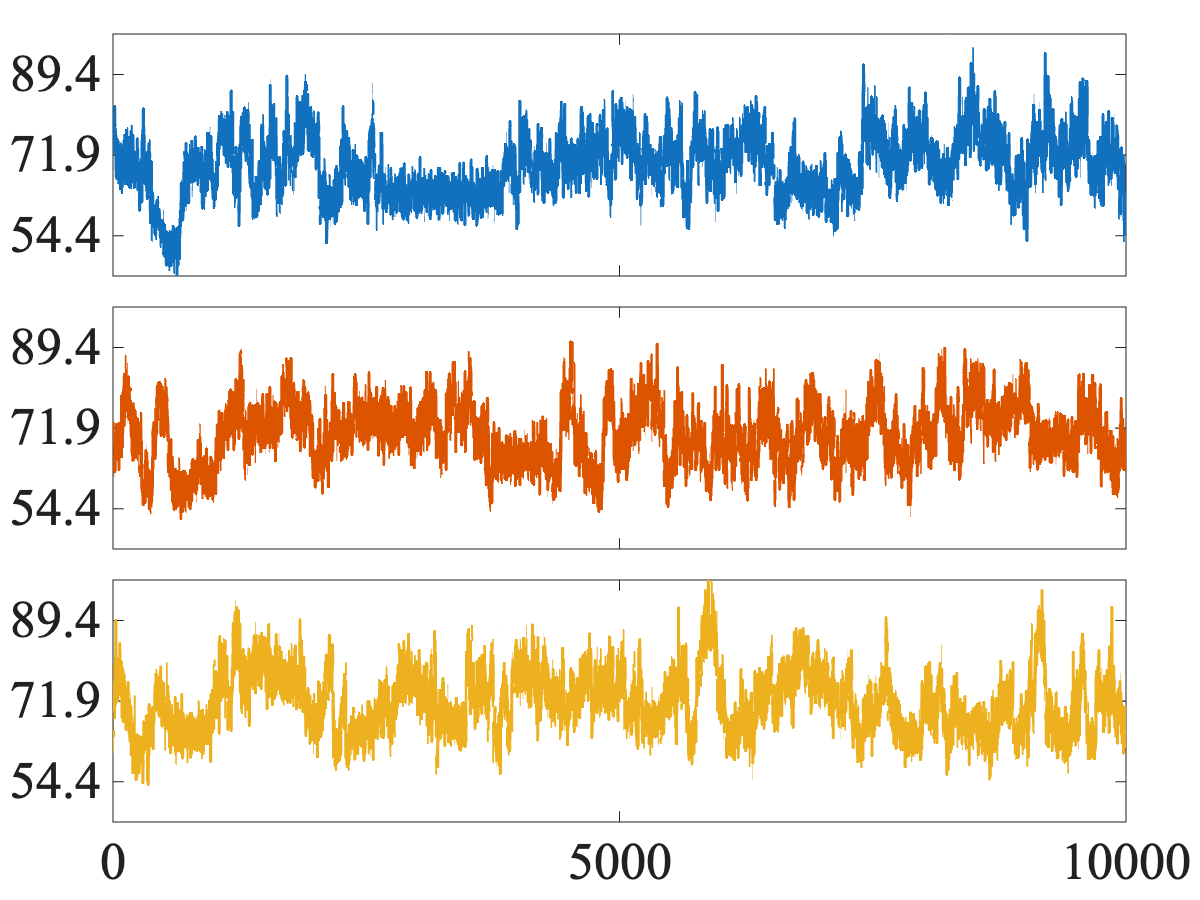}} & {\includegraphics[width=.285\linewidth]{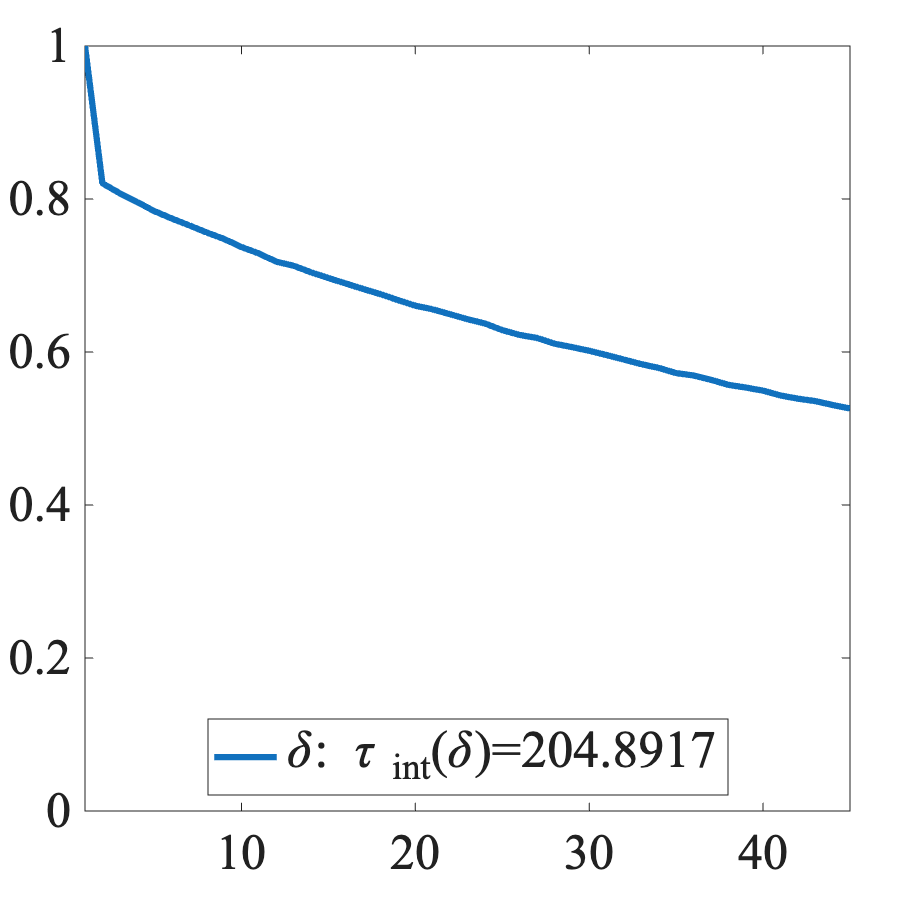}} 
\\
    \raisebox{1.8cm}{1,000} & \raisebox{.05cm}{\includegraphics[width=.37\linewidth]{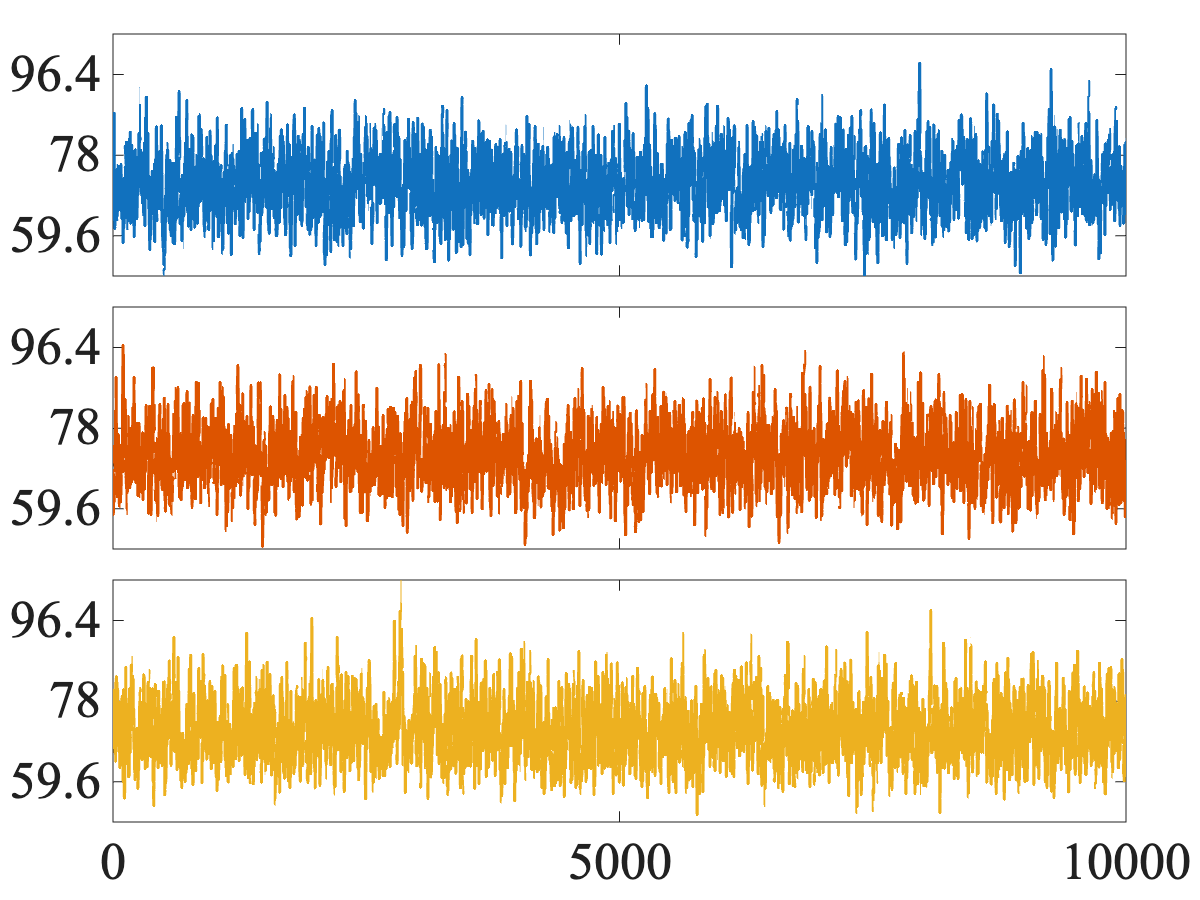}} & {\includegraphics[width=.285\linewidth]{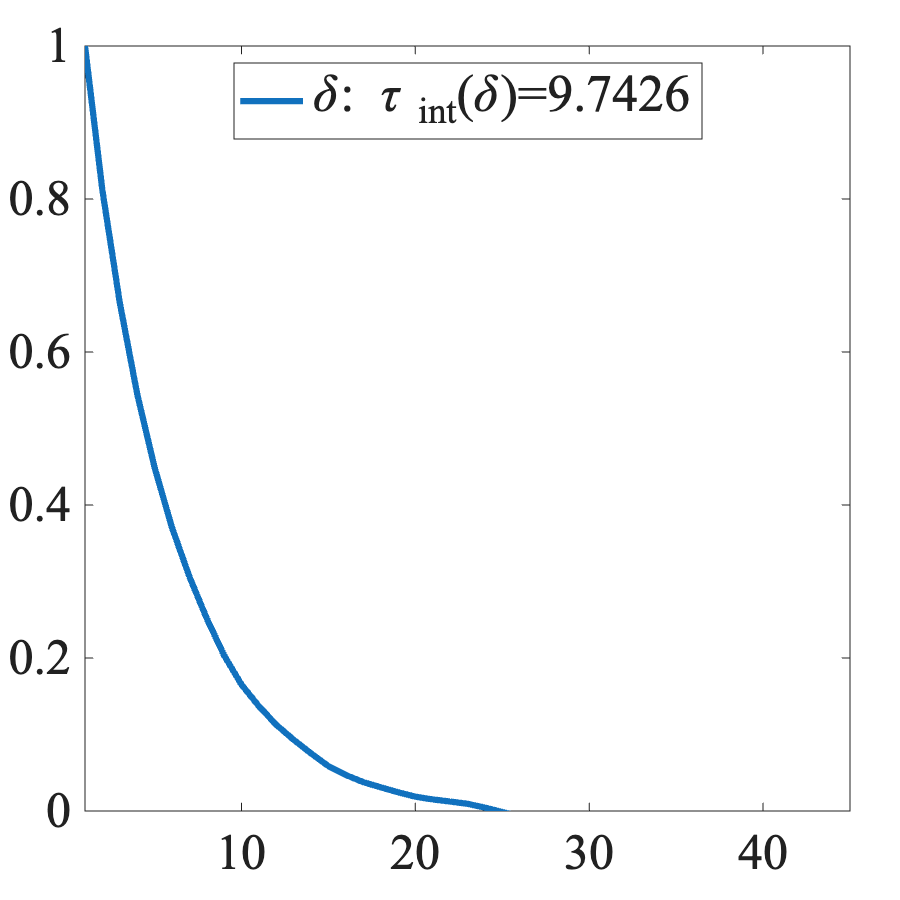}} 
\\
\end{tabular}
\caption{We provide three $\delta$ chains for $T=20,000$ samples using \Cref{alg:gibbs_genGK} on the \texttt{PRseismic} test problem, along with the \jmc{ACF}. Here, $k$ denotes the rank of $\widehat\bfGamma$.}
\label{tab:delchain}
\end{figure}

\begin{figure}[bthp]
    \centering
    \begin{tabular}{cc}
        $k=500$ & $k=1,000$  \\
        \includegraphics[width=.4\textwidth]{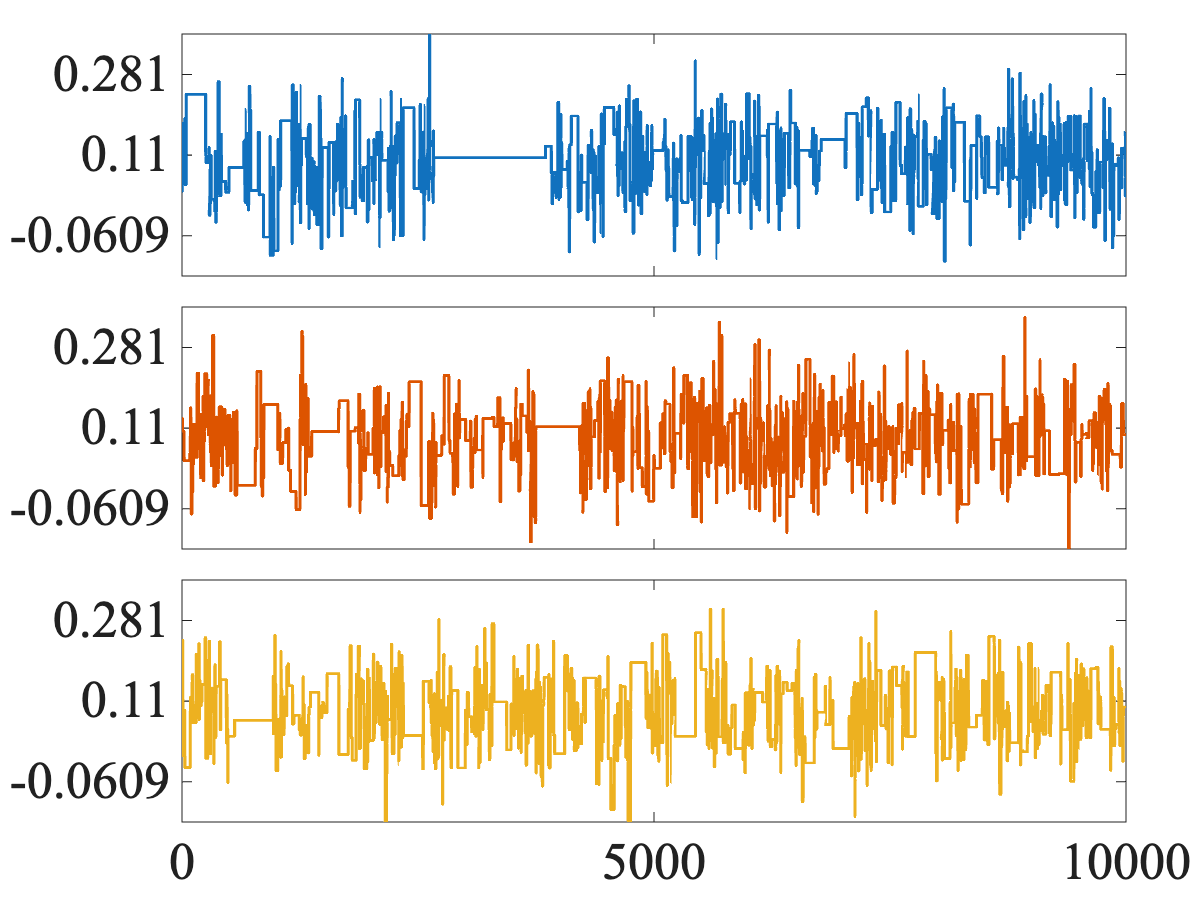} 
        & \includegraphics[width=.4\textwidth]{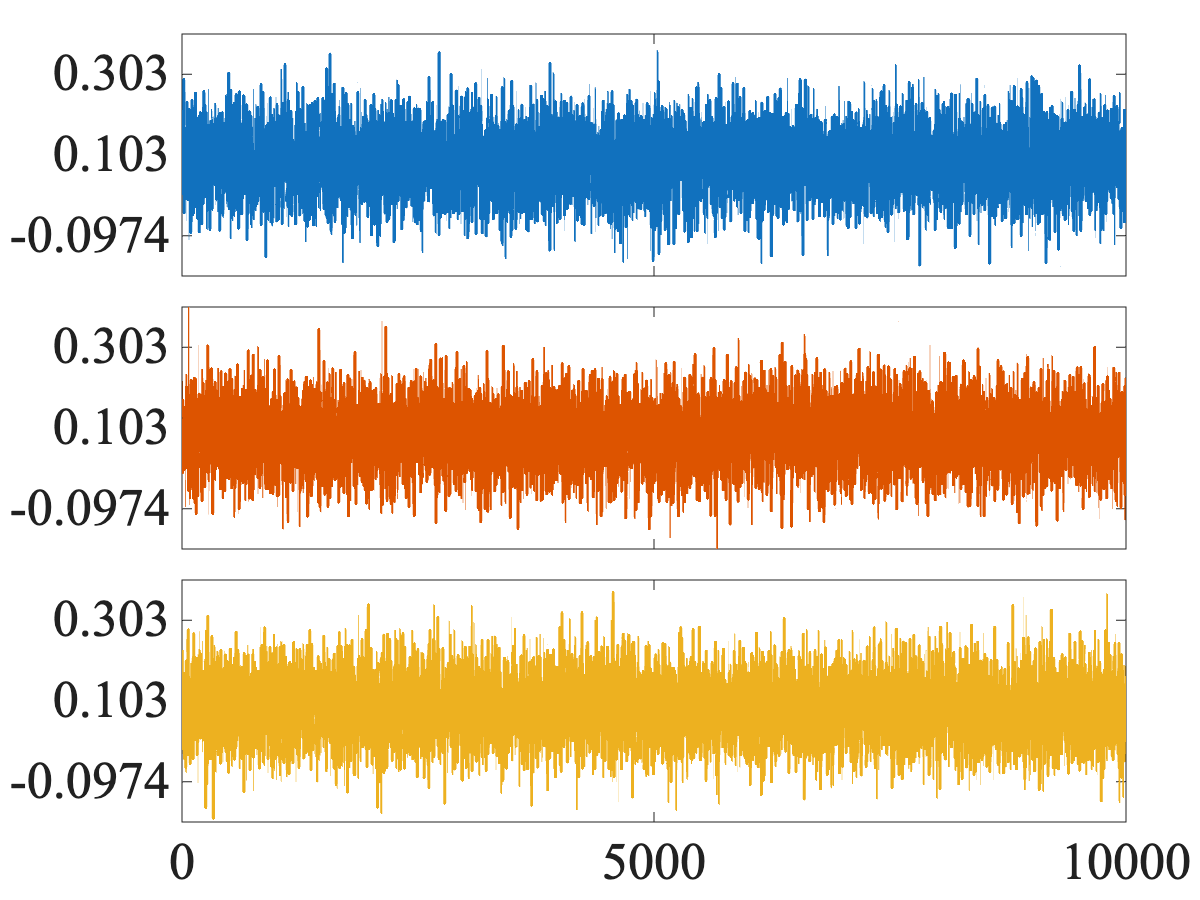}
    \end{tabular}
    \caption{The trace plot of a random point $x_j$ in $\bfx$ using \Cref{alg:gibbs_genGK} for the \texttt{PRseismic} test problem where $k$ denotes the rank of $\widehat\bfGamma$. }
    \label{fig:seismic_tracex}
\end{figure}

The trace plots and the estimated integrated ACF for the $\lambda$ and $\delta$ chains are provided in \Cref{tab:lamchain} and \Cref{tab:delchain}
respectively. For the $\lambda$ chain, each of the ACF quickly decays to 0, meaning there is little correlation between samples. This can also be observed through the trace plot of $\lambda$. 
Compared to $\lambda$, the $\delta$ chain has a significantly smaller ESS for each rank as seen in \Cref{tab:seismic_diag} and requires a higher rank to exhibit a fast ACF decay.  \rev{For $k=500$, the ESS of a random element $x_i$ of $\bfx$ was found to be 636 with $\tau_{\rm int}\approx47.135$, the p-value is 0.865, and the 95\% confidence interval is $[-0.051,\ 0.232]$. For $k=1,000$, the ESS of $x_i$ was found to be 29,020 with $\tau_{\rm int}\approx1.033$, the p-value is 0.997, and the 95\% confidence interval is $[-0.059,\ 0.229]$. } The trace plot \rev{provided} in \Cref{fig:seismic_tracex} \rev{of the corresponding element} $x_i$ \rev{ and the p-value} give a good indication that the chain has little correlation and is in equilibrium for $k=1,000$. The same conclusion may not be reached for $k=500$.

\begin{table}[]
    \centering
    \begin{tabular}{|c|c|c|ccc|ccc|}
    \hline
    \multicolumn{9}{|c|}{$k=500$}
    \\ \hline 
        \multicolumn{1}{|c|}{} 
        & \textbf{Acc.}
        & \multirow{2}{*}{\textbf{Time (s)}}
        & \multicolumn{3}{c|}{$\lambda$} 
        & \multicolumn{3}{c|}{$\delta$}
        \\ 
        &\textbf{Rate}&& \textbf{p-value} & \textbf{ESS} & $\widehat{R}$ & \textbf{p-value} & \textbf{ESS} & $\widehat{R}$
        \\ \hline
        \textbf{genGK} & 10\% & 1,179 
        & 0.995 & 6,514 & 1.002 
        & 0.928 & 146.42 & 1.019 
        \\ \hline
        \textbf{tSVD} & 32\% & 4,108 
        & 0.999 & 13,028 & 0.999
        & 0.997 & 820 & 1.005
        \\ \hline
        \textbf{rSVD} & 0.5\% & 4,148  
        & 0.906 & 16 & 1.253 
        & 0.84 & 15 & 1.334
        \\ \hline\hline
        \multicolumn{9}{|c|}{$k=1,000$}
        \\ \hline 
        \multicolumn{1}{|c|}{} 
        & \textbf{Acc.}
        & \multirow{2}{*}{\textbf{Time (s)}}
        & \multicolumn{3}{c|}{$\lambda$} 
        & \multicolumn{3}{c|}{$\delta$}
        \\ 
        &\textbf{Rate}&& \textbf{p-value} & \textbf{ESS} & $\widehat{R}$ & \textbf{p-value} & \textbf{ESS} & $\widehat{R}$
        \\ \hline
        \textbf{genGK} &  98\% & 1,346
        & 0.998 & 30,767 & 0.999
        & 0.993 & 3,079 & 1.001
        \\ \hline
        \textbf{tSVD} & 99\% & 5,069 
        & 0.999 & 20,963 & 0.999
        & 0.988 & 3,313 & 1.001
        \\ \hline
        \textbf{rSVD} & 97\% & 5,086 
        & 0.999 & 20,316 & 0.999
        & 0.997 & 3,454 & 1.0001 \\
        \hline
    \end{tabular}
    \caption{Diagnostic results from Algorithm \ref{alg:gibbs_genGK} on the \texttt{PRseismic} test problem using genGK, tSVD, rSVD for approximations of rank $k=500$ and $k=1,000$. The Acc. Rate gives the acceptance rate from the Metropolis-Hastings step, the compute time is provided in seconds, the p-value gives the probability of accepting or rejecting the Geweke null hypothesis, the ESS provides the effective sample size, and $\widehat{R}$ gives the Gelman-Rubin statistic.}
    \label{tab:seismic_diag}
\end{table}

\subsection{Atmospheric inverse modeling}
\label{sub:atmospheric}

Now we consider an atmospheric transport problem where $\bfA\in\bbR^{98,880\times 3,222}$ represents a forward atmospheric transport model from NOAA's CarbonTracker-Lagrange project \cite{GeostatInvModel,DataReductInvModel} modeling simulations from the weather research and forecasting stochastic time-inverted Lagrangian transport model \cite{LagrangTransport,LagrangTransportModel} sampled at 3,222 grid locations covering North America. The goal is to approximate $\bfx\in\bbR^{3,222}$, the vectorized true average fluxes at the grid locations. The observations $\bfb\in\bbR^{98,880}$ are sampled from OCO-2 during July through mid-August 2015. For this problem, we do not use realistic CO$_2$ emissions. Instead, $\bfx_{\rm true}$ is a randomly generated vectorized emission map that is used to produce the observations $\bfb$. We add Gaussian white noise corresponding to a $50\%$ noise level to the observations, i.e., $\frac{\sigma\norm[2]{\bfxi}}{\norm[2]{\bfA\bfx_{\rm true}}}=0.50$ where $\sigma$ is the standard deviation and $\bfxi\sim\calN\left(\bf0,\bfI\right)$. Using this setup, the hyperparameter associated with the noise, $\lambda$, should \jmc{be} $\frac{1}{\sigma^2}\approx 16.51$ for our generated $\bfxi$. 

Let $\bfS\in\{0,1\}^{3,222\times 11,900}$ be the sampling matrix that extracts the 3,222 grid locations over North America from the entire $11,900 \times 11,900$ grid. Then, the prior covariance matrix $\bfQ\in\bbR^{3,222\times 3,222}$ can be created by sampling a Mat\'ern \jmc{covariance matrix, $\bfQ_M \in \bbR^{11,900\times 11,900}$}, with $\nu=2.5$ and $\ell = 0.05$, \jmc{as}
\begin{equation}
    \bfQ = \bfS\bfQ_M\bfS\t.
\end{equation}

\begin{figure}[bthp]
    \begin{tabular}{ccc}
         \textbf{True} & 
         \textbf{Mean} & \textbf{Variance}  \\
         \includegraphics[width=.32\textwidth]{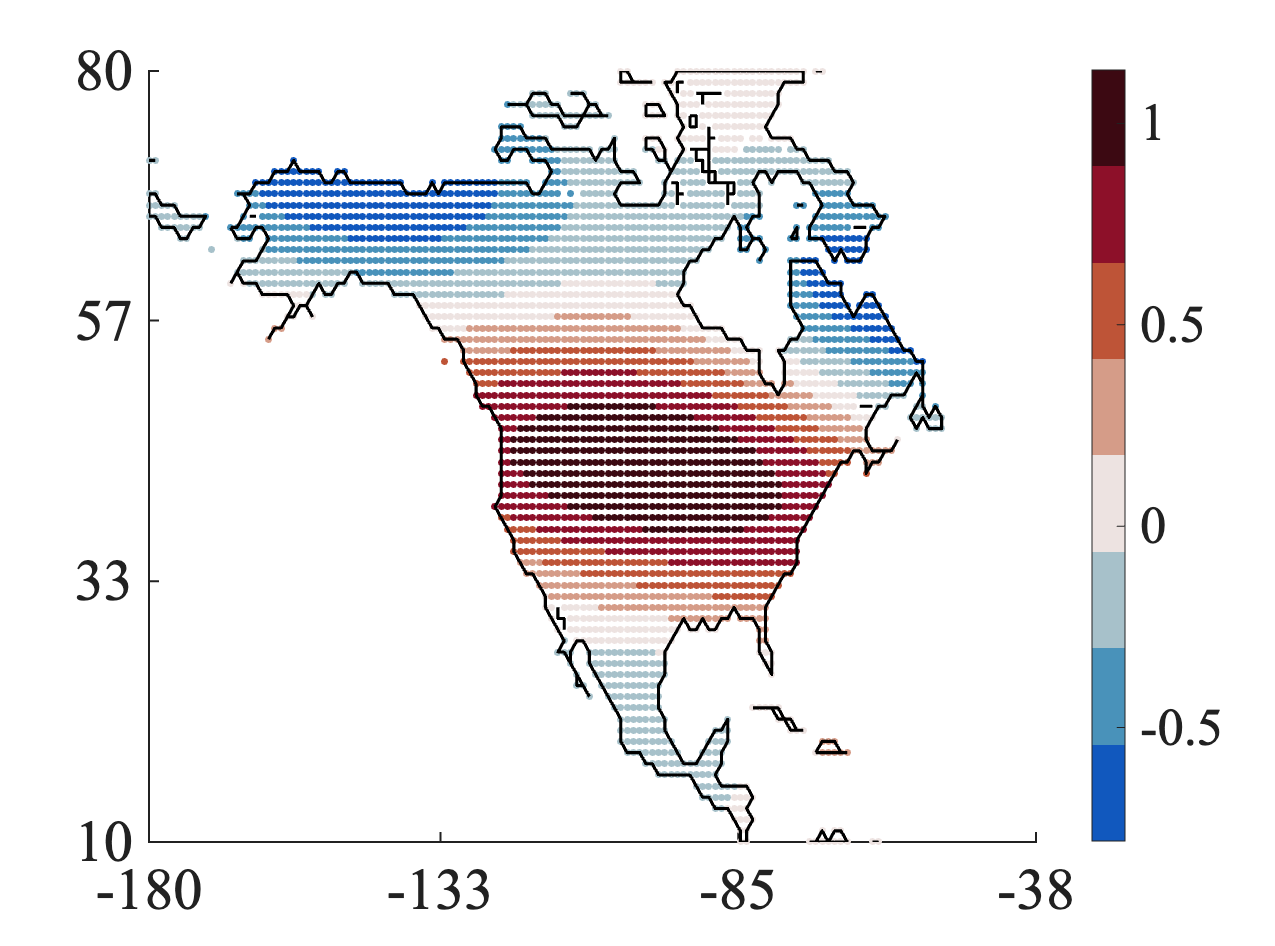}& 
         \includegraphics[width=.32\textwidth]{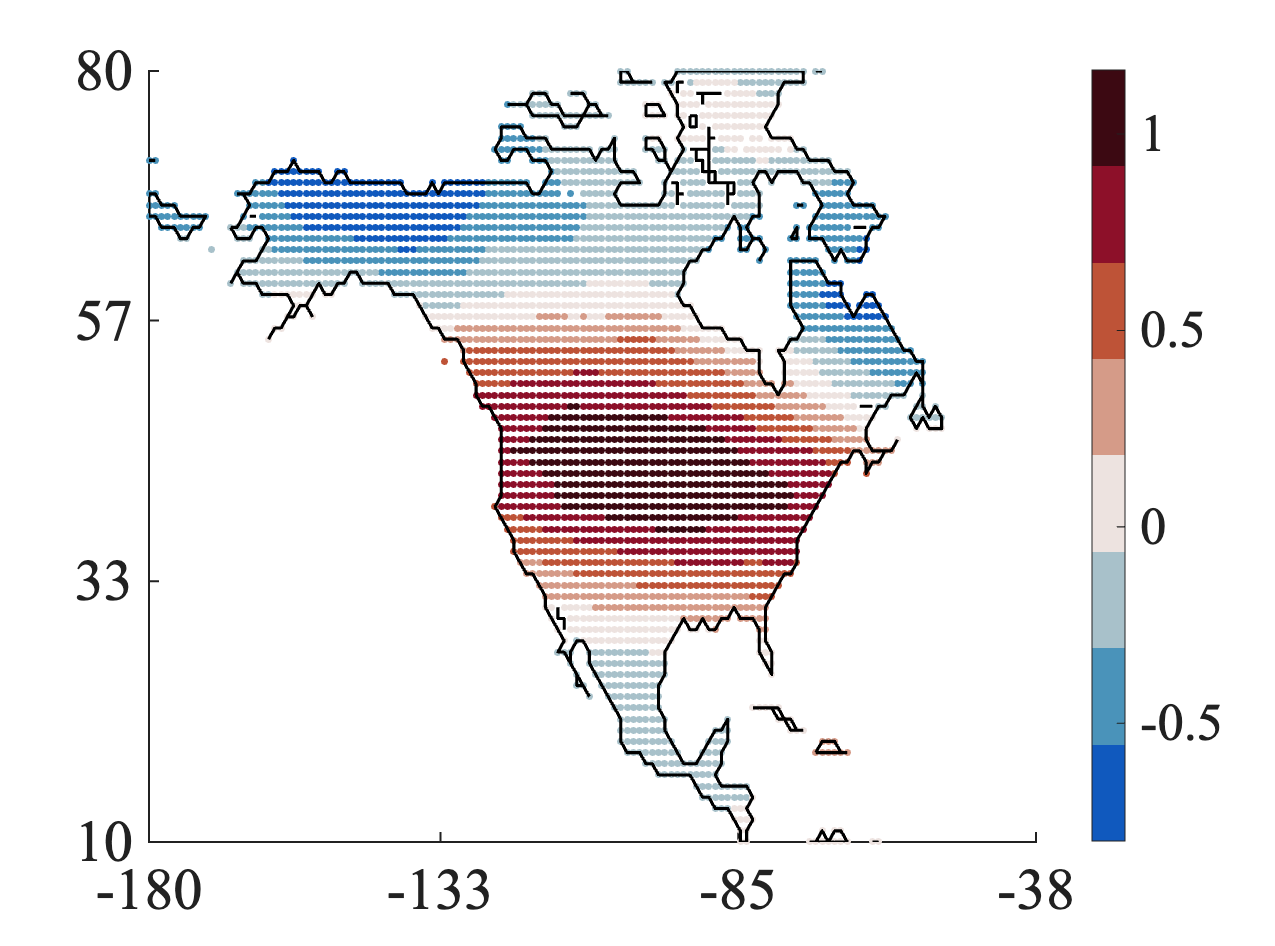} & \includegraphics[width=.32\textwidth]{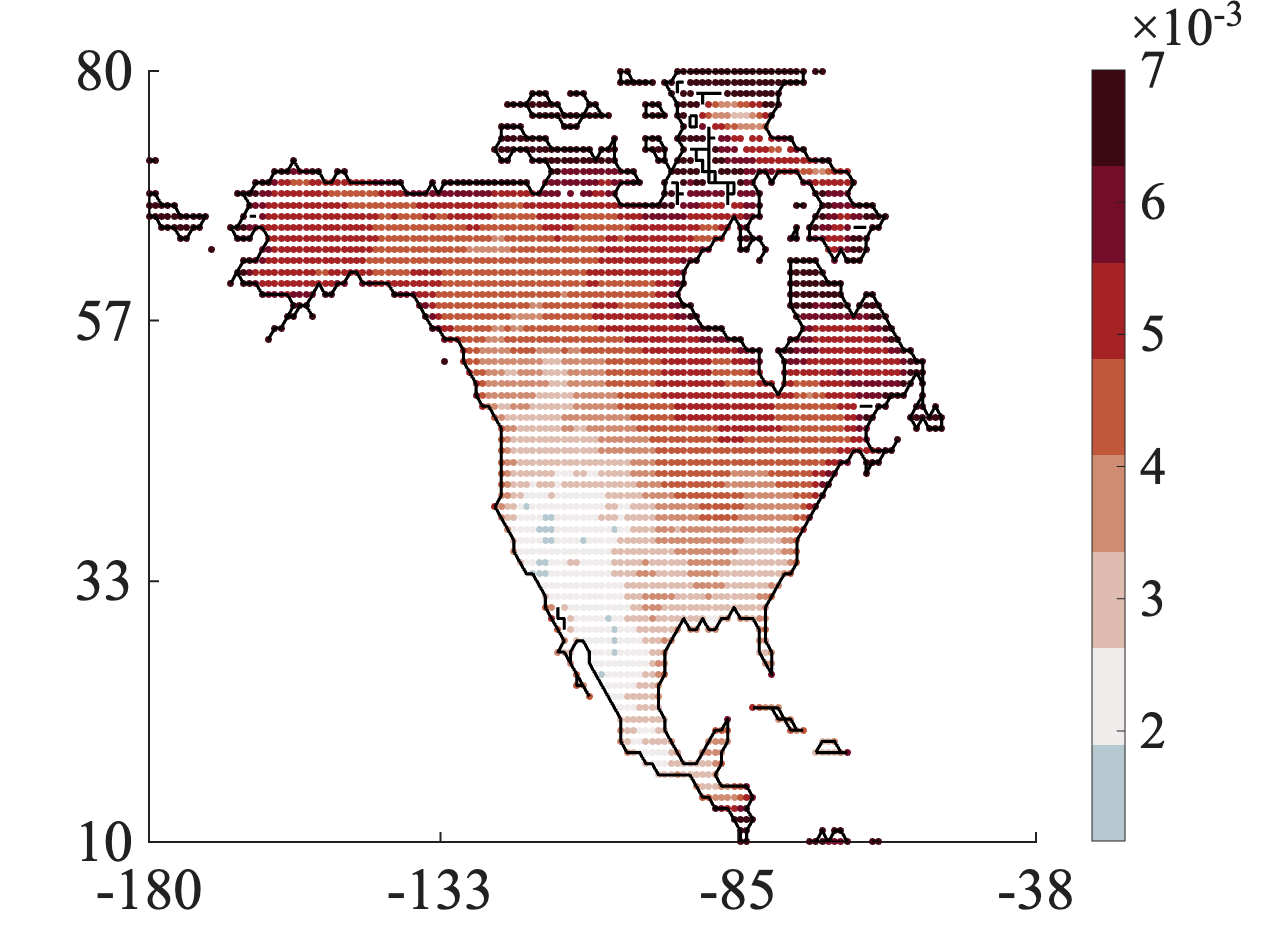} \\ &
         $\lambda$ Distribution & $\delta$ Distribution \\ &
         \includegraphics[width=.32\textwidth]{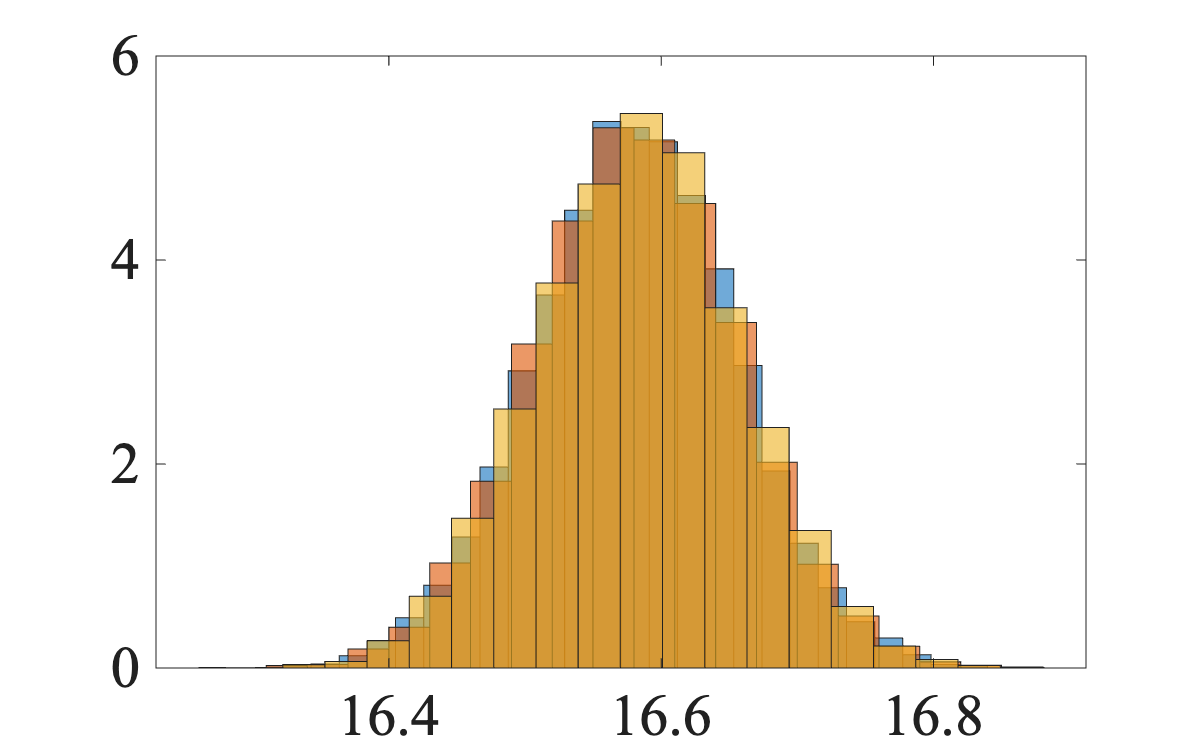} & \includegraphics[width=.32\textwidth]{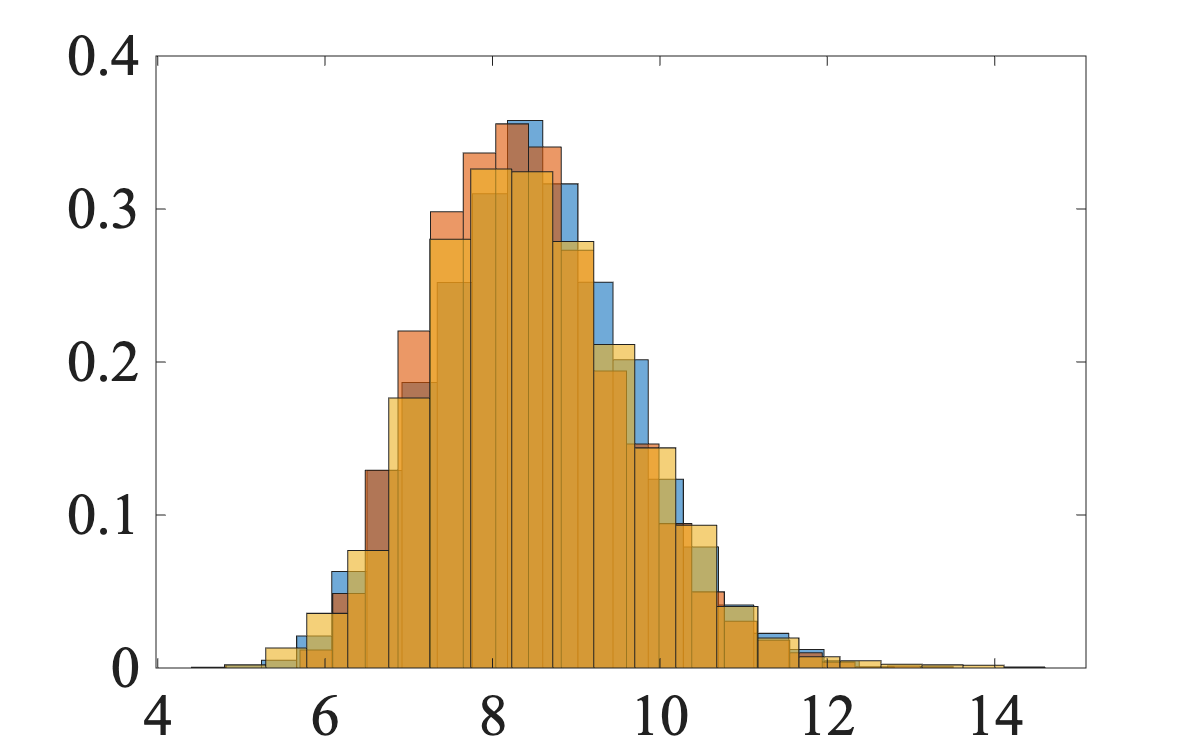}
    \end{tabular}
    \caption{For the atmospheric transport problem, we provide the true image and results for $T=20,000$ samples using \Cref{alg:gibbs_genGK} where $\widehat{\bfGamma}$ has a rank of $k=750$. The mean and variance are taken over the accepted samples after burn-in. The $\lambda$ and $\delta$ distributions are normalized histograms containing all draws from $\pi_\lambda$ and $\pi_\delta$ after 50\% burn-in.}
    \label{fig:atmos}
\end{figure}

\begin{figure}[bthp]
    \centering
    \includegraphics[width=.45\textwidth]{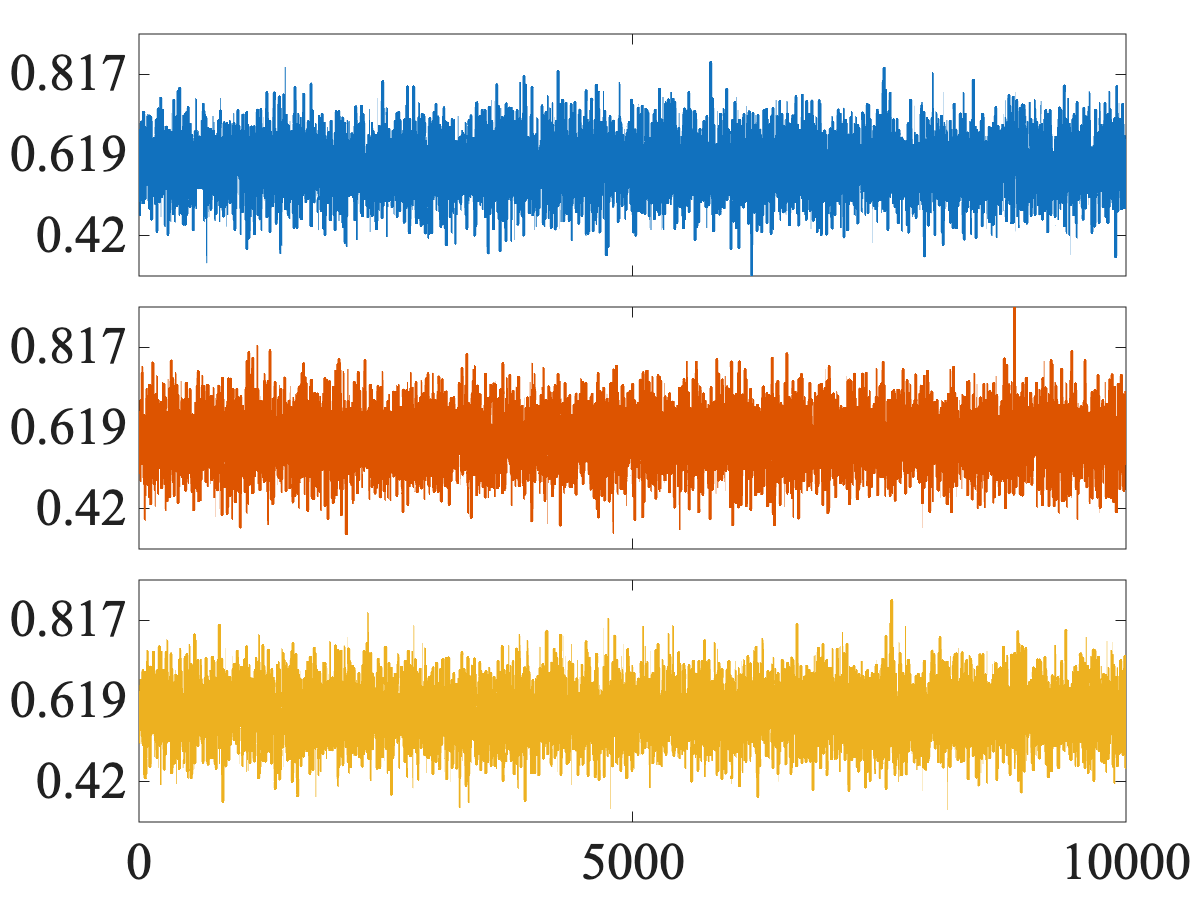}
    \caption{The trace plot of a random point $x_j$ in $\bfx$ using \Cref{alg:gibbs_genGK} for the atmospheric transport problem \jmc{for three chains} where $\widehat{\bfGamma}$ has a rank of $k=750$.}
    \label{fig:atmo_trace}
\end{figure}

\begin{figure}[bthp]
\centering
\begin{tabular}{ccc}
\textbf{Param.} & \textbf{Samples} & \textbf{\jmc{ACF}} \\  

\raisebox{1.8cm}{$\lambda$} & \raisebox{.05cm}{\includegraphics[width=.37\linewidth]{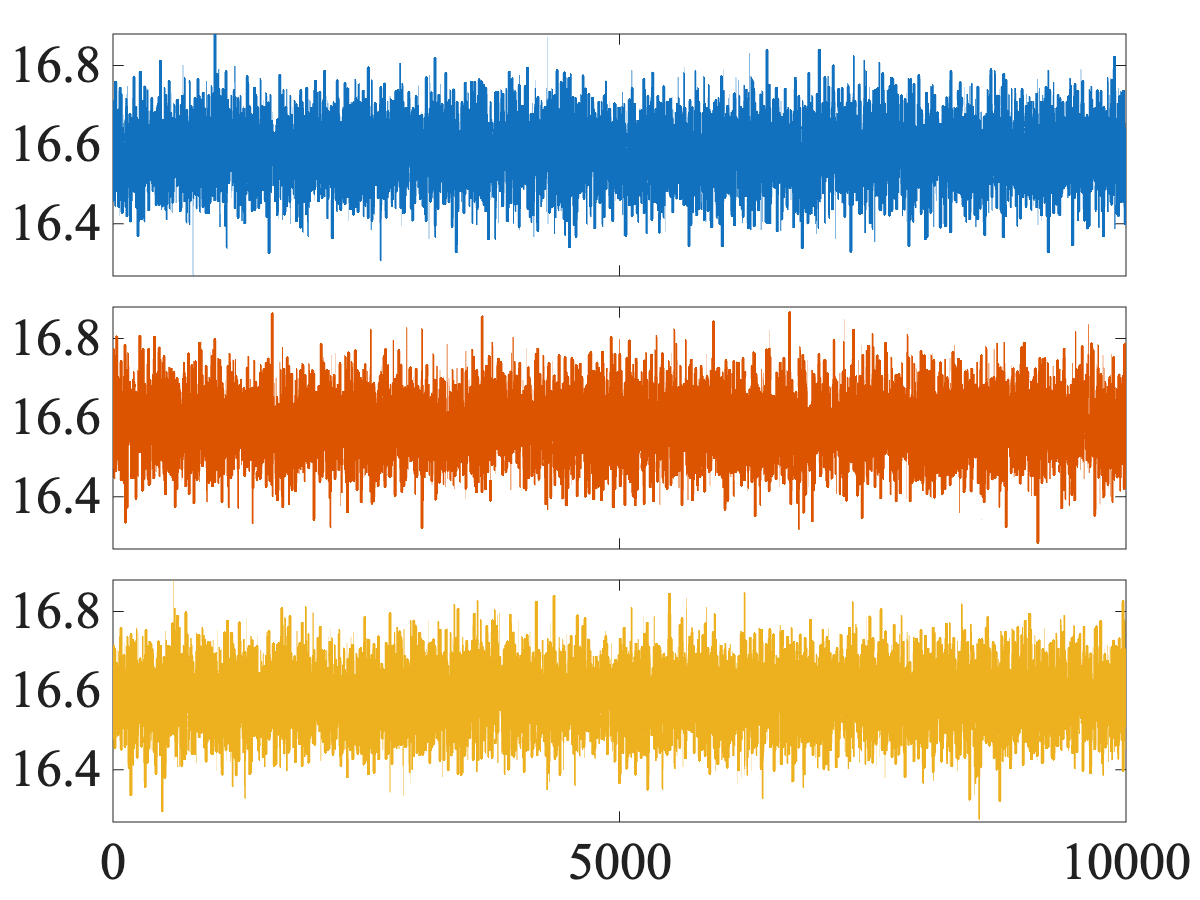}} & \includegraphics[width=.285\linewidth]{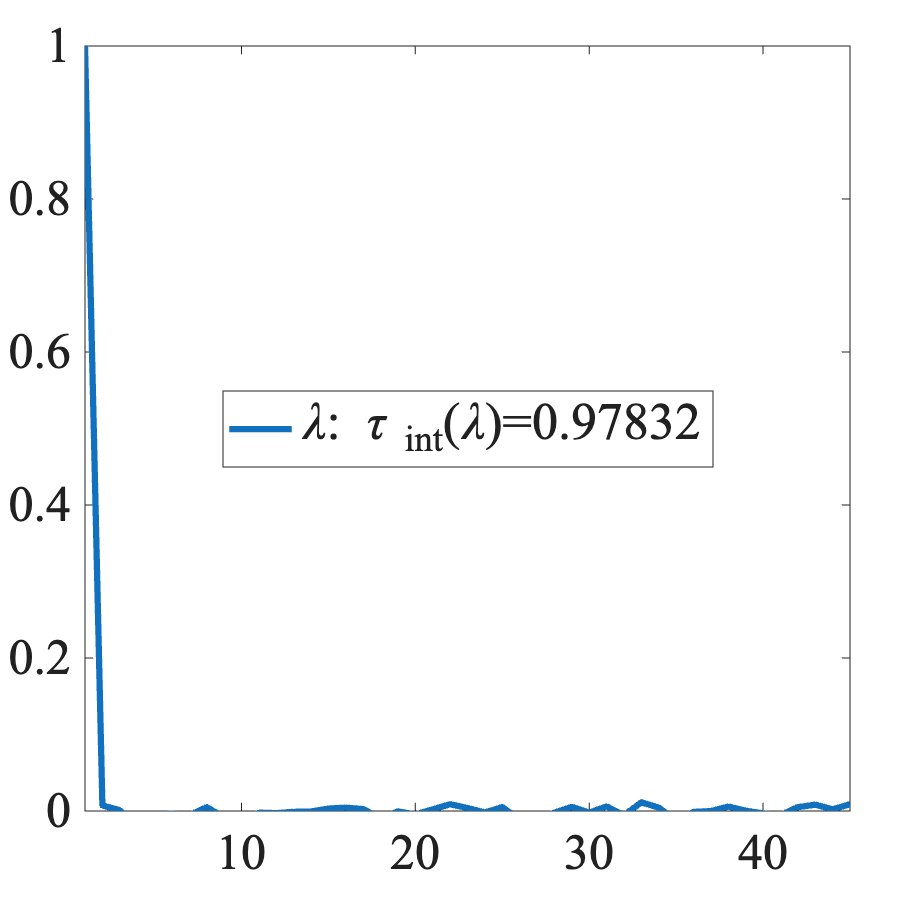} 
\\
\raisebox{1.8cm}{$\delta$} & \raisebox{.05cm}{\includegraphics[width=.37\linewidth]{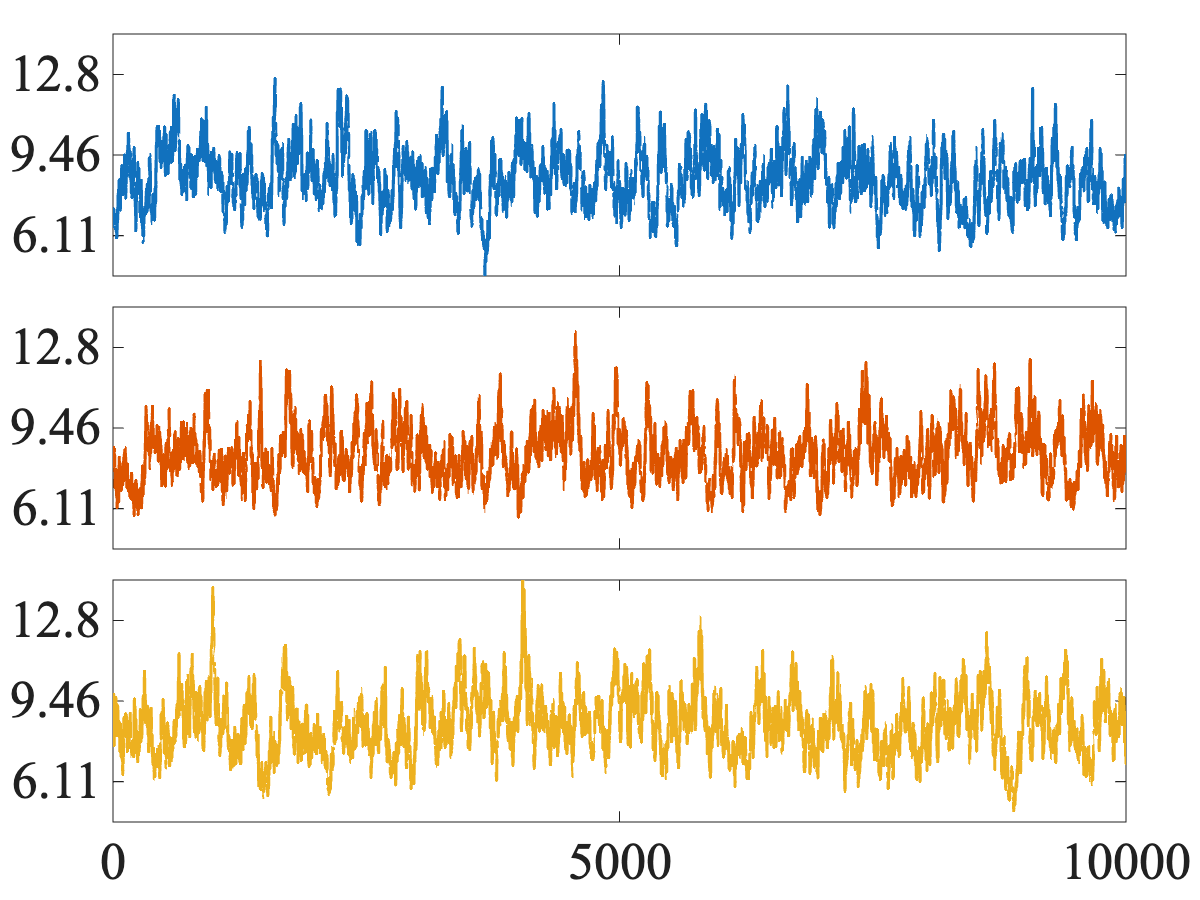}} & \includegraphics[width=.285\linewidth]{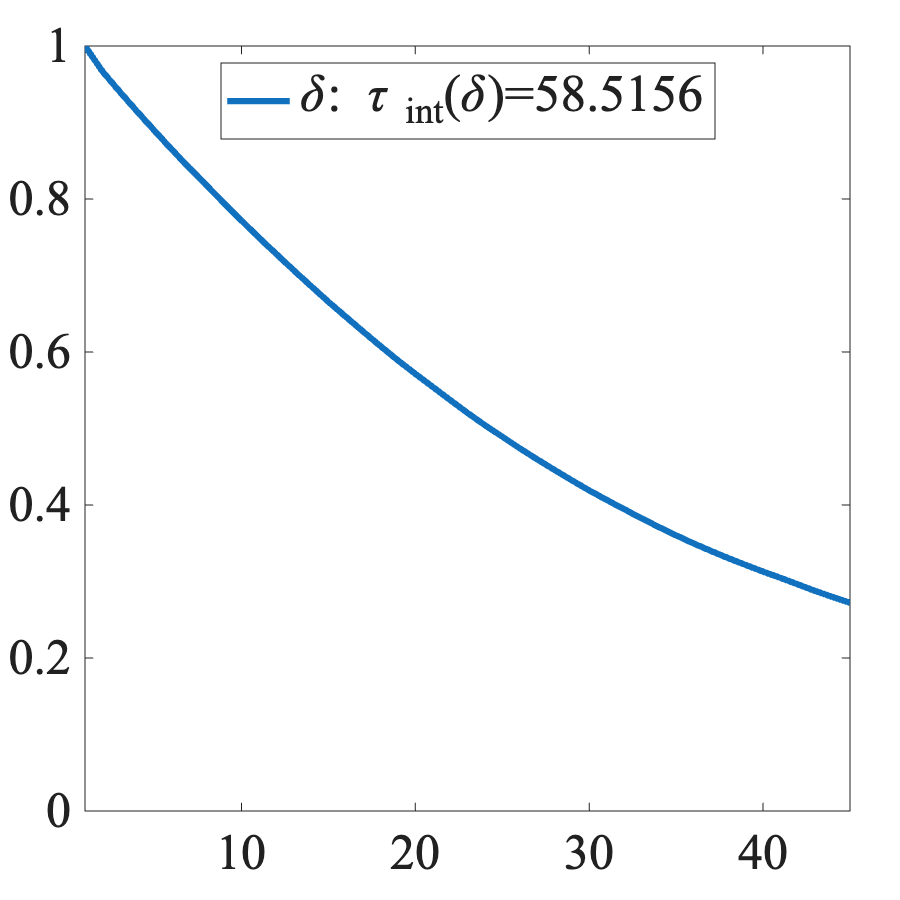}
\\
\end{tabular}
\caption{The $\lambda$ and the $\delta$ chain for $T=20,000$ samples using Algorithm \ref{alg:gibbs_genGK} on the atmospheric transport problem where $\widehat{\bfGamma}$ has a rank of $k=750$. We provide a trace plot of the sample values of $\lambda^j$ and $\delta^j$ \jmc{for three chains}, 
and we provide
the corresponding \jmc{ACF}.}
\label{tab:atmochain}
\end{figure}

We compute $T=20,000$ samples using the Metropolis-Hastings within Gibbs method with proposal sampling using the genGK approximation (\Cref{alg:gibbs_genGK}) with rank $k=750$. In \Cref{fig:atmos} we provide the true image along with the mean and variance images of the samples after burn-in and the normalized distributions of computed hyperparameters $\lambda$ and $\delta$. It should be noted that when plotting the variance of the solution, the colorbar is set to have a maximum value of $0.07$ to give a better representation of the variance throughout the grid by excluding the high variances found only at the boundaries. 
This setup produced an acceptance rate of $86\%$ and it was found that the $\lambda$ chain \jmc{has} a $95\%$ confidence interval of $[16.44,\ 16.73]$. 
The p-value of $\lambda$ was $0.999$ and the the p-value of $\delta$ was 0.968,
indicating strong evidence that the chains are in equilibrium.  Moreover, the values of $\widehat{R}$ were both less than 1.01 ($1.0001$ for $\lambda$ and $1.001$ for $\delta$).
Finally, the ESS for $\lambda$ was $30,647$, so most of the accepted samples are independent, while the ESS for $\delta$ was $512$, so there is higher correlation between the $\delta$ samples. The trace plots and the estimated integrated ACFs for the $\lambda$ and $\delta$ chains are provided in \Cref{tab:atmochain}. The ACF for $\lambda$ also provides evidence that there is little to no correlation in the chain. \rev{The ESS of a random element $x_i$ of $\bfx$ was found to be 22,334 with $\tau_{\rm int}\approx1.343$, the p-value is 0.997, and the 95\% confidence interval is $[0.461,\ 0.721]$. The trace plot in \Cref{fig:atmo_trace} for the corresponding element $x_i$ and the p-value give good indication that the chain has little correlation and is in equilibrium. The relative reconstruction error between the mean of all accepted samples and the ground truth is $0.09$.}


\subsection{Dynamic photoacoustic tomography}
\label{sub:PAT}
Finally, we consider a dynamic photoacoustic tomography test problem, based on the \texttt{PRseismic} example in IRTools \cite{gazzola2019ir}. For this example, 20 $64\times 64$ true images were generated using two Gaussians moving in different directions. Such problems require a spatiotemporal prior, and we use a Kronecker product $\bfQ = \bfQ_t \kron \bfQ_s \in\bbR^{81,920 \times 81,920}$ where $\bfQ_t\in\bbR^{20\times 20}$ and $\bfQ_s\in\bbR^{4,096\times 4,096}$ are temporal and spatial priors corresponding to Mat\'ern kernels with $\nu=2.5, \ \ell=0.1$ and $\nu=0.5, \ \ell=0.25$ respectively. We note that although the Kronecker product structure can be exploited for efficient computations, there are many scenarios where exact factorizations of $\bfQ$ are not computationally feasible, e.g.,  covariance matrices defined by three-dimensional covariance kernels or kernels defined on unstructured spatial grids.

The linear problem can be formed by combining the subproblem at each time point $i$ as follows,
$$
\bfx = \begin{bmatrix} \bfx_1 \\ \vdots \\ \bfx_{20} \end{bmatrix}, \quad \bfA = \begin{bmatrix}
    \bfA_1 & & \\ & \ddots & \\ & & \bfA_{20}
\end{bmatrix}, \quad \text{and} \quad \bfb = \begin{bmatrix}
    \bfb_1 \\ \vdots \\ \bfb_{20}
\end{bmatrix},
$$
where $\bfA_i \in \bbR^{1,638 \times 4,096}$ is a seismic projection matrix corresponding to 18 equally spaced angles between $i$ and $340+i$ for $i=1,\ldots,20$, $\bfx_i\in\bbR^{4,096}$ is the vectorized image at $i$, and $\bfb_i \in\bbR^{1,638}$ is the simulated projection data. We add Gaussian white noise corresponding to a $2\%$ noise level to the observations, i.e., $\frac{\sigma\norm[2]{\bfxi}}{\norm[2]{\bfA\bfx_{\rm true}}}=0.02$ where $\sigma$ is the standard deviation and $\bfxi\sim\calN\left(\bf0,\bfI\right)$. Using this setup, the hyperparameter associated with the noise, $\lambda$, should be $\frac{1}{\sigma^2}\approx 2.02 \times 10^4$.
Given the small size of $\bfQ_t$, we can obtain the Cholesky factorization $\bfQ_t^{-1} = \bfG_t\t\bfG_t$. Now we can define a preconditioner of the form $\bfG=\bfG_t \kron \bfG_s$ with $\bfG_s$ being the Cholesky factorization of $(-\Delta)^{\gamma}$ for $\gamma\geq 1$ where $\Delta$ is the Laplacian operator discretized using finite difference. 

\begin{figure}[h]
    \begin{subfigure}[]{1\textwidth}
    \centering
    \includegraphics[width=1\textwidth]{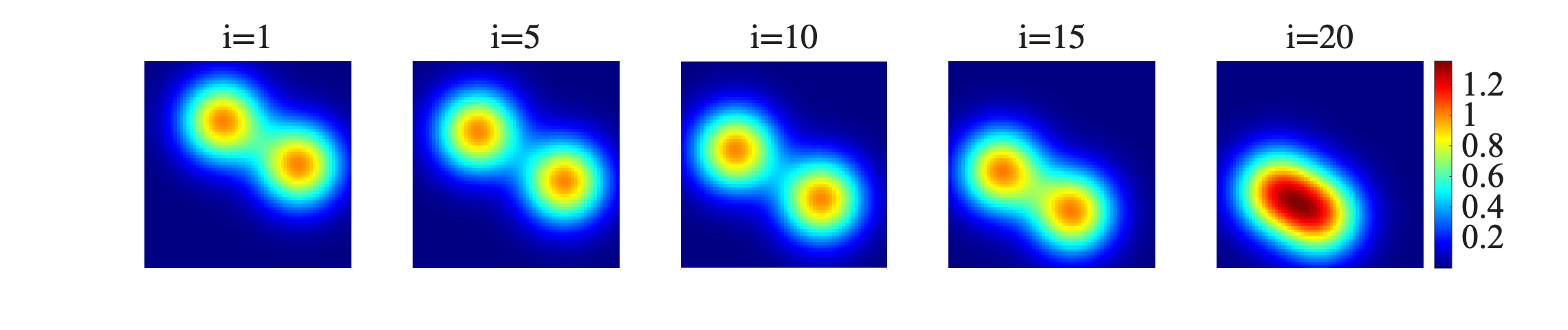}    
    \end{subfigure} \hfill
    \begin{subfigure}[]{1\textwidth}
    \centering
    \includegraphics[width=1\textwidth]{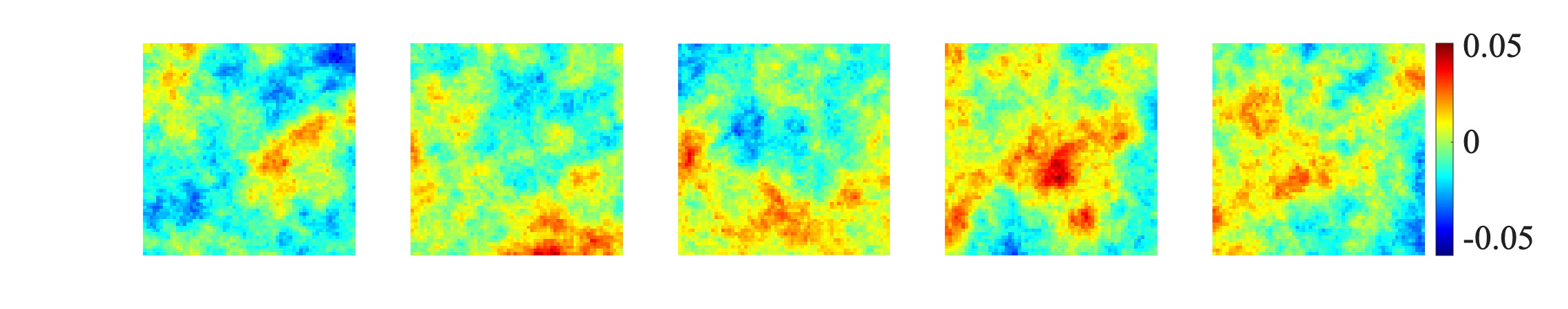}    
    \end{subfigure} \hfill 
    \begin{subfigure}[]{1\textwidth}
    \centering
    \includegraphics[width=1\textwidth]{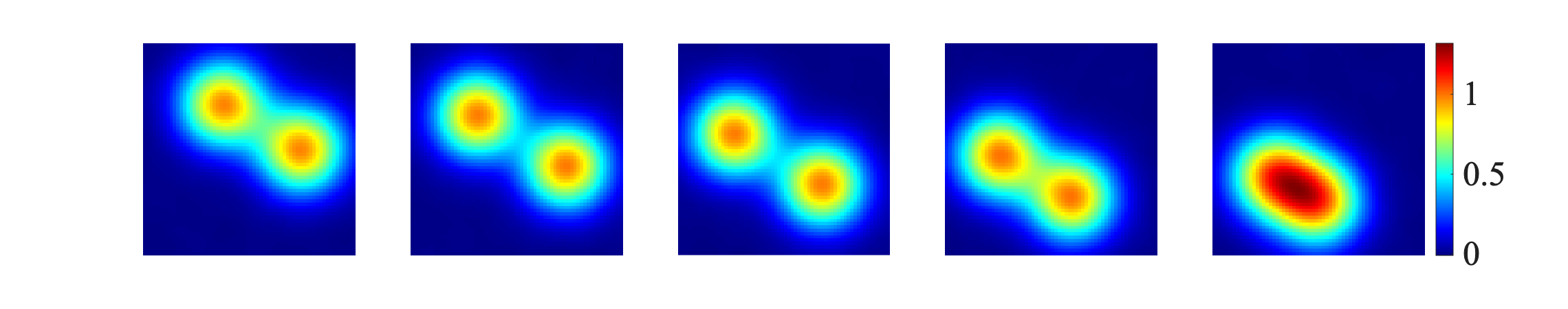}    
    \end{subfigure} 
    \begin{subfigure}[]{1\textwidth}
    \centering
    \includegraphics[width=1\textwidth]{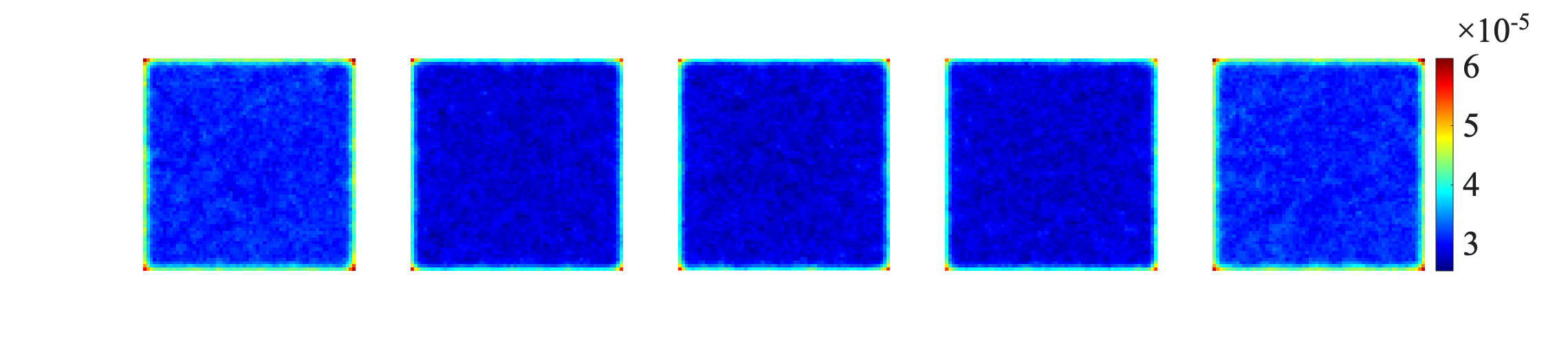}    
    \end{subfigure} 
    \caption{For the dynamic photoacoustic tomography test problem, the true images at time points $i=1,5,10,15,$ and $20$ are provided in the top row. The results correspond to Algorithm \ref{alg:gibbs_precond} with an approximate mean obtained using $k=200$ genGK iterations. Second row: Random samples from the prior. Third row: The means of the accepted samples after burn-in.
    Bottom: The variances of the accepted samples after burn-in.}
    \label{fig:PAT500}
\end{figure}

\begin{figure}[bthp]
    \centering
    \begin{tabular}{cc}
     $\lambda$ Distribution & $\delta$ Distribution  \\
     \includegraphics[width=.35\textwidth]{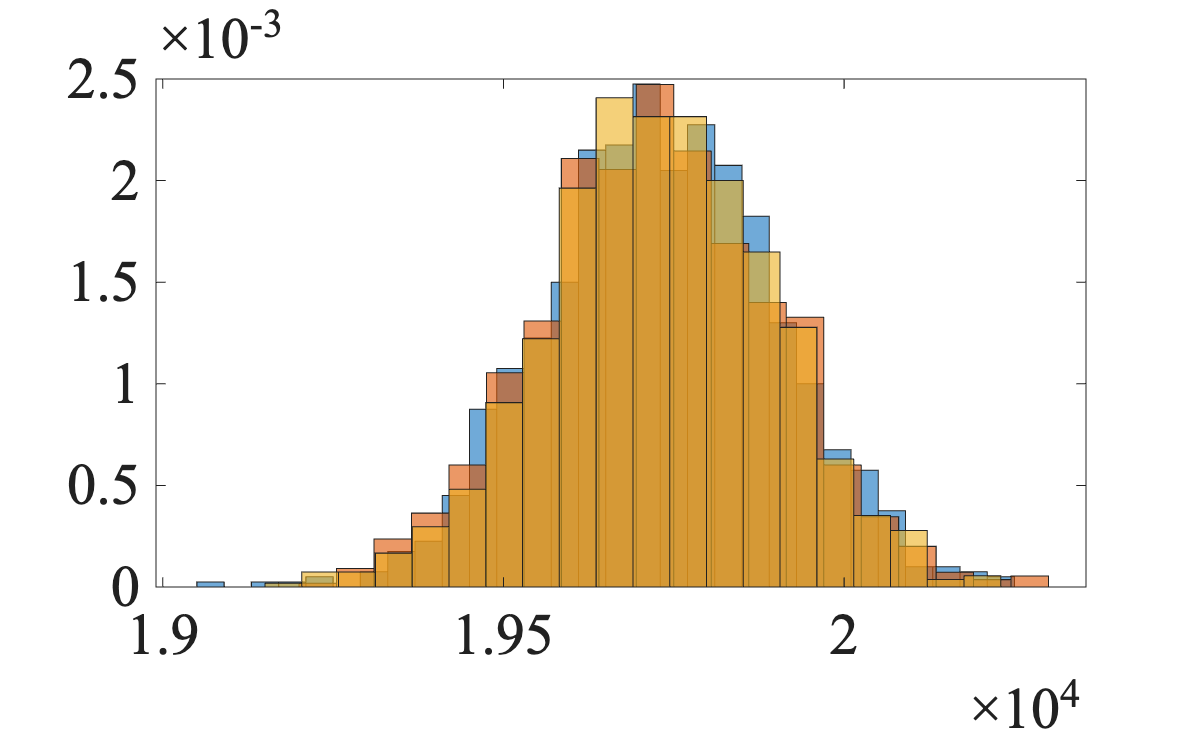}    & 
     {\includegraphics[width=.35\textwidth]{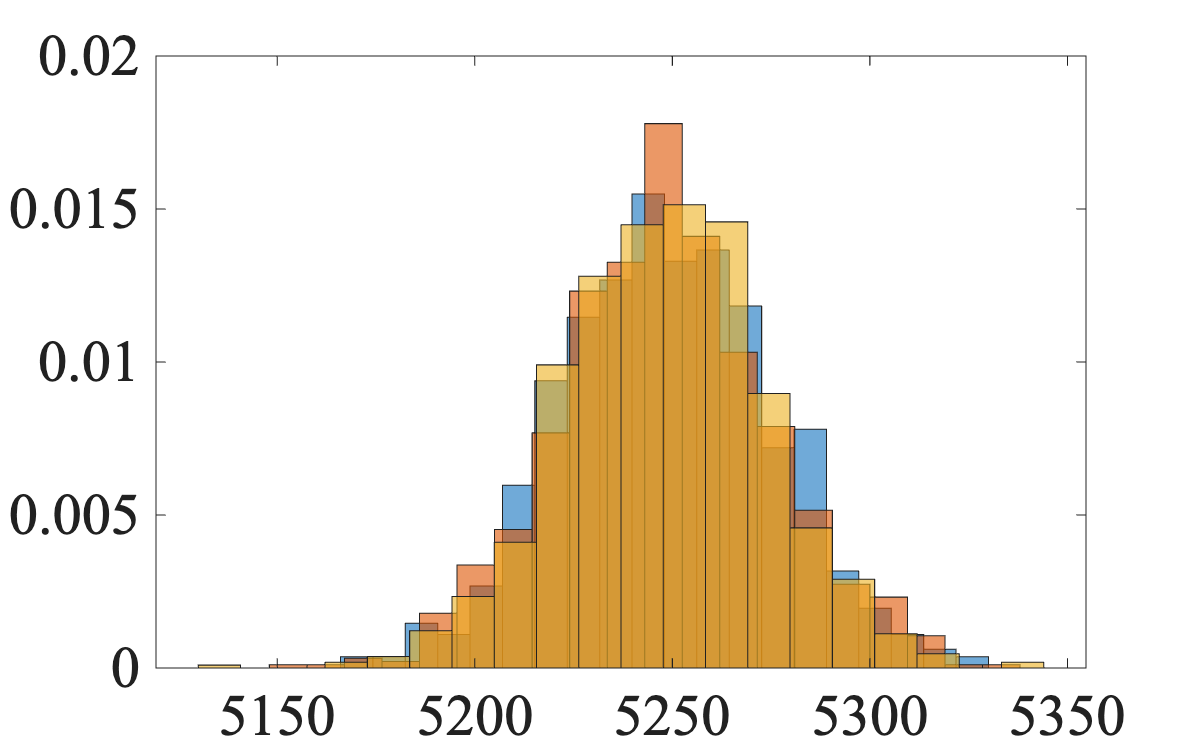}}
    \end{tabular}
    
    \caption{The $\lambda$ and $\delta$ distributions are normalized histograms containing all draws from $\pi_\lambda$ and $\pi_\delta$ after 50\% burn-in from Algorithm \ref{alg:gibbs_precond} for the dynamic photoacoustic tomography test problem with an approximate mean using $k=200$ genGK iterations, \jmc{for three chains}. }
    \label{fig:PAT500_dist}
\end{figure}

\begin{figure}[bthp]
\centering
\begin{tabular}{ccc}
\textbf{Param.} & \textbf{Samples} & \hspace{-.1cm}\textbf{\jmc{ACF}} \\  

\raisebox{1.8cm}{$\lambda$} & \raisebox{.05cm}{\includegraphics[width=.37\linewidth]{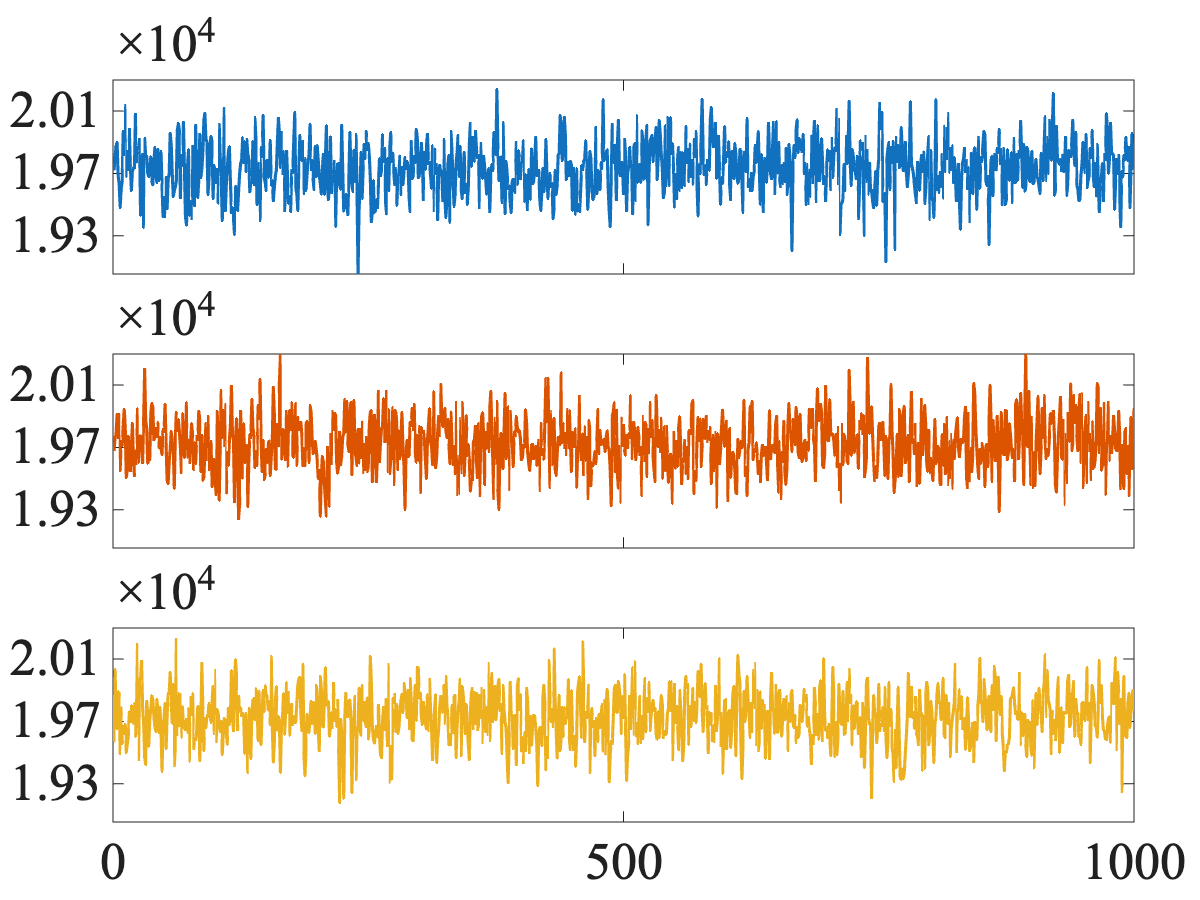}} & {\includegraphics[width=.285\linewidth]{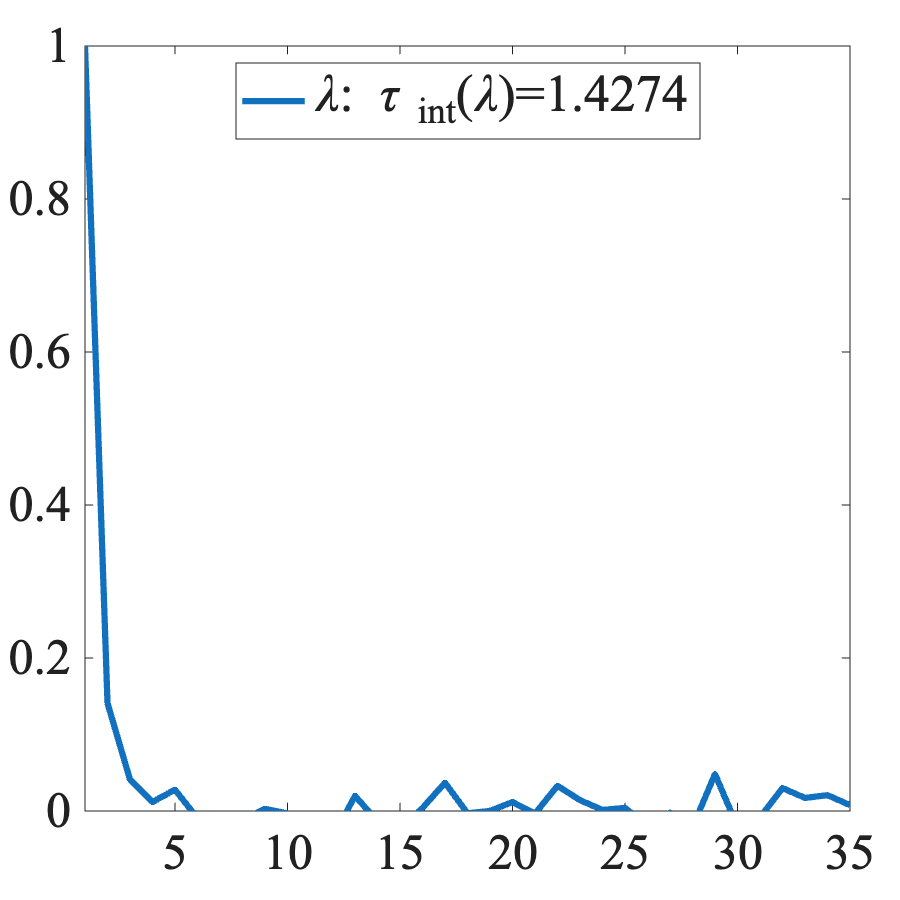}} 
\\
\raisebox{1.8cm}{$\delta$} & \raisebox{.05cm}{\includegraphics[width=.37\linewidth]{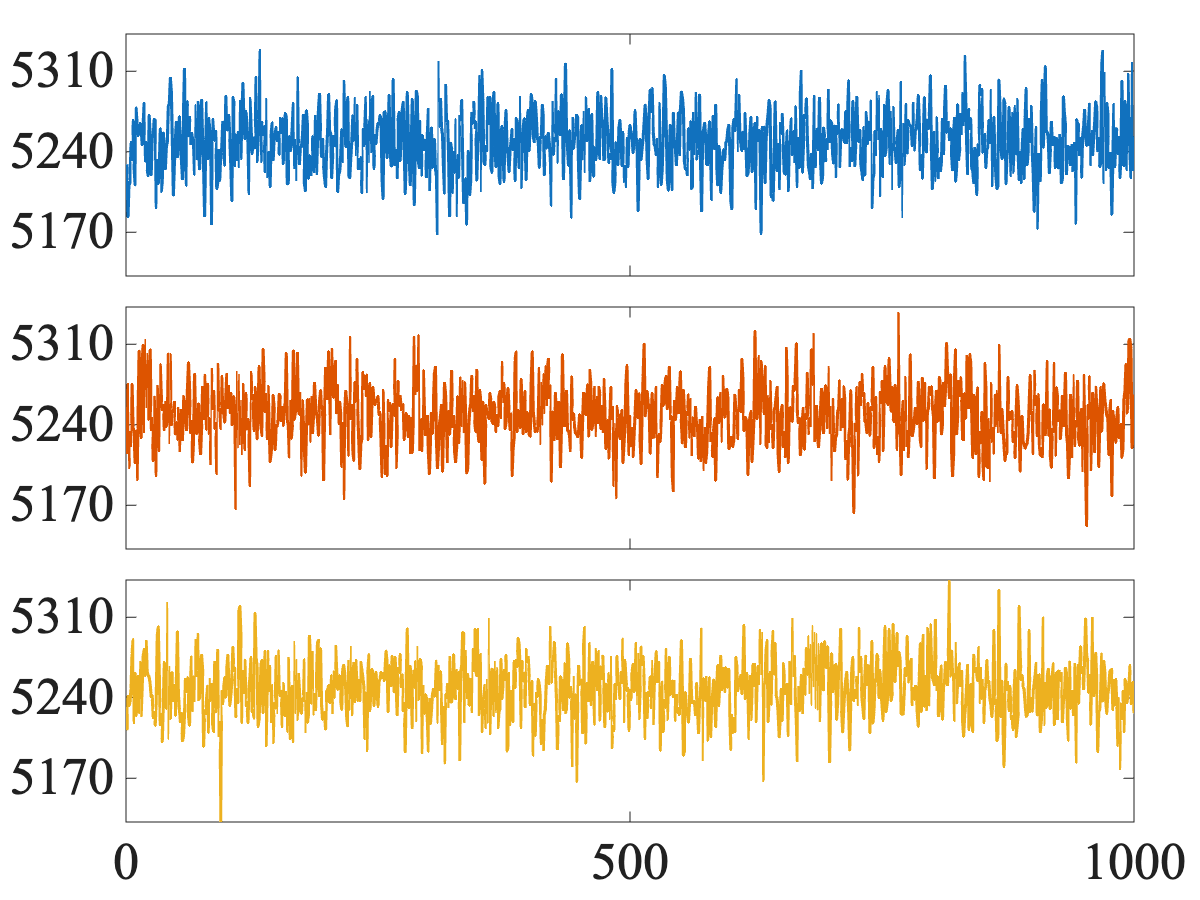}} & {\includegraphics[width=.285\linewidth]{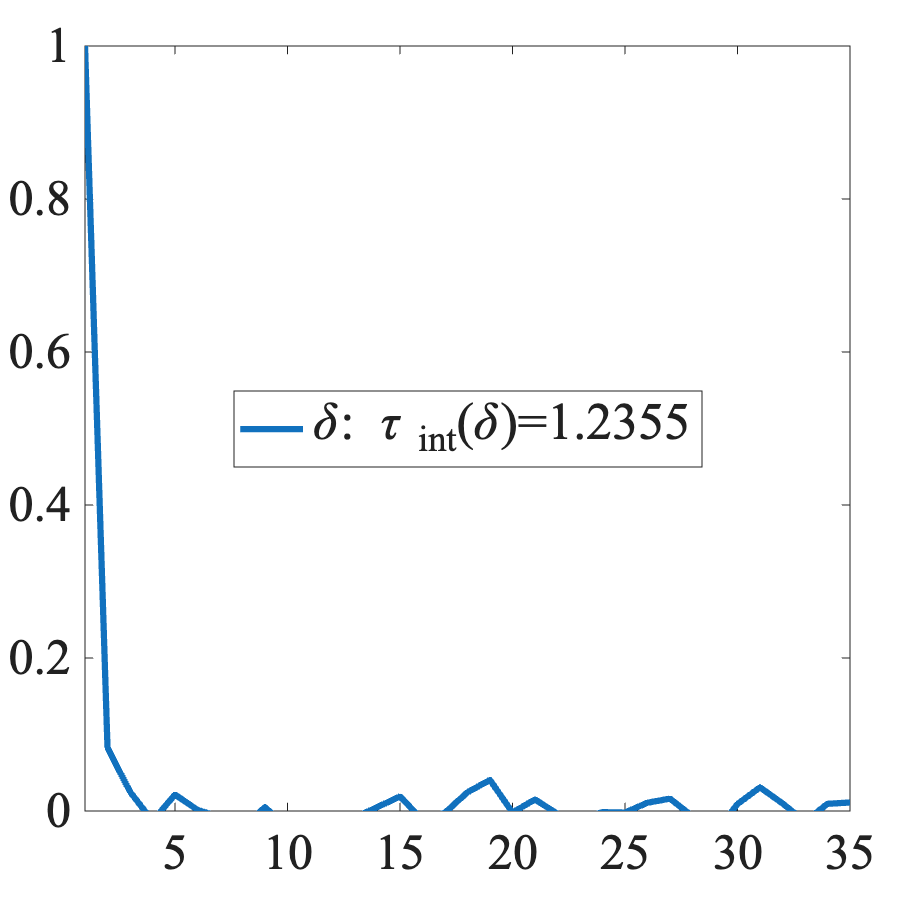}}
\\
\end{tabular}
\caption{The $\lambda$ and the $\delta$ chain for $T=2,000$ samples using Algorithm \ref{alg:gibbs_precond} on the photoacoustic tomography test problem with an approximate mean using $k=200$ genGK iterations. We provide a trace plot of the sample values of $\lambda^j$ and $\delta^j$ \jmc{for three chains}, and
corresponding \jmc{ACF}.}
\label{tab:genGK_dynamic}
\end{figure}

\begin{figure}[bthp]
    \centering
    \includegraphics[width=.45\textwidth]{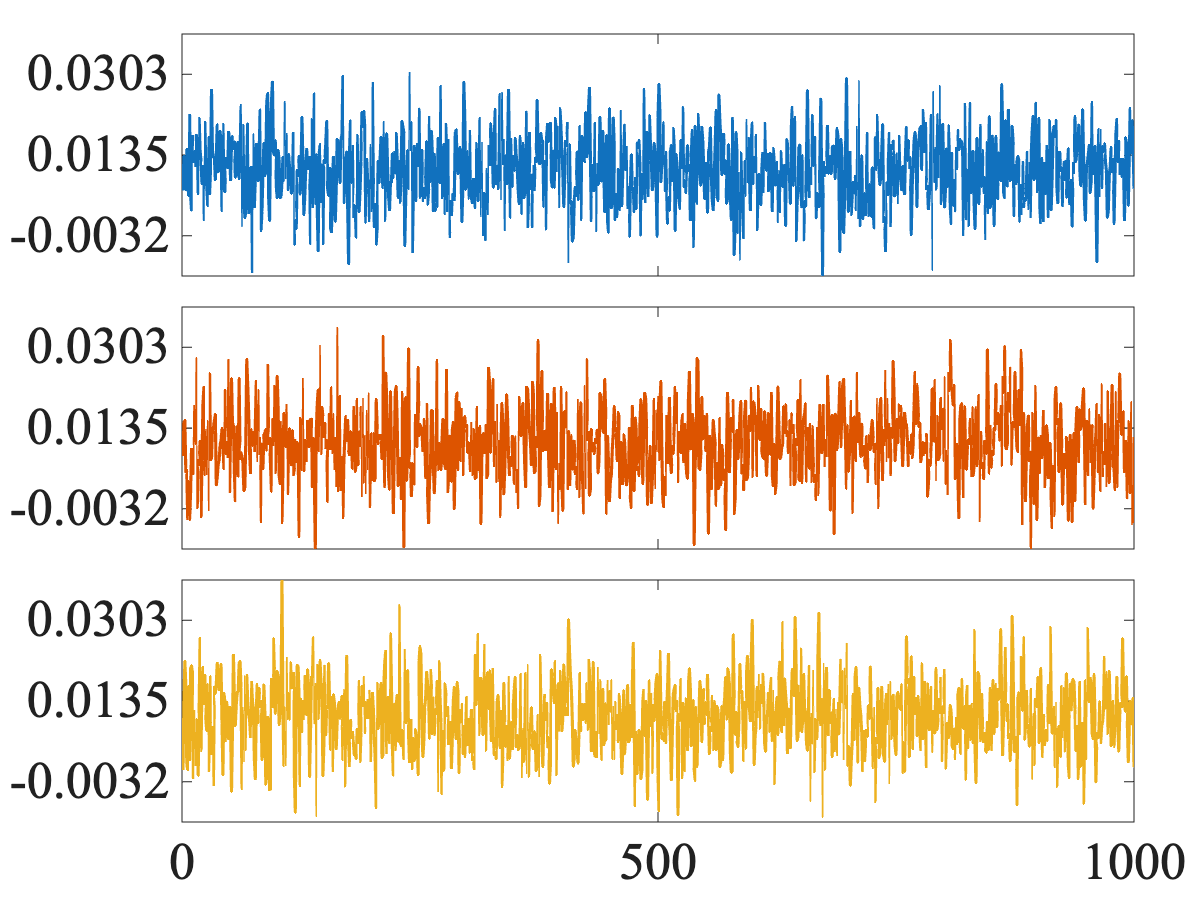}
    \caption{The trace plot of a random point $x_j$ in $\bfx$ using \Cref{alg:gibbs_precond} for a photoacoustic tomography test problem with an approximate mean using $k=200$ genGK iterations, \jmc{for three chains}.}
    \label{fig:pat_trace}
\end{figure}

The goal of this problem is to obtain a sequence of images from a sequence of projection datasets. We used a genGK approximation to the target distribution, but due to the large rank of this problem, the proposal failed to accept a sufficient number of samples without taking $k$ to be very large. Thus, we used the preconditioned Lanczos method from \Cref{alg:gibbs_precond}. In \Cref{fig:PAT500} we provide the true images alongside the results for $T=2,000$ samples. The approximate mean, $\bfx_k$, is obtained using $k=$200 genGK iterations. Additionally, the corresponding $\lambda$ and $\delta$ distributions are given in \Cref{fig:PAT500_dist}.  We also provide a sample from the prior.
We observe that the method produced a $100\%$ acceptance rate. While this is higher than that of the previous examples, it is expected, given that an exact factorization of $\bfGamma_{\rm cond}$ is used.
While the $\lambda$ distribution approaches the true value, it is outside the $95\%$ confidence interval of $[19,399, \ 20,056]$. The p-values of $\lambda$ and $\delta$ (both $0.999$) show strong evidence that the chains are in equilibrium, and values of $\widehat{R}$ \jmc{are} less than 1.01 (both 0.999). For $\lambda$ and $\delta$, an ESS of $2,102$ and $2,428$ respectively, show that most of the accepted samples are independent for both chains. For this example, the ESS for $\lambda$ and $\delta$ are comparable meaning \Cref{alg:gibbs_precond} may produce $\delta$ chains which are less correlated than those from \Cref{alg:gibbs_genGK}. In \Cref{tab:genGK_dynamic} we provide the trace plots and the estimated integrated ACF for the $\lambda$ and $\delta$ chains. The ACFs corresponding to $\lambda$ and $\delta$ quickly decay to 0, \rev{indicating that} the chains produced by \Cref{alg:gibbs_precond} are highly uncorrelated. \rev{The ESS of a random element $x_i$ of $\bfx$ was found to be 2,822 with $\tau_{\rm int}\approx1.063$, the p-value is 0.997, and the 95\% confidence interval is $[-0.00398,\ 0.02495]$.} The trace plot in \Cref{fig:pat_trace} for \rev{the corresponding} element $x_i$ \rev{and the p-value} give a good indication that the chain has little correlation and is in equilibrium.  The relative reconstruction error between the mean of accepted samples and the ground truth was $0.01.$

\section{Conclusions}
\label{sec:conclusions}
This paper provides an approach to perform \jmc{UQ} for large-scale hierarchical Bayesian inverse problems, where Metropolis-Hastings independence sampling within Gibbs is used to overcome the limitations of \jmc{MCMC} methods.  We consider two proposal distributions, both of which are based on \jmc{genGK} methods. First, \jmc{we define a proposal distribution using a low-rank genGK approximation to the conditional covariance matrix, where we can reuse the genGK matrices to efficiently draw proposal samples.} Second, for matrices where a low-rank approximation is not sufficient, we describe a preconditioned Lanczos-based method to efficiently perform computations with the square-root of the conditional covariance to draw proposal samples. In addition to comparisons with existing approaches for small problems, we demonstrate the performance of these methods on a variety of large-scale inverse problems, including atmospheric models and dynamic tomography problems.
Although it is known that Gibbs and Metropolis-Hastings within Gibbs sampling have potential drawbacks, few alternatives exist for efficient sampling for hierarchical large-scale inverse problems. The described approaches provide a computationally convenient approach for sampling, where \jmc{genGK} approximations can be used for proposal sampling. Future work includes further investigations, e.g., for scenarios with strong correlations and for other prior models, as well as incorporating further convergence diagnostics such as the multivariate version of the ESS.

\backmatter

\section*{Statements and Declarations}



\begin{appendices}
\section{Low-Rank representation of square-root covariance} \label{sec:lowrank}

First note that operations with the square-root matrix $\bfQ^{1/2}$ can be performed efficiently, using a Lanczos algorithm that only requires mat-vecs with $\bfQ$, see \cite{UQlargebayesian} and references therein. Next recall after $k$ iterations of the genGK bidiagonalization process, we have  
\begin{equation*}
    \bfQ^{1/2}\bfA\t\bfR^{-1}\bfA\bfQ^{1/2} \approx \bfQ^{1/2}\bfV_{k}\bfB_{k}\t\bfB_{k}\bfV_{k}\t\bfQ^{1/2}.
\end{equation*}
Let $\bfQ^{(1)}\bfR^{(1)} = \bfB_{k}$ and $\bfQ^{(2)}\bfR^{(2)} = \bfQ^{1/2}\bfV_{k}(\bfR^{(1)})\t$ be QR decompositions, then
\begin{align*}
    \bfQ^{1/2}\bfV_{k}\bfB_{k}\t\bfB_{k}\bfV_{k}\t\bfQ^{1/2} 
    &= \bfQ^{(2)}\bfR^{(2)}(\bfR^{(2)})\t(\bfQ^{(2)})\t.
\end{align*}
Taking the eigenvalue decomposition $\bfR^{(2)}(\bfR^{(2)})\t = \bfW_k \mathbf{\Theta}_{k}\bfW_{k}\t$ and assigning $\bfP_k = \bfQ^{(2)}\bfW_k$, we get the low-rank representation,
\begin{equation}
    \bfQ^{1/2}\bfA\t\bfR^{-1}\bfA\bfQ^{1/2} \approx \bfP_k\mathbf{\Theta}_{k}\bfP_k \quad \text{where} \quad \mathbf{\Theta}_{k} = \diag{\theta_{1},\ldots,\theta_{k}}.
\end{equation}

To compute the inverse square-root $\left( \delta\bfI + \lambda\bfQ^{1/2}\bfA\t\bfR^{-1}\bfA\bfQ^{1/2} \right)^{-1/2}$, following \cite{ghattasinfinitebayes}, we use the Woodbury identity to get
\begin{align*}
    \left( \delta\bfI + \lambda\bfQ^{1/2}\bfA\t\bfR^{-1}\bfA\bfQ^{1/2} \right)^{-1} 
    &= \delta^{-1}\left(\bfI + \bfP_{k}(\lambda/\delta)\mathbf{\Theta}_k \bfP_k \t \right)^{-1} \\
    &= \delta^{-1}\left(\bfI - \bfP_{k}\left(\bfI_k + (\lambda/\delta)^{-1}\mathbf{\Theta}_k^{-1} \right)^{-1}\bfP_{k}\t\right).
\end{align*}

 Then we can compute a square-root matrix as
 \begin{align}
    \delta^{-1/2} \left(\bfI + \bfP_{k}(\lambda/\delta)\mathbf{\Theta}_k \bfP_k \t \right)^{-1/2}
     &= \delta^{-1/2}\left(\bfI_{n} - \bfP_{k}\bfD_{k}\bfP_{k}\t\right)
 \end{align}
 where $\bfD_k \equiv \bfI_k - \left((\lambda/\delta)\mathbf{\Theta}_k + \bfI_k\right)^{-1/2}$.

\section{Prior norm for genGK proposal}
To evaluate $\norm[\bfQ^{-1}]{\bfx^{t-1}}^2$, we use an equivalent expression that avoids computations with $\bfQ^{-1}$. Let $\bfx_k=\bfQ\bfV_k\bfz_k$ be the genGK approximation used as the mean of the distribution from which $\bfx^{t-1}$ was drawn and $\bfD_k$ be the matrix used to form the corresponding $\widehat\bfGamma_{\rm cond}^{1/2}$.
Then, the norm can be written as
\begin{align*}
    & \norm[\bfQ^{-1}]{\bfx^{t-1}}^2 \\
    & = \left(\bfQ\bfV_{k}\bfz_k + \frac{\bfQ^{1/2}}{\sqrt{\delta}}\left(\bfI - \bfP_{k}
\bfD_{k}\bfP_{k}\t\right)\bfxi\right)\t \bfQ^{-1}\left(\bfQ\bfV_{k}\bfz_k + \frac{\bfQ^{1/2}}{\sqrt{\delta}}\left(\bfI - \bfP_{k}
\bfD_{k}\bfP_{k}\t\right)\bfxi\right) \\
&= \left(\bfQ^{1/2}\bfV_{k}\bfz_k + \delta^{-1/2}\left(\bfI - \bfP_{k}
\bfD_{k}\bfP_{k}\t\right)\bfxi\right)\t \left(\bfQ^{1/2}\bfV_{k}\bfz_k + \delta^{-1/2}\left(\bfI - \bfP_{k}
\bfD_{k}\bfP_{k}\t\right)\bfxi\right) \\
&= \left((\bfV_{k}\bfz_k)\t\bfQ^{1/2} + \delta^{-1/2}\bfxi\t\left(\bfI - \bfP_{k}
\bfD_{k}\bfP_{k}\t\right)\right) \left(\bfQ^{1/2}\bfV_{k}\bfz_k + \delta^{-1/2}\left(\bfI - \bfP_{k}
\bfD_{k}\bfP_{k}\t\right)\bfxi\right)
\end{align*}
so finally we have
\begin{equation}
    \norm[\bfQ^{-1}]{\bfx^{t-1}}^2 = \left(\bfV_k\bfz_k\right)\t(\bfx_k +2\widehat\bfGamma_{\rm cond}^{1/2}\bfxi) + \bfxi\t\delta^{-1}\left(\bfI- \bfP_k\left((\bfD_k)^2-2\bfD_k\right)\bfP_k\t\right)\bfxi.
\end{equation} 

\section{Derivation of acceptance ratio and prior norm for \Cref{sub:MH_genGK}}\label{sec:accept_norm_precond}

Let $\bfx_k$ be the genGK approximation used as the mean of the distribution from which $\bfx^{t-1}$ was drawn.
Consider the log of the full acceptance ratio,
\begin{align*}
&\log\jmc{\rho}_2(\bfx^{t-1},\bfx^{\star}) \\ &\quad = -\frac{1}{2}\Big( (\bfx^{\star} - \bfx_{\rm cond})\t \bfGamma_{\rm cond}^{-1}(\bfx^{\star} - \bfx_{\rm cond}) + (\bfx^{t-1} - \bfx_{k})\t \bfGamma_{\rm cond}^{-1}(\bfx^{t-1} - \bfx_{k})  \\
& \qquad - (\bfx^{t-1} - \bfx_{\rm cond})\t \bfGamma_{\rm cond}^{-1}(\bfx^{t-1} - \bfx_{\rm cond}) + (\bfx^{\star} - \bfx_{k})\t \bfGamma_{\rm cond}^{-1}(\bfx^{\star} - \bfx_{k}) \Big) \\
& \quad = \bfx_{\rm cond}\t\bfGamma_{\rm cond}^{-1}\bfx^{\star} + \bfx_{k}\t\bfGamma_{\rm cond}^{-1}\bfx^{t-1} - \bfx_{\rm cond}\t\bfGamma_{\rm cond}^{-1}\bfx^{t-1} - \bfx_{k}\t\bfGamma_{\rm cond}^{-1}\bfx^{\star},
\end{align*}
after expanding and canceling like terms. Using $\bfx_{\rm cond} = \lambda\bfGamma_{\rm cond}\bfA\t\bfR^{-1}\bfb$, $\bfGamma_{\rm cond}^{-1} = \lambda\bfA\t\bfR^{-1}\bfA + \delta\bfQ^{-1}$, and $\bfx_k = \bfQ\bfV_k\bfz_k$, the equation reduces to
\begin{align*}
\log\jmc{\rho}_2(\bfx^{t-1},\bfx^{\star}) &= \lambda(\bfA\t\bfR^{-1}\bfb)\t\bfx^\star + (\lambda\bfx_k\t\bfA\t\bfR^{-1}\bfA + \delta(\bfV_k\bfz_k)\t)\bfx^{t-1} \\
& \quad - \lambda(\bfA\t\bfR^{-1}\bfb)\t\bfx^{t-1} - (\lambda\bfx_k\t\bfA\t\bfR^{-1}\bfA + \delta(\bfV_k\bfz_k)\t)\bfx^{\star} \\
&=\left( \bfx^{\star} - \bfx^{t-1} \right)\t \left( \lambda\bfA\t\bfR^{-1}\left(\bfb - \bfA\bfx_k\right) - \delta\bfV_k\bfz_k\right)\\
&= \left( \bfx^{\star} - \bfx^{t-1} \right)\t \left( \lambda\left(\gamma_1\jmc{\zeta}_1\bfv_1 - \bfV_{k+1}\begin{bmatrix}
    \bfB_k\t \\ \jmc{\zeta}_{k+1}\bfe_{k+1}\t
\end{bmatrix} \bfB_k\bfz_k\right)- \delta\bfV_k\bfz_k\right).
\label{eq:accept_precond}
\end{align*}
Note that 
\begin{equation}
    \bfV_{k+1}\begin{bmatrix}
    \bfB_k\t \\ \jmc{\zeta}_{k+1}\bfe_{k+1}\t
\end{bmatrix} \bfB_k\bfz_k = \bfV_k\bfb_k\t\bfB_k\bfz_k + \jmc{\zeta}_{k+1}\gamma_{k+1}\left( \bfe_{k}\t\bfz_k \right)\bfv_{k+1}.
\end{equation}
Next, to avoid $\bfQ^{-1}$ when computing $\norm[\bfQ^{-1}]{\bfx^{t-1}}^2$, notice that the norm can be equivalently written as
\begin{align*}
\norm[\bfQ^{-1}]{\bfx^{t-1}}^2 
&= \left(\bfx^{t-1}\right)\t\bfQ^{-1}\left(\bfx_k + \bfS_F\bfxi\right) \\ 
& = \left(\bfx^{t-1}\right)\t\bfQ^{-1}\left(\bfQ\bfV_k\bfz_k + \lambda^{-1/2}\bfQ\bfG\t\left(\bfG\bfF\bfG\t\right)^{-1/2}\bfxi\right) 
\end{align*}
where $\bfF = \frac{\delta}{\lambda}\bfQ+\bfQ\bfA\t\bfR^{-1}\bfA\bfQ$. Distributing $\bfQ^{-1}$ results in
\begin{equation}
    \norm[\bfQ^{-1}]{\bfx^{t-1}}^2 = \left(\bfx^{t-1}\right)\t\left( \bfV_k\bfz_k + \lambda^{-1/2}\bfG\t\left(\bfG\bfF\bfG\t\right)^{-1/2}\bfxi\right)
\end{equation}
which reuses the previously computed $\bfV_k\bfz_k$ and $\bfG\t\left(\bfG\bfF\bfG\t\right)^{-1/2}\bfxi$.

\section{SVD and rSVD approximations of the conditional covariance matrix} \label{sec:svdapprox}
Assume the factorization $\bfQ=\bfL\t\bfL$ is accessible. Then, for fixed $\lambda$ and $\delta$,
\begin{equation*}
\bfGamma_{\rm cond} = \left(\delta\bfL^{-1}\bfL^{-\top} + \lambda\bfA\t\bfR^{-1}\bfA\right)^{-1} = \bfL\t\left(\delta\bfI + \lambda\bfL\bfA\t\bfR^{-1}\bfA\bfL\t \right)^{-1}\bfL.
\end{equation*}
To form a low-rank approximation of $\bfH = \bfL\bfA\t\bfR^{-1}\bfA\bfL\t$, we consider two approaches. 

One approach is to use the truncated SVD (tSVD) where a rank $k$ approximation is formed by computing $\bfR^{-1/2}\bfA\bfL\t\approx\widehat\bfU_k\bfSigma_k\widehat\bfV_k$ where $\widehat\bfU_k\in\bbR^{m\times k}$ and $\widehat\bfV_k\in\bbR^{n \times k}$ are orthonormal matrices containing the first $k$ left and right singular vectors of $\bfR^{-1/2}\bfA\bfL\t$ respectively and $\bfSigma_k\in\bbR^{k\times k}$ is a diagonal matrix with the first $k$ largest singular values. Then $\bfH\approx \widehat\bfV_k\bfSigma^2_k\widehat\bfV_k\t$. 

Another approach is to use the randomized SVD (rSVD) algorithm which begins by using a random matrix $\bfOmega\in\bbR^{n \times (k+p)}$ whose entries are realizations of i.i.d. standard Gaussian random variables to approximate the column space of $\bfH$, i.e., $\bfY=\bfH\bfOmega$. Here $k$ is the target rank and $p$ is an oversampling parameter usually taken to be a small integer ($p=5$ or $p=10$). Compute a thin-QR factorization, $\bfY=\widehat\bfQ\widehat\bfR$. Now an approximation of $\bfH$ is given by
$$
  \bfH \approx \widehat\bfQ\widehat\bfQ\t\bfH\widehat\bfQ\widehat\bfQ\t 
  = \widehat\bfQ \Tilde{\bfU}\Tilde{\bfSigma} \Tilde{\bfU}\t\widehat\bfQ\t 
  = \Tilde{\bfV}\Tilde{\bfSigma}\Tilde{\bfV}\t
$$
where $\widehat\bfQ\t\bfH\widehat\bfQ = \Tilde{\bfU}\Tilde{\bfSigma} \Tilde{\bfU}\t$ is an eigendecomposition and $\Tilde{\bfV} = \widehat\bfQ\Tilde{\bfU}$.

For both methods, the approximation of $\bfH$ can be used to form a low-rank representation of $\bfGamma_{\rm cond}$ following a similar approach to the one in \Cref{sub:sampling_prop}. Using the rSVD approximation as a proposal sampler was considered in \cite{brown2018low}.

The mean is then found using the approximate covariance matrix by computing
$$\bfx_{\rm cond} \approx \widehat{\bfGamma}_{\rm cond}(\lambda\bfA\t\bfR^{-1}\bfb) $$ which is used to draw a sample from the conditional distribution.

\end{appendices}


\bibliography{references}

@book{ghoshbook,
  title={An Introduction to {B}ayesian Analysis},
  author={Ghosh, J. K. and Delampady, M. and Samanta, T.},
  publisher={Springer},
  address={New York},
  year={2006}}

@article{saibaba2021randomized,
  title={Randomized approaches to accelerate {MCMC} algorithms for {B}ayesian inverse problems},
  author={Saibaba, Arvind K and Prasad, Pranjal and De Sturler, Eric and Miller, Eric and Kilmer, Misha E},
  journal={Journal of Computational Physics},
  volume={440},
  pages={110391},
  year={2021},
  publisher={Elsevier}
}

@incollection{dashti2015bayesian,
  title={The {B}ayesian approach to inverse problems},
  author={Dashti, Masoumeh and Stuart, Andrew M},
  booktitle={Handbook of uncertainty quantification},
  pages={1--118},
  year={2015},
  publisher="Springer International Publishing",
  address="Cham"
}

@book{ghanem2017handbook,
  title={Handbook of uncertainty quantification},
  author={Ghanem, Roger and Higdon, David and Owhadi, Houman and others},
  volume={6},
  year={2017},
  publisher={Springer},
  address ={New York}
}

@article{adler2019deep,
  title={Deep posterior sampling: {U}ncertainty quantification for large scale inverse problems},
  author={Adler, Jonas and {\"O}ktem, Ozan},
  year={2019}
}

@article{lan2022scaling,
  title={Scaling up {B}ayesian uncertainty quantification for inverse problems using deep neural networks},
  author={Lan, Shiwei and Li, Shuyi and Shahbaba, Babak},
  journal={SIAM/ASA Journal on Uncertainty Quantification},
  volume={10},
  number={4},
  pages={1684--1713},
  year={2022},
  publisher={SIAM}
}

@article{flath2011fast,
  title={Fast algorithms for {B}ayesian uncertainty quantification in large-scale linear inverse problems based on low-rank partial {H}essian approximations},
  author={Flath, H Pearl and Wilcox, Lucas C and Ak{\c{c}}elik, Volkan and Hill, Judith and van Bloemen Waanders, Bart and Ghattas, Omar},
  journal={SIAM Journal on Scientific Computing},
  volume={33},
  number={1},
  pages={407--432},
  year={2011},
  publisher={SIAM}
}

@article{martin2012stochastic,
  title={A stochastic {N}ewton {MCMC} method for large-scale statistical inverse problems with application to seismic inversion},
  author={Martin, James and Wilcox, Lucas C and Burstedde, Carsten and Ghattas, Omar},
  journal={SIAM Journal on Scientific Computing},
  volume={34},
  number={3},
  pages={A1460--A1487},
  year={2012},
  publisher={SIAM}
}

@book{nagydeblurring,
  title={Deblurring Images: Matrices, Spectra, and Filtering},
  author={Hansen, Per Christian and Nagy, James G. and O’Leary, Dianne P.},
  publisher = {SIAM},
  address = {Philadelphia},
  year = {2006}
}

@article{calvetti2019hierachical,
  title={Hierachical {B}ayesian models and sparsity: $\ell_2$-magic},
  author={Calvetti, Daniela and Somersalo, Erkki and Strang, Alexander},
  journal={Inverse Problems},
  volume={35},
  number={3},
  pages={035003},
  year={2019},
  publisher={IOP Publishing}
}

@article{lindbloom2025efficient,
  title={Efficient sparsity-promoting {MAP} estimation for {B}ayesian linear inverse problems},
  author={Lindbloom, Jonathan and Glaubitz, Jan and Gelb, Anne},
  journal={Inverse Problems},
  volume={41},
  number={2},
  pages={025001},
  year={2025},
  publisher={IOP Publishing}
}

@article{calvetti2015hierarchical,
  title={A hierarchical {K}rylov--{B}ayes iterative inverse solver for {MEG} with physiological preconditioning},
  author={Calvetti, Daniela and Pascarella, Annalisa and Pitolli, Francesca and Somersalo, Erkki and Vantaggi, Barbara},
  journal={Inverse Problems},
  volume={31},
  number={12},
  pages={125005},
  year={2015},
  publisher={IOP Publishing}
}

@article{agapiou2014analysis,
  title={Analysis of the {G}ibbs sampler for hierarchical inverse problems},
  author={Agapiou, Sergios and Bardsley, Johnathan M and Papaspiliopoulos, Omiros and Stuart, Andrew M},
  journal={SIAM/ASA Journal on Uncertainty Quantification},
  volume={2},
  number={1},
  pages={511--544},
  year={2014},
  publisher={SIAM}
}

@article{ascolani2024scalability,
  title={Scalability of {M}etropolis-within-{G}ibbs schemes for high-dimensional {B}ayesian models},
  author={Ascolani, Filippo and Roberts, Gareth O and Zanella, Giacomo},
  journal={arXiv preprint arXiv:2403.09416},
  year={2024}
}

@article{glaubitz2025efficient,
  title={Efficient sampling for sparse {B}ayesian learning using hierarchical prior normalization},
  author={Glaubitz, Jan and Marzouk, Youssef},
  journal={arXiv preprint arXiv:2505.23753},
  year={2025}
}

@article{lindbloom2025priorconditioned,
  title={Priorconditioned sparsity-promoting projection methods for deterministic and {B}ayesian linear inverse problems},
  author={Lindbloom, Jonathan and Pasha, Mirjeta and Glaubitz, Jan and Marzouk, Youssef},
  journal={arXiv preprint arXiv:2505.01827},
  year={2025}
}

@article{calvetti2024computationally,
  title={Computationally efficient sampling methods for sparsity promoting hierarchical {B}ayesian models},
  author={Calvetti, Daniela and Somersalo, Erkki},
  journal={SIAM/ASA Journal on Uncertainty Quantification},
  volume={12},
  number={2},
  pages={524--548},
  year={2024},
  publisher={SIAM}
}

@article{fox2016fast,
  title={Fast sampling in a linear-{G}aussian inverse problem},
  author={Fox, Colin and Norton, Richard A},
  journal={SIAM/ASA Journal on Uncertainty Quantification},
  volume={4},
  number={1},
  pages={1191--1218},
  year={2016},
  publisher={SIAM}
}

@article{calvetti2025subspace,
  title={Subspace Splitting Fast Sampling from {G}aussian Posterior Distributions of Linear Inverse Problems},
  author={Calvetti, Daniela and Somersalo, Erkki},
  journal={arXiv preprint arXiv:2502.05703},
  year={2025}
}

@article{sanz2025hierarchical,
  title={Hierarchical {B}ayesian inverse problems: A high-dimensional statistics viewpoint},
  author={Sanz-Alonso, Daniel and Waniorek, Nathan},
  journal={SIAM Review},
  volume={67},
  number={3},
  pages={543--575},
  year={2025},
  publisher={SIAM}
}

@article{agrawal2022variational,
  title={A variational inference approach to inverse problems with gamma hyperpriors},
  author={Agrawal, Shiv and Kim, Hwanwoo and Sanz-Alonso, Daniel and Strang, Alexander},
  journal={SIAM/ASA Journal on Uncertainty Quantification},
  volume={10},
  number={4},
  pages={1533--1559},
  year={2022},
  publisher={SIAM}
}

@article{brown2018low,
  title={Low-rank independence samplers in hierarchical {B}ayesian inverse problems},
  author={Brown, D Andrew and Saibaba, Arvind and Vall{\'e}lian, Sarah},
  journal={SIAM/ASA Journal on Uncertainty Quantification},
  volume={6},
  number={3},
  pages={1076--1100},
  year={2018},
  publisher={SIAM}
}

@Misc{amsmath,
  author =	 {{American Mathematical Society}},
  title =	 {User's Guide for the \texttt{amsmath} Package
                  (Version 2.0)},
  url =		 {ftp://ftp.ams.org/pub/tex/doc/amsmath/amsldoc.pdf},
  urldate =	 {2015-07-30},
  year =	 2002}

@book{gelman2013bayesian,
  title={Bayesian Data Analysis, Third Edition},
  author={Gelman, A. and Carlin, J.B. and Stern, H.S. and Dunson, D.B. and Vehtari, A. and Rubin, D.B.},
  isbn={9781439840955},
  lccn={2013039507},
  series={Chapman \& Hall/CRC Texts in Statistical Science},
  url={https://books.google.com/books?id=ZXL6AQAAQBAJ},
  year={2013},
  publisher={Taylor \& Francis},
  address = {Florida}
}

@book{calvetti2007introduction,
  title={An Introduction to Bayesian Scientific Computing: Ten Lectures on Subjective Computing},
  author={Calvetti, Daniela and Somersalo, Erkki},
  volume={2},
  year={2007},
  publisher={Springer Science \& Business Media},
   address = {New York}
}

@book{calvetti2023bayesianbook,
  title={Bayesian Scientific Computing},
  author={Calvetti, Daniela and Somersalo, Erkki},
  year={2023},
  publisher={Springer Nature},
   address = {Switzerland},
    doi={10.1007/978-3-031-23824-6}
}

@book{bardsley2018computational,
  title={Computational Uncertainty Quantification for Inverse Problems},
  author={Bardsley, Johnathan M},
   volume={19},
  year={2018},
   address = {Philadelphia},
  publisher={SIAM}
}

@article{donoho1995noising,
  title={De-noising by soft-thresholding},
  author={Donoho, David L},
  journal={IEEE Transactions on Information Theory},
  volume={41},
  number={3},
  pages={613--627},
  year={1995},
  publisher={IEEE}
}

@article{BuiThanh2014,
   doi = {10.1137/120894877},
   url = {https://doi.org/10.1137/120894877},
  year = {2014},
   month = jan,
  publisher = {Society for Industrial {\&} Applied Mathematics ({SIAM})},
  volume = {2},
  number = {1},
  pages = {203--222},
  author = {Tan Bui-Thanh and Omar Ghattas},
  title = {An analysis of infinite dimensional {B}ayesian inverse shape acoustic scattering and its numerical approximation},
  journal = {{SIAM}/{ASA} Journal on Uncertainty Quantification}
}

@article{gazzola2019ir,
  title={{IR Tools}: a {MATLAB} package of iterative regularization methods and large-scale test problems},
  author={Gazzola, Silvia and Hansen, Per Christian and Nagy, James G},
  journal={Numerical Algorithms},
  volume={81},
  number={3},
  pages={773--811},
  year={2019},
  publisher={Springer}
}

@article{ghattasinfinitebayes,
	author = {Bui-Thanh, Tan and Ghattas, Omar and Martin, James and Stadler, Georg},
	journal = {SIAM Journal on Scientific Computing},
	number = {6},
	pages = {A2494-A2523},
	title = {A Computational Framework for Infinite-Dimensional {B}ayesian Inverse Problems Part {I}: The Linearized Case, with Application to Global Seismic Inversion},
	volume = {35},
	year = {2013}}

@article{UQlargebayesian,
author = {Saibaba, Arvind and Chung, Julianne and Petroske, Katrina},
year = {2020},
month = {08},
pages = {},
title = {Efficient {K}rylov subspace methods for uncertainty quantification in large {B}ayesian linear inverse problems},
volume = {27},
journal = {Numerical Linear Algebra with Applications},
doi = {10.1002/nla.2325}
}

@article{ChungSaibabaHybrid2017,
	author = {Chung, Julianne and Saibaba, Arvind K.},
	journal = {SIAM Journal on Scientific Computing},
	number = {5},
	pages = {S24-S46},
	title = {Generalized Hybrid Iterative Methods for Large-Scale {B}ayesian Inverse Problems},
	volume = {39},
	year = {2017}}

@article{saibaba2019efficient,
  title={Efficient marginalization-based {MCMC} methods for hierarchical {B}ayesian inverse problems},
  author={Saibaba, Arvind K and Bardsley, Johnathan and Brown, D Andrew and Alexanderian, Alen},
  journal={SIAM/ASA Journal on Uncertainty Quantification},
  volume={7},
  number={3},
  pages={1105--1131},
  year={2019},
  publisher={SIAM}
}

@article{chung2024computational,
  title={Computational methods for large-scale inverse problems: a survey on hybrid projection methods},
  author={Chung, Julianne and Gazzola, Silvia},
  journal={SIAM Review},
  volume={66},
  number={2},
  pages={205--284},
  year={2024},
  publisher={SIAM}
}

@article{arioli2013generalized,
  title={Generalized {G}olub--{K}ahan bidiagonalization and stopping criteria},
  author={Arioli, Mario},
  journal={SIAM Journal on Matrix Analysis and Applications},
  volume={34},
  number={2},
  pages={571--592},
  year={2013},
  publisher={SIAM}
}

@article{hall2024efficient,
  title={Efficient iterative methods for hyperparameter estimation in large-scale linear inverse problems},
  author={Hall-Hooper, Khalil A and Saibaba, Arvind K and Chung, Julianne and Miller, Scot M},
  journal={Advances in Computational Mathematics},
  volume={50},
  number={6},
  pages={1--33},
  year={2024},
  publisher={Springer}
}

@article{pasha2023computational,
  title={A computational framework for edge-preserving regularization in dynamic inverse problems},
  author={Pasha, Mirjeta and Saibaba, Arvind K and Gazzola, Silvia and Espa{\~n}ol, Malena I and {de Sturler}, Eric},
  journal={Electronic Transactions on Numerical Analysis},
  volume={58},
  pages={486--516},
  year={2023},
  publisher={Kent State University}
}

@article{chung2018efficient,
  title={Efficient generalized {G}olub--{K}ahan based methods for dynamic inverse problems},
  author={Chung, Julianne and Saibaba, Arvind K and Brown, Matthew and Westman, Erik},
  journal={Inverse Problems},
  volume={34},
  number={2},
  pages={024005},
  year={2018},
  publisher={IOP Publishing}
}

@article{cho2022computationally,
  title={Computationally efficient methods for large-scale atmospheric inverse modeling},
  author={Cho, Taewon and Chung, Julianne and Miller, Scot M and Saibaba, Arvind K},
  journal={Geoscientific Model Development},
  volume={15},
  number={14},
  pages={5547--5565},
  year={2022},
  publisher={Copernicus GmbH}
}

@Book{MCstatmethods,
  author =	 {Christian P. Robert and George Casella},
  title =	 {Monte Carlo Statistical Methods},
  publisher =	 {Springer},
address = {New York},
  year =	 {2004},
  edition =	 {2nd}}

@Article{DataReductInvModel,
AUTHOR = {Liu, X. and Weinbren, A. L. and Chang, H. and Tadi\'c, J. M. and Mountain, M. E. and Trudeau, M. E. and Andrews, A. E. and Chen, Z. and Miller, S. M.},
TITLE = {Data reduction for inverse modeling: an adaptive approach v1.0},
JOURNAL = {Geoscientific Model Development},
VOLUME = {14},
YEAR = {2021},
NUMBER = {7},
PAGES = {4683--4696},
URL = {https://gmd.copernicus.org/articles/14/4683/2021/},
DOI = {10.5194/gmd-14-4683-2021}
}

@Article{GeostatInvModel,
AUTHOR = {Miller, S. M. and Saibaba, A. K. and Trudeau, M. E. and Mountain, M. E. and Andrews, A. E.},
TITLE = {Geostatistical inverse modeling with very large datasets: an example from the Orbiting Carbon Observatory 2 ({OCO}-2) satellite},
JOURNAL = {Geoscientific Model Development},
VOLUME = {13},
YEAR = {2020},
NUMBER = {3},
PAGES = {1771--1785},
URL = {https://gmd.copernicus.org/articles/13/1771/2020/},
DOI = {10.5194/gmd-13-1771-2020}
}

@article{LagrangTransport,
author = {Lin, J. C. and Gerbig, C. and Wofsy, S. C. and Andrews, A. E. and Daube, B. C. and Davis, K. J. and Grainger, C. A.},
title = {A near-field tool for simulating the upstream influence of atmospheric observations: The Stochastic Time-Inverted Lagrangian Transport ({STILT}) model},
journal = {Journal of Geophysical Research: Atmospheres},
volume = {108},
number = {D16},
pages = {},
doi = {https://doi.org/10.1029/2002JD003161},
url = {https://agupubs.onlinelibrary.wiley.com/doi/abs/10.1029/2002JD003161},
eprint = {https://agupubs.onlinelibrary.wiley.com/doi/pdf/10.1029/2002JD003161},
year = {2003}
}

@article{LagrangTransportModel,
    author = {Nehrkorn, T. and Eluszkiewicz, J. and Wofsy, Steven C. and Lin, John C. and Gerbig, C. and Longo, M. and Freitas, S.},
    title = {Coupled weather research and forecasting–stochastic time-inverted lagrangian transport ({WRF–STILT}) model},
    journal = {Meteorology and Atmospheric Physics},
    year = {2010},
    volume = {107},
pages = {51-64},
url = {https://doi.org/10.1007/s00703-010-0068-x},
doi = {10.1007/s00703-010-0068-x}    
}

@ARTICLE{geman1984Gibbs,
  author={Geman, Stuart and Geman, Donald},
  journal={IEEE Transactions on Pattern Analysis and Machine Intelligence}, 
  title={Stochastic Relaxation, Gibbs Distributions, and the Bayesian Restoration of Images}, 
  year={1984},
  volume={PAMI-6},
  number={6},
  pages={721-741},
  doi={10.1109/TPAMI.1984.4767596}}

@article{Metropolis_Rosenbluth_Rosenbluth_Teller_Teller_1953, title={Equation of state calculations by Fast Computing Machines}, volume={21}, DOI={10.1063/1.1699114}, number={6}, journal={The Journal of Chemical Physics}, author={Metropolis, Nicholas and Rosenbluth, Arianna W. and Rosenbluth, Marshall N. and Teller, Augusta H. and Teller, Edward}, year={1953}, month={Jun}, pages={1087–1092}}

@article{Hastings_1970, title={Monte {C}arlo sampling methods using Markov chains and their applications}, volume={57}, DOI={10.2307/2334940}, number={1}, journal={Biometrika}, author={Hastings, W. K.}, year={1970}, month={Apr}, pages={97}}

@article{Vats_Knudson_2021, title={Revisiting the {G}elman–{R}ubin diagnostic}, volume={36}, DOI={10.1214/20-sts812}, number={4}, journal={Statistical Science}, author={Vats, Dootika and Knudson, Christina}, year={2021}, month=nov }

@article{gelman_rubin_diag_1992,
author = {Andrew Gelman and Donald B. Rubin},
title = {{Inference from Iterative Simulation Using Multiple Sequences}},
volume = {7},
journal = {Statistical Science},
number = {4},
publisher = {Institute of Mathematical Statistics},
pages = {457 -- 472},
year = {1992},
doi = {10.1214/ss/1177011136},
URL = {https://doi.org/10.1214/ss/1177011136}
}

@article{petra2014computational,
  title={A computational framework for infinite-dimensional Bayesian inverse problems, {P}art {II}: Stochastic {N}ewton {MCMC} with application to ice sheet flow inverse problems},
  author={Petra, Noemi and Martin, James and Stadler, Georg and Ghattas, Omar},
  journal={SIAM Journal on Scientific Computing},
  volume={36},
  number={4},
  pages={A1525--A1555},
  year={2014},
  publisher={SIAM}
}

@Inbook{Sokal1997,
author="Sokal, A.",
editor="DeWitt-Morette, Cecile
and Cartier, Pierre
and Folacci, Antoine",
title="Monte Carlo Methods in Statistical Mechanics: Foundations and New Algorithms",
bookTitle="Functional Integration: Basics and Applications",
year="1997",
publisher="Springer US",
address="Boston, MA",
pages="131--192",
isbn="978-1-4899-0319-8",
doi="10.1007/978-1-4899-0319-8_6",
url="https://doi.org/10.1007/978-1-4899-0319-8_6"
}

@article{Calvetti2020sparsehyperpriors,
doi = {10.1088/1361-6420/ab4d92},
url = {https://doi.org/10.1088/1361-6420/ab4d92},
year = {2020},
month = {jan},
publisher = {IOP Publishing},
volume = {36},
number = {2},
pages = {025010},
author = {Calvetti, Daniela and Pragliola, Monica and Somersalo, Erkki and Strang, Alexander},
title = {Sparse reconstructions from few noisy data: analysis of hierarchical Bayesian models with generalized gamma hyperpriors},
journal = {Inverse Problems}
}

\end{document}